\documentclass{article}[10pt]
\usepackage{graphicx} 
\usepackage{amsmath, amssymb ,amsthm, amsfonts, amsgen, color, mathtools}
\usepackage{graphicx, hyperref}
\usepackage{bbm}
\usepackage{tikz}
\usepackage[colorinlistoftodos]{todonotes}
\usepackage[bottom]{footmisc}
\usepackage{mathabx}
\numberwithin{equation}{section}
\definecolor{vg}{rgb}{0.0, 0.26, 0.15}

\newcommand{\Nb}{{\mathbb{N}}}
\newcommand{\R}{{\mathbb{R}}}

\newcommand{\A}{{\mathcal{A}}}
\newcommand{\beq}{\begin{equation}}
\newcommand{\eeq}{\end{equation}}

\newcommand{\weaklystar}{\overset{\ast}{\rightharpoonup}}

\newcommand{\Fcal}{{\mathcal{F}}}

\newcommand{\Mcal}{{\mathcal{M}}}
\newcommand{\Lcal}{{\mathcal{L}}}

\newcommand{\eps}{{\epsilon}}
\newcommand{\diff}{{\rm diff}}

\renewcommand{\eps}{\varepsilon}

\def\e{{\varepsilon}}
\def\O{{\Omega}}

\def\Q{\mathcal Q}
\def\A{\mathcal A}
\def\B{\mathcal B}
\def\L{\mathcal L}
\def\Lin{\mbox{Lin}}
\def\symA{\mathbb{A}}
\def\symB{\mathbb{B}}
\def\Leb{\mathcal{L}}
\def\Q{\mathcal{Q}}

\def\XXint#1#2#3{{\setbox 0=\hbox{$#1{#2#3}{\int}$}
\vcenter{\hbox{$#2#3$}}\kern-.5\wd0}}

\definecolor{vg}{rgb}{0.0, 0.26, 0.15}
\definecolor{coolblack}{rgb}{0.0, 0.18, 0.39}

\newcommand{\res}{\mathop{\hbox{\vrule height 7pt width .5pt depth 0pt
\vrule height .5pt width 6pt depth 0pt}}\nolimits}

\newcommand{\weakcon}{\overset{\ast}{\rightharpoonup}}

\newtheorem{Theorem}{Theorem}[section]
\newtheorem{Lemma}[Theorem]{Lemma}
\newtheorem{Proposition}[Theorem]{Proposition}
\newtheorem{Corollary}[Theorem]{Corollary}

\newtheorem{Remark}[Theorem]{Remark}
\newtheorem{Definition}[Theorem]{Definition}

\def\f1min{f_1^{\rm min}}

\usepackage{color}
\usepackage{natbib}
\setcitestyle{numbers,square,semicolon}

\title{Relaxation of variational problems in the space of functions with bounded $\B$-variation: interaction with measures and lack of concentration phenomena}
\author{Lorenza D'Elia, Elvira Zappale}

\AtEndDocument{\bigskip{\footnotesize%
 Lorenza D'Elia\\ \textsc{Institute of Analysis and Scientific Computing, TU Wien, Wiedner Hauptstraße 8-10, 1040 Vienna, Austria} \par  
  \textit{E-mail address}: \texttt{lorenza.delia@tuwien.ac.at} \par
  \textit{ORCID}: \href{https://orcid.org/0000-0002-2809-5553}{0000-0002-2809-5553}
  \par
  \addvspace{\medskipamount}
  Elvira Zappale\\\textsc{Dipartimento di Scienze di Base e Applicate per l’Ingegneria, Sapienza - Universit\`{a} di Roma,
via Antonio Scarpa, 16, 00161 Roma, Italy,  and
CIMA, Universidade de \'Evora, Portugal} \par
  \textit{E-mail address}: \texttt{elvira.zappale@uniroma1.it}
\par
\textit{ORCID}:\href{https://orcid.org/0000-0001-7419-300X}{0000-0001-7419-300X}
}
}

\date{}

\begin{document}

\maketitle

\begin{abstract}  

We prove an integral representation result for variational functionals in the space $BV^{\B}$ of functions with bounded $\B$-variation where $\B$ denotes a $k$-th order, $\mathbb C$-elliptic, linear homogeneous differential operator. This result has been used as a key tool to get an explicit representation of relaxed energies with linear growth which lead to limiting generic measures. According to the space dimension and the order of the operator, concentration phenomena appear and an explicit interaction is featured.
These results are complemented also with Sobolev-type counterparts.  As a further application, a lower semicontinuity result in the space of fields with $p(\cdot)$-bounded $\B$- variation has also been obtained.

\medskip

{\scriptsize
}
\end{abstract}

\section{Introduction}
A central question in Calculus of Variations is the investigation of the appearance of oscillation and concentration effects. The latter mainly pops up in problems with no compactness and convexity-type properties, which are in general due to the non-reflexivity of the underlying function spaces. This, therefore, calls for a suitable extension procedure, the so-called {\it relaxation} (see \cite{DM}), to a space where there is compactness of bounded sequences. In the present manuscript, we carry out the relaxation of the following integral functionals  
\begin{align}\label{G}G(u,v)\coloneqq\int_\Omega \mathfrak{g}(x,u(x),v(x),\B u(x))dx=\int_\Omega g_1(u(x),v(x))dx + \int_\Omega g_2(x,\mathcal B u(x)) dx,
\end{align} where $\Omega$ is an open and bounded subset of $\R^N$ and  $\mathcal B$ represents a $k$-th order homogeneous $\mathbb C$-elliptic, linear partial differential operator, with constant coefficients, on $\R^N$ from $V$ to
$W$ ($V$ and $W$ being Euclidean vector spaces on $\mathbb R$ or $\mathbb C$). Throughout the paper, $\B$ takes the form
\begin{align}\label{Budef}
\mathcal B u=\sum_{|\beta|=k}B_\beta \partial^\beta u, \;\;\;
\mbox{for } u\in \mathcal{D'}(\R^N; V),
\end{align}
with $B_\beta \in {\rm Lin}(V, W)$. Here $\beta \in \mathbb N_0^N$ is a multi-index with  order $|\beta| = \beta_1 +\dots +\beta_N$ and $\partial _\beta$ represents the distributional derivative $\partial_1^{\beta_1}\dots
\partial_N^{\beta_N}$
(see subsection \ref{PDO} for precise definitions and properties).
In \eqref{G},  the function $u$ lies in the function spaces of Sobolev-type $W^{\B, 1}(\Omega)$ associated to the operator $\mathcal B$ (see Subsection \ref{FS} for details) and $v \in L^1(\Omega;\mathbb R^l)$ and, in the sequel, when no ambiguity could be created, the $x$ dependence on the densities and functions will be dropped in the integrals defining $G$.

 Concretely, we consider the case where the densities $g_1$ and $g_2$ have a form similar to the one in  \cite{KKZ23}, i.e., $g_1(u,v)=f_1(u)f_2(v)$ and $g_2(x,\xi)= h(\xi)$. Therefore, the functional $G$ in \eqref{G} turns into 
\begin{equation}
\label{GG1}
G_1(u,v)=\int_\Omega\left( f_1(u(x))f_2(v(x)) + h(\B u(x))\right)dx.
\end{equation}     
where $f_1: V\to \mathbb R$, $f_2 :\mathbb R^l\to \mathbb R$, and $h: W \to \mathbb R$ are continuous functions such that 
		
		\begin{itemize}
			\item[$(H_1)$] there exist $C_2>C_1>0$ such that  for every $a \in \mathbb R^d$
			\begin{equation}\nonumber C_1\leq f_1(a) \leq C_2;\end{equation}
			\item[$(H_2q)$] there exists $K>0$, such that  for every $b \in \mathbb R^l$
			\begin{equation}\nonumber
			K^{-1}|b|^q\leq f_2(b)\leq K(1+ |b|^q)\;\;\;\hbox{for } 1\leq q\in  \mathbb R^+;
			\end{equation} 
			\item[$(H_3p)$] there exists $\kappa >0$ such that for every $\xi \in \mathbb W$
			\begin{equation}\nonumber
			\kappa^{-1}|\xi|^p \leq h(\xi)\leq \kappa(1 +|\xi|^p) \;\;\;\hbox{for } 1\leq p\in  \mathbb R^+;
			\end{equation}
            \item[$(H_4)$] the strong recession function of $h$ exists in  $\mathcal{R}(\B)$, where the latter set represents the span of the image cone associated with $\B$ (see \eqref{imagecone} below and Subsection \ref{subsecint} for recession functions and related properties and symbols).
		\end{itemize}
The case when $\mathcal B$ is the gradient $\nabla$ has been extensively studied by many authors, see, e.g.,  \cite{HKW, KKK, KK, KKZ23} and the bibliography therein.

Different integral representations for the relaxed functional $\mathcal G_1$ will emerge according to assumptions on the densities $g_1, g_2$. More specifically, the core of our investigation is the analysis of the linear-growth models, i.e., models where Assumptions $(H_{2q})$ and $(H_{3p})$  with  $p=q=1$ hold. This alongside the particular form of $g_1(u,v)= f_1(u)f_2(v)$ might feature  a "strong" interaction term between $u$ and $v$.  As in the case $\B=\nabla$ (see \cite{KKZ23}), it is possible to highlight oscillation effects that can be caused by the non-convexity of the energy density in $\B u$ and $v$ with 
possible concentrations of minimizing sequences related to the lack of weak compactness of $W^{B,1}(\O)\times L^1(\O;\R^l)$. 
Therefore, an extension procedure for the functional $G_1$ to a space with compactness properties is needed. The compactness property is retrieved in the  $BV$-type spaces $BV^\B(\O)$ associated to the operator $\B$ (see Subsection \ref{FS} for details) for the variable $u$ and the Radon measures $\mathcal{M}(\bar\Omega;\R^l)$ on $\bar\Omega$ for the variable $v$. In the present manuscript, we are going to provide an explicit integral representation of the  
{\it relaxed} functional $\overline{\mathcal G}_1$ of $G_1$ given by
 \begin{align}\label{G1relax}
		\begin{aligned}
		\overline{\mathcal G}_1(u,v):=
		\inf\left\{\,\liminf_{j\to +\infty}
		G_1(u_j,v_j)
		\,\left|\,
		\begin{array}{l}
		(u_j,v_j)\in W^{\B,1}(\Omega)\times L^1(\Omega;\mathbb R^l),\\
		(u_j,v_j)\overset{\ast}{\rightharpoonup} (u,v)\hbox{ in } BV^\B(\Omega)\times \mathcal M(\bar\Omega;\mathbb R^l)
		\end{array}
		\right.\right\}.
		\end{aligned}
		\end{align}
 The main challenge is represented by the fact that  the integral representation of $\overline{\mathcal G}_1$ in our framework has an explicit dependence on $\B$, and in particular on its order $k$. In fact, for $k$ sufficiently large, compared with the dimension space $N$, no concentration effect and interaction between $u$ and $v$ appear.
More precisely, for $k<N$, we prove the following result.
  \begin{Theorem}\label{thm:relk<N,N>=2}
   Let $N\geq 2$ and $\Omega$ be a bounded Lipschitz domain in $\R^N$.  Let $\B$ be a linear partial differential operator, $\mathbb{C}$-elliptic of order $k$, with $k< N$, as in \eqref{Budef}. Let $f_1, f_2, h$ be functions satisfying assumptions ${\rm (H_1)}, {\rm (H_2q)}, {\rm (H_3p)}$ and ${\rm (H_4)}$, for $p=q=1$. Assume also that $\overline {\mathcal G}_1\in \mathcal G \mathcal M(\Omega)$ (see Definition \eqref{GM}). Then, we have that 
        \begin{align}\nonumber
            \overline{\mathcal G}_1(u,v)=&
				\int_\Omega d\mathcal Q_\B h(\mathcal B u)(x)  
			+\int_\Omega  g\left(u(x), \frac{d v^a}{d {\mathcal L}^N}(x)\right)\,dx
				+\int_{\bar\Omega} \f1min \,
				(f_2^{\ast\ast})^\infty\left(\frac{d v^s }{d|v^s|}(x)\right)\,d|v^s|(x),\nonumber
        \end{align}
        for any $u\in BV^{\B}(\Omega)$ and $v\in\Mcal(\Omega; \R^l)$. Here, with an abuse of notation, we have identified 
        ${\B u}$ with the Radon–Nykod\'{y}m derivative of $\B u$ with respect to $\L^N$, $\mathcal Q_\B h$ denotes the $\B$-quasiconvexification of $h$,  $d\mathcal Q_\B h$ denotes   $\displaystyle{\int_{\Omega} \Q_{\B}h(\B u(x))dx + \int_{\Omega} (\Q_{\B}h)^\infty\left(\frac{d\B u}{d |\B^s u|}(x)\right)d |\B^s u|(x)}$, $f_2^{\ast\ast}$ is the bipolar function of $f_2$, and $g$ is given by
              \begin{equation}
              \label{defg}
                  g(a,b) := \min\left\{ f_1(a)f_2^{\ast\ast}(b_1) + f_1^{\min}(f_2^{\ast\ast})^\infty(b_2) : b_1, b_2\in\R^l, b_1 + b_2=b  \right\},
              \end{equation}
        and $f_1^{\min} := \inf_{a\in\R^N} f_1(a) \geq C_1.$ The superscript $\infty$ denotes the recession function 
        (see Subsection \ref{subsecABquasiconvex} and \eqref{exiQBfinfty} for definitions and properties).
   \end{Theorem}

The importance of detecting the phenomena enlighten by Theorem \ref{thm:relk<N,N>=2} relies on the fact that it might appear in the modeling of impulse control, where interactions of oscillations, concentrations and discontinuities form naturally in minimization procedures associated to functionals as in \eqref{GG1}. This is mainly due to the fact that discontinuous functions do not belong to the predual space of Radon measures. Furthermore, the concentration  effect described in Theorem \ref{thm:relk<N,N>=2} is intimately connected with 
nonreflexivity of underlying spaces $L^1$ and $W^{\B,1}$ (or simply $W^{k,1}$).

Theorem \ref{thm:relk<N,N>=2}  represents the $\B$-variant of \cite[Theorem 1.1]{KKZ23}, and the obtained integral representation is completely analogous. However, contrary to the case when $\B=\nabla$, we impose a restriction on the order of the operator $\B$. This is also the cause of the technical assumption on  $\overline{\mathcal G}_1$, trivially satisfied by all first-order operators and when $BV^\B= BH$. Such a technical assumption has its root in the global method for relaxation, which turns out to be a fundamental tool to carry out our analysis. Diverse and direct approaches might be used, playing around with the explicit structure of the operators $\B$. For example, when $N=1$, the global method is unnecessary since a blow-up approach (introduced in \cite{FM92}) can be adopted. 
However, we believe that the global method approach is more general and schematic and this is reason why a global method for relaxation for $\mathbb C$-elliptic linear partial differential operators $\B$ of any order is the subject of an ongoing research.

Theorem \ref{thm:relk<N,N>=2} is complemented by the following result. 

 \begin{Theorem}
      \label{thm:relKgeq N}
      Let $\Omega$, $f_1, f_2, h$ and $\B$ be as in Theorem \ref{thm:relk<N,N>=2} and assume that $k \geq N$ ($k>N$ if $N\leq 2$), and $p=q=1$. 
      Let  $G_1$ and $\overline{\mathcal G_1}(u,v)$ as in \eqref{GG1} and \eqref{G1relax}, respectively. Then  
                 \begin{align}\nonumber
                     \overline{\mathcal G_1}(u,v) &= \int_{\overline{\Omega}} f_1(u(x))df_2^{\ast\ast}(v)(x)+\int_\Omega d\mathcal Q_{\B}h (\B u)(x). 
                 \end{align}
  \end{Theorem}

Both Theorems \ref{thm:relk<N,N>=2} and \ref{thm:relKgeq N} rely on an intermediate integral representation result which plays a crucial role in our analysis. It is key to  understanding the relaxation of functional $G_1$ as in \eqref{GG1}, when $f_1=f_2=0$, i.e. providing an integral representation for the {\it relaxed} functional of 
$W^{\B,1}(\Omega)  \ni u \mapsto G_\B(u):=\int_\Omega g_2(x, \B u(x)) dx$. We prove in Section \ref{secRel} the following result.

\begin{Theorem}
\label{thm:RelaxationB} Let $\B$ be a $k$-th order, $\mathbb C$-elliptic, linear partial differential operator as in \eqref{Budef}.
    Let $g_2:\Omega\times W\to [0,\infty)$ be a continuous integrand such that
    \begin{itemize}
       \item[$(H_{1 \B})$] $g_2$ is Lipschitz continuous in the second variable uniformly in $x$;
       \item[$(H_{2 \B})$] $g_2$ has linear growth at infinity in the second argument, i.e., there exists $M>0$ such that
           \begin{equation}
               \notag
               |g_2(x, \xi)|\leq M(1+|\xi|), \qquad\mbox{for all } x\in\Omega, \xi\in W;
           \end{equation}
       \item[$(H_{3 \B})$] there exists a modulus of continuity $\omega:[0,\infty)\to [0,\infty)$ such that
              \begin{equation}
                  \notag
                  |g_2(x, \xi) - g_2(y, \xi)|\leq \omega(|x-y|)(1+ |\xi|),\qquad\mbox{for all } x,y\in\Omega,\;\; \xi\in W;
              \end{equation}
              \item[$(H_{4 \B})$] the strong recession function $(g_2)^{\infty}(x, \xi)$ exists for any $(x, \xi) \in \Omega\times \mathcal R({\B})$.
   \end{itemize}
 Let  \begin{equation}
                      \label{GB}
                      G_{\B}(u):= \int_{\Omega} g_2(x, \B u(x))dx,
                  \end{equation}
                  and \begin{align}
        \label{GBcal}
        \overline{\mathcal{G}}_{\B}(u):= \inf \biggl\{ \liminf_{j\to\infty} G_{\B}(u_j): W^{\mathcal B, 1}(\Omega)     \ni u_j\to u 
        \hbox{ in } W^{k-1,1}(\Omega; V),
    \B u_j\mathcal{L}^n \weakcon \B u \hspace{0.3cm} \mbox{in } \Mcal(\Omega;W)    \biggr\}.
    \end{align}
   Then, $\overline{\mathcal{G}}_{\B}$ is characterized by 
    \begin{equation}
              \notag
    \overline{\mathcal{G}}_{\B}[ u] = \int_\Omega Q_{\B}g_2\left( x, {\B u}(x) \right)dx + \int_\Omega (Q_{\B}g_2)^\#\left( x, {d\B^s u\over d|\B^s u|}(x)\right)d|\B^s u|(x),
          \end{equation}
    where $Q_{\B}g_2$ is the $\mathcal B$- quasi-convex envelope of $g_2(x,\cdot)$, for every $x \in \Omega$ and $(Q_{\B}g_2 )^\#$ is the generalized upper recession function (see \eqref{fdown} below for the definition) of $Q_{\B}g_2$.
\end{Theorem}

This result, in the specific case when $\B=\nabla$, has been proved by the pioneering papers \cite{AD} and \cite{FMARMA}. Related results under more general assumptions are contained in \cite{BFMglob}, \cite{KR10} and \cite{RS}.
For $\B=\mathcal E$, we refer to \cite{BFT} as well as \cite{CFVG20}, and to \cite{ADC}, \cite{FHP} and \cite{H} in the case $\B=D^2$.
Results dealing with general first-order operators $\B$ under suitable convexity assumptions are showed in \cite{BDG20}, while the relaxation under the weaker convergence of $\B u_j\mathcal L^N \overset{\ast}{\rightharpoonup} \B u$ in $\mathcal M(\Omega;W)$ can be found in \cite{APR20} and \cite{AR20}.
Our proof expands the results contained in these two latter papers, taking advantage of the properties of the operators $\B$ and its {\it annihilator} $\mathcal A$ and the embedding results (see Section \ref{pre} for an overview). We point out that the analysis carried out in the present paper strongly relies on the interplay between the operators $\A$ and $\B$. More specifically, $\A$ is supposed to have a $\mathbb C$-elliptic potential $\B$. Therefore, many operators $\A$ are left out of our analysis. For example, $\A={\rm div}$ is excluded since its potential $\B$ is the curl operator.  We refer to \cite{CG22} for the importance of this assumption in related contexts.

It is worth noticing that our analysis may encompass more general energy densities $\mathfrak{g}(x,u,v, \B u)$ than $ f_1(u)f_2(v)+ h(\B u)$. Indeed, one may replace $h$ by $g_2$ being as in Theorem \ref{thm:RelaxationB} as well as by $g_2(x,u, \mathcal E u)$ in the specific case $\B=\mathcal E$. In this latter case, we need to exploit the results in \cite{CFVG20}, thus leading to the same type of formula in \cite[Remark 3.1]{KKZ23} for the symmetrized gradient with $K_f$ in \cite[Remark 3.1]{KKZ23} replaced by $g$ in \cite[eq. (2.6)]{CFVG20}. Moreover, in view of the embeddings of $W^{\B,1}(\Omega)$ into $W^{k-1,1}(\Omega)$ which hold when $\B$ is a $\mathbb C$-elliptic operator, one may consider more general $\mathfrak g$ in $G$. For example, it is possible to introduce a dependence also on $\nabla^2 u, \dots, \nabla^{k-1} u$. Of course, the analysis would be more complete, featuring interaction between the measure $v$ and the highest derivative $\nabla^{k-1} u$ (according to the order $k$ and $N$) and with no substantial differences in techniques, but for the sake of brevity and clarity, we have preferred to work with the energy density depending only on $\B u$.

As for $\B=D$ (see \cite{KKZ23}),  when $N=1$ and for first-order operators a separate analysis is required.  
This is due to the fact that points in dimension $N$ greater than $2$ have $1$-capacity zero, which is not the case for the real line. 
We now state the following $\B$-variant of \cite[Theorem 1.2]{KKZ23}.

 \begin{Theorem}
      \label{thm:relN=1}
      Let $\Omega=(\alpha, \beta)$ be a bounded open interval, with $\alpha<\beta$. Let $\B$ be the $\mathbb C$-elliptic order partial differential operator of order  $k=1$ of the form 
       \begin{equation}
                \label{def:ellipticoperN1}
                \B u = B{d^ku\over dx^k}, \;\;\; u:(\alpha,\beta)\subset\R \to V.
            \end{equation}
      Let $f_1, f_2, h$ be the functions satisfying assumptions ${\rm (H_1)}, {\rm(H_2q)}$ and $ {\rm(H_3p)}$, with $p=q=1$.  
      Let $ G_1(u,v)$ and $\overline{\mathcal G_1}(u,v)$ be the functionals  given by \eqref{GG1} and \eqref{G1relax}. Then,
                 \begin{align}\label{rep1}
                     \overline{\mathcal G_1}(u,v) &= \int_{\Omega} f_1(u(x))df_2^{\ast\ast}(v^{\rm diff})(x)+\int_\Omega dh^{\ast\ast}(\B u^{\rm diff})(x) + \sum_{x\in S^0}f^0_h\left(u(x^{+}), u(x^{-}), v^0( \{x\})\right)\notag\\
                     &\quad + \inf_{z\in V}f^0_h\left(u(\alpha^{+}), z, v^0(\{\alpha\}) \right) + \inf_{z\in V}f^0_h\left(z, u(\beta^{-}), v^0(\{\beta\}) \right),
                 \end{align}
    where $f_2^{\ast\ast}$ and $h^{\ast\ast}$ are the convex hulls of $f_2$ and $h$ respectively and $f^0_h: V\times V\times \R^d\to \R$ is defined by 
                 \begin{equation}
                     \label{deff0W-RelN=1}
                     f^0_h(a^{+}, a^{-}, b):= \inf_{\begin{array}{ll}
		u \in W^{\B,1}((-1,1)), \\
		v\in L^1((-1,1); \R^l)\\
		u(-1)=a^-, u(1)=a^+,\\
		\int_{-1}^{1}v dx= b
\end{array}} \left\{ \int_{-1}^1 \left(f_1(u(x))(f_2^{\ast\ast})^\infty(v(x)) + (h^{\ast\ast})^\infty(\B u(x))  \right)dx    \right\}.
                 \end{equation}
 Here, $u(\alpha^{+})$ and $u(\beta^{-})$ denote the traces of $u$ at the boundary and similarly $u(x^{+})$ and $u(u^{-})$ denote the traces of $u$ from the right and from the left at an interior point $x$.
  \end{Theorem}
 
Note that the symbol $W^{\B,1}$ in \eqref{deff0W-RelN=1} is fictitious, since $\B$ is $\mathbb C$-elliptic. Indeed, in this case, the space $W^{\B,1}(\alpha, \beta)$ coincides with $W^{1,1}(\alpha,\beta)$, for every $\alpha, \beta \in \mathbb R$. 

 In higher dimensions, when the order $k$ of the operator $\B$ is small enough, i.e. the convergence of minimizing sequences for $u$ is not strong enough, it is possible to optimize $f_1(u)$ locally around points that the measure  $v$ might charge, without a great expense of energy $h$. This is the concentration effect which might appear even if the measure $v$ is absolutely continuous with respect to $\mathcal L^N$. On the contrary, also in dimension $N=1$, if $\B$ has order $k>1$, no concentration effect is detected due to the boundedness of the approximating sequences for $u$ and the pointwise properties of $u$. 

 To complete the investigation performed in the present manuscript, in Section  \ref{relBFLhigh} we also consider functionals $G$ as in \eqref{G} in the superlinear case, i.e. $\mathfrak g$ has superlinear growth in the $v$ and $\B u$ variables, respectively. In particular, when $1<p=q$, there are no concentration phenomena, in view of well-known equi-integrability results, see \cite{FL, FMP} and also \eqref{repkpN}. Furthermore, we also present variants of Theorems \ref{thm:relk<N,N>=2}, \ref{thm:relKgeq N}  and \ref{thm:relN=1} including the case $p>1$, $u \in  W^{\B, p}$ and assuming for the density $h$ the validity of Assumption $(H_3p)$, with $p>1$. This is the subject of subsection \ref{Sobolevconc}. In details, the proof of Theorem \ref{thm:relk<N,N>=2} can essentially be followed
step by step. The only differences are that Proposition \ref{proprel**} can be avoided relying directly on Theorem \ref{the:relaxBFL},  (which is, per se, a new relaxation result in the framework of $W^{\mathcal A,p}(\Omega)$-type spaces) and the proof of the upper bound. 
In the latter, the construction of
the function $\varphi$ (the same as in the proof of Theorem \ref{thm:relk<N,N>=2}, cf. \eqref{def:varphi}) now has to be $W^{k,p}$ to still allow the cut-off arguments.  This is possible only for $p \leq N$ if $\B$ is of order $1$ and for $p=1$ if $\B$ is of order $k\geq N$, as for large $p$ any function in $W^{1,p}$ is locally bounded by embedding. The relaxation for $p > N\geq 2$ and $\B$ of order $k>N$ is more similar to Theorem \ref{thm:relKgeq N}, since the variable $u$ gets frozen due to the uniform convergence. 

 The results presented in Section \ref{relBFLhigh} extend the ones in the $\mathcal A$-free setting, contained in \cite{BFL00} reformulating them in terms of a suitable potential $\B$, related to the annihilator $\mathcal A$ (see \cite{R19, V13}).
Moreover, in our analysis, both the operators $\mathcal A$ and $\mathcal B$ need not to be necessarily of first order.  The questions addressed in Section \ref{relBFLhigh} are not motivated by the mere academic curiosity, but the models considered therein find application in many contexts. Indeed, if $\B$ is of the type $(\nabla, {\rm Id})$ acting on $u, v$, our analysis falls in the scope of modeling materials where both elastic deformation and chemical composition play a role (see \cite{CRZ1, CRZ2, FKP2}). We leave to a forthcoming research general relaxation results for the cases with mixed growth in $\B u$ and $v$: for instance, one may consider models with linear and superlinear behaviors or different superlinear behaviors. We point out that some of the results in Section \ref{relBFLhigh} are still valid under the assumption of mixed growth.  We refer to \cite{CZA,CZE, FKP1} for the case $\B=\nabla$, $p=1$ and $q\leq \infty$. The references \cite{RZ1, RZ} contain relaxation results with respect to the$L^1$-weak topology and weak$^*$ in the sense of the measures and with $p=1$ and $q\ge 1$. 
In the context of dimension reduction,  we mention for the superlinear case \cite{BFMbend, KR, LDR} where there is no explicity dependence on $u$ and for the linear case \cite{BZZ, FMZ25} where $\B$ either coincides with $\nabla$ or is a more general operator appearing in the context of non simple materials.  Finally, we refer to \cite{DD,FKP2} for problems modeling  magnetoelasticity and \cite{CL} for  models of particle inclusions in elastic matrices.

\smallskip

 The plan of the paper is the following: Section \ref{pre} contains notation and preliminary results about partial differential operators, function spaces, integrands and their envelopes. Section \ref{secRel} is devoted to the proof of Theorem \ref{thm:RelaxationB}.
This latter theorem will be exploited in Section \ref{Sect:Applications} to deduce Theorems \ref{thm:relk<N,N>=2}, \ref{thm:relKgeq N}, \ref{thm:relN=1} as well as a lower semicontinuity result in the framework of the space with variable exponents (cf. Theorem \ref{thm:asBHHZsci}). 
The parallel analysis in the Sobolev setting is the target of Section \ref{relBFLhigh}. Finally, the Appendix collects some results and proofs which will be used throughout the paper.

\section{Preliminaries}\label{pre}

{{\bf Notation.}} Throughout the manuscript, $d, j, k, K, l, m, n, N$ denote elements in $\mathbb N$.  We denote the unit cube of $\R^N$ by $Y^N$ while the unit torus in $\R^d$ by $T^d$.
 $C^\infty (T^N; \R^m)$ is the space of smooth periodic functions defined in $\mathbb R^N$ with values in $\mathbb R^m$. We denote by $B_r(x_0)$ the ball of radius $r$ and center $x_0$. 
For an open set $\Omega$ of $\R^N$, $\mathcal O_r(\overline\Omega)$ denotes the family of all subsets of $\overline\Omega$ which are open with respect to the relative topology of $\overline\Omega$.
For $A\in \mathcal O_r(\overline\Omega)$, $\mathcal M(A;\mathbb R^l)$ is the dual space of the continuous functions $\varphi:A\to \mathbb R^l$ with $\varphi=0$ on $\partial A\cap\Omega$.  For every $\mu \in\mathcal M(A;\mathbb R^l)$, 
we have the Radon-Nikodym decomposition $\mu= \mu^a + \mu^s= {d\mu\over d\mathcal L^N}+ \mu^s$, where $\mu^a={d\mu\over d\mathcal L^N}$ is the so-called Radon-Nikodym derivates of $\mu$ with respect to $\mathcal L^N$. Sometimes, where no ambiguity arises, we simply identify $\mu^a$ with $\mu$. By $\mu_j \overset{*}{\rightharpoonup} \mu$, as $j\to +\infty$ we denote the weak* convergence of the sequence of measures $(\mu_k)_k \subset \mathcal M(A;\mathbb R^l)$ towards $\mu \in \mathcal M(A;\mathbb R^l)$. 

Set $\mathbb K=\{\R, \mathbb C\}$.  $V,W,X$ are vector spaces of finite dimensions on $\mathbb K$. By $\Mcal(\Omega; W)$ we denote the set of finite $W$-valued Radon measures on $\Omega$.
$\Lin (W, X)$ denotes  the space of linear maps  $ W \to X$ and $X\odot^l W$ stands for the space of $X$-valued symmetric $l$-linear maps on $W$. This is naturally the space of the $l$-th gradients, i.e., $D^l f(x) \in X \odot^l W$ for $f \in C^l(W, X)$,  $x \in W$. The tensor product $a\otimes b$, with $a,b\in \R^N$, is defined by  $a \otimes b =(a_ib_j)$. In particular, $\otimes ^la := a\otimes \dots \otimes a$, where $a$ appears $l$ times on the right-hand side. We also have that $\widehat \nabla^lf(\xi) =\hat f(\xi) \otimes^l \xi \in X\odot^l W$ for $f\in \mathcal S(W, X)$, $\xi \in W$, with  $S(W, X)$ being the Schwartz space. With $\mathcal{D}(\R^N;V)$, we denote the space $C_0^\infty(\R^N;V)$ and $\mathcal{D}'(\R^N;V)$ is the dual of $\mathcal{D}(\R^N;V)$.  With the aim of emphasizing classical  derivatives, we will replace the symbol $D$ by $\nabla$. Eventually,  for a function $f: \Omega\times \R^l\to \R$ with recession function $f^\infty(x, \cdot)$ (see \eqref{finfty}) in the second variable and for a measure $\mu\in \mathcal M (\Omega; \R^l)$, we define
     \begin{equation}
         \label{def:df}
         df(x, \mu)(x) \coloneqq f\left(x, {d\mu^a\over d \mathcal L^N}(x)\right) dx + f^\infty\left(x, {d\mu^s\over d |\mu^s|}(x)\right) d|\mu^s|(x). 
     \end{equation}

\subsection{Partial Differential Operators}\label{PDO}

We consider $k$-homogeneous (i.e., homogeneous and of order $k$)\; partial differential operators $\mathcal B$ (or $\mathcal A$)  on $\R^N$ from $V$ to
$W$ (or from $W$ to $X$) with constant coefficients of the form
\eqref{Budef}, i.e.,
\begin{align*}
\mathcal B u=\sum_{|\beta|=k}B_\beta \partial^\beta u, \;\;\;
\mbox{for } u\in \mathcal{D'}(\R^N; V),
\end{align*}
with $B_\beta \in {\rm Lin}(V, W)$, $\beta \in \mathbb N_0^N$ being a multi-index with  order $|\beta| = \beta_1 +\dots +\beta_N$ and $\partial _\beta$ representing the distributional derivative $\partial_1^{\beta_1}\dots
\partial_N^{\beta_N}$. Denoting by $\B^\ast$ the formal $L^2$-adjoint operator of $\B$ given by 
        \begin{equation}
            \notag
            \B^\ast := (-1)^k\sum_{|\alpha|=k} B_\alpha^\ast \partial^\alpha, 
        \end{equation}
the definition of $\mathcal B u$ in \eqref{Budef} is easily understood in the sense of distributions (see \cite{APR20}). 

We say that the operator $\mathcal{B}$ has  {\it constant rank} if  there exists a positive integer $r$ such that
\begin{align}\notag
{\rm rank} \mathbb B(\xi ) = r\,\,\, \hbox{ for all } \xi \in  \mathbb K^N \setminus\{0\},
\end{align}
where the tensor-valued $k$-homogeneous polynomial
\begin{align}\label{mathbbAdef}
\mathbb B(\xi ) := \sum_{|\beta|=k}
\xi^\beta B_\beta   \in   \Lin(V; W),\;\;\; \xi  \in \mathbb K^N,
\end{align}
is the principal symbol associated to the operator $\mathcal B$. Here, $\xi^\beta  := \xi_1^{\beta_1}\dots \xi_N^{\beta_N}$.
In other words, $\mathcal B$ has {\it constant rank} if ${\rm dim}({\rm Im}(\mathbb B(\xi)))$ is a constant $r$ independent of $\xi \in  \mathbb R^N \setminus \{0\}$.
The symbol $\mathbb{B}$ as in \eqref{mathbbAdef} induces a bilinear pairing $\otimes_{\symB}: V\times \mathbb K^N\to W$ defined by 
           \begin{equation}\notag
               v\otimes_{\symB} z := \mathbb B(z)v= \sum_{|\beta|=k} z^\beta B_{\beta} v, \qquad \mbox{for } v\in V, z\in \mathbb K^N. 
           \end{equation}
Note that the map $\otimes_{\B}$ extends the standard tensor product notation to more general differential operators $\mathcal  B$ than the first-order derivative $D$. In this latter case $\otimes_{\mathbb B}=\otimes$. For all $\varphi\in C^\infty_{\rm c}(\R^N)$ and $v\in C^\infty(\R^N; V)$,
      \begin{equation}
          \notag
          \B(\varphi v) = \varphi\B v + [\symB, \varphi](v),
      \end{equation}
  where $[\symB, \varphi]$ is called the commutator of $\B$ on $\varphi$ acting on distributions $\eta\in \mathcal{D}'(\R^d, W)$ as
     \begin{equation*}
         [\symB, \varphi](\eta)= \B(\varphi\eta) - \varphi\B\eta.
     \end{equation*}
 Note that $[\symB, \varphi]$ is a partial differential operator of order $(k-1)$ from $W$ to $X$, whose coefficients are smooth and depend only on the ones of the principal symbol $\symB$ and the first $k$ derivatives of $\varphi$, see \cite[eq. (34) and (35)]{AR21} for more details and related formulas.
       
We, now,  recall other useful properties of the operators $\B$ which will be crucial in the sequel. 
 \begin{Definition}\label{Kellipticitydef}
   Let $\mathbb K\in \{\mathbb R,\mathbb C\}$. We say that a $k$-homogeneous partial differential operator $\B$ on $\R^N$ from $V$ to $W$ of the form \eqref{Budef} 
      is {\it $\mathbb K$-elliptic} if $\ker(\symB [\xi])=\{0\}$ for any $\xi\in\mathbb K^N\setminus\{0\}$.  Equivalently there exists a positive constant $C$ such that
         $|\symB (\xi)v| \geq C |\xi|^k|v|$ for all $\xi \in \mathbb K^N$ and all $v \in V$;

 \end{Definition}
If $\mathcal B$ is a $1$-homogeneous partial differential operator $\mathbb K$-elliptic, from Definition \ref{Kellipticitydef}, it follows that
 there exist two constants
$0 < k_1 \leq k_2 < \infty$ such that
\[k_1|v| |\xi| \leq |v \otimes_{\mathbb B} \xi| \leq k_2|v| |\xi|
\] for all $v \in V$ and $z \in {\mathbb K}^N$ (see \cite[Lemma 2.3]{BDG20}). We have that operators which are $\mathbb C$-elliptic are $\mathbb R$-elliptic but the converse is not true (see, e.g., \cite[Example 2.2 (c)]{BDG20}). In addition, a $\mathbb C$ (or $\mathbb R$)-elliptic operator satisfies the constant rank properties by the Rank-Nullity Theorem.  

The importance of constant rank operators is explained in the next property (see also \cite[Lemma 3.2]{CG22} for first-order operators) which will be crucial later on.

\begin{Proposition} \label{prop:ExistenceofAandB}
\begin{itemize} 
    \item[{\rm (a)}]{\rm {\cite[Proposition 4.2]{V13}}} Let $\B$ be a   $\mathbb R$-elliptic partial differential operator as in \eqref{Budef}. Then, there exist $l\in\mathbb{N}$, a real finite dimensional vector space $X$ and a $l$-th order $X$-valued {\it constant-rank} differential operator $\A$ of the form 
          \begin{equation}
              \notag
              \A v := \sum_{|\alpha|=l } A_\alpha \partial^\alpha v, \qquad v:\R^N\to W,
          \end{equation}
    being $A_\alpha \in \Lin (W; X)$ such that the associated Fourier symbol complex
           \begin{equation}
               \notag
               V \xrightarrow{\symB(\xi)} W  \xrightarrow{\symA(\xi)} X 
           \end{equation}
    is exact at $W$, i.e., 
    \begin{equation}\label{4R19}
    {\rm Ker}(\symA(\xi))= {\rm Im} (\symB(\xi)), \hbox{ for any }\xi\in\R^N\setminus\{0\}.
\end{equation}
    \item[{\rm (b)}]{\rm {\citetext{\citealp[Theorem 1]{R19} and see also \citealp{R24}}}}  Let $\A$ be a linear, homogeneous differential operator with constant coefficients on $\R^N$,  with principal symbol $\symA$. Then, $\A$ has constant rank if and only if there exists a linear, homogeneous differential operator $\B$ with constant coefficients on $\R^N$ (with principal symbol $\symB$) such that  \eqref{4R19} holds. 
   
\end{itemize}
\end{Proposition}
Adopting the terminology of \cite{R19, V13}, we refer to $\A$ as the {\it annihilator} for $\B$ and $\B$ as the {\it potential} for $\A$. 
Moreover, the operator $\A$ in Proposition \eqref{prop:ExistenceofAandB} (a) may not be  $\mathbb C$-elliptic.

\begin{Remark}\rm
The orders of the operators $\A$ and $\B$ considered in Proposition \ref{prop:ExistenceofAandB} might be different from each other, as showed in the following examples.
   \begin{itemize}
   \item $\B u= 1/2(Du + Du^\intercal)$ and $\A v=\mbox{curlcurl}v$ (see, e.g., \cite[Example 3.10(e)]{FM99}), where $\B$ is a first order operator and $\A $ is of second order. 
   \item $\B u = D^2 u$, with $u :\Omega\subset \mathbb R^2 \to \mathbb R$ and $\mathcal A =\mbox{curl}$ with  symbol $\mathbb A $ defined on the space of symmetric $2\times 2$ matrices 
   as follows
   $\mathbb A \xi:=\sum_{\alpha=1}^2 A_\alpha \xi^\alpha   \in W= \mathbb R^{3}\equiv \mathbb R^{2 \times 2}_{\rm sym}$, where the latter space coincides with the space of symmetric bilinear maps from $\mathbb R^2 \times \mathbb R^2$ and with
     \begin{equation*}
         A_1\coloneqq \left[
   \begin{array}{llll}
        0 & 0 & 0 &0  \\
        0 & -1 & 0 &0  \\
        0 & 0 & 0 &0  \\
        0 & 0 & 0 & -1  
   \end{array}
   \right] \,\,\,\,\, \hbox{and} \,\,\,\,\,
   A_2 \coloneqq\left[\begin{array}{llll}
   0&0&0&0\\
   1&0&0&0\\
   0&0&0&0\\
   0&0&1&0
   \end{array}\right].
     \end{equation*}
  Let $v:\mathbb R^2 \to \mathbb R^4_{\rm sym}$ be given by $v=\left[\begin{array}{llll}
   v_{11}, \; v_{21}, \; v_{12}, \; v_{22}
   \end{array}\right]^T$, with $v_{12}=v_{21}$.
   In this setting, $\A$ is a first-order operator and $\B$ is a second-order one, (see, e.g., \cite[Example 3.10(d)]{FM99}). In fact, $\A v=0$ is equivalent to $$\left\{
   \begin{array}{ll}
   \frac{-\partial v_{21}}{\partial x_1}+ \frac{\partial v_{11}}{\partial x_2}=0,\\
   \frac{\partial v_{12}}{\partial x_2}-\frac{\partial v_{22}}{\partial x_1}=0.
   \end{array}
   \right.$$
  This guarantees that $v=\left[
   \begin{array}{llll}
   \frac{\partial u^1}{\partial x_1}, \;
   \frac{\partial u^1}{\partial x_2},\;
   \frac{\partial u^2}{\partial x_1},\;
   \frac{\partial u^1}{\partial x_2}
   \end{array}\right]^T,$
   which together with the fact that $v_{12}=v_{21}$, provides
   that $u^1=\left[\begin{array}{ll}\frac{\partial^2 u}{\partial x_1^2},\;\frac{\partial^2 u}{\partial x_1\partial x_2}\\
   \end{array}\right]^T$
   and $u^2= \left[\begin{array}{ll}\frac{\partial^2 u}{\partial x_1\partial x_2},\;
   \frac{\partial^2 u}{\partial x_2^2}
   \end{array}\right]^T.$
   \end{itemize}
\end{Remark}

Lastly, we recall the notion of {\it wave cone} (also named {\it characteristic cone} or {\it $\Lambda$ cone}) associated to $\mathcal A$, which plays a fundamental
role for the study of $\mathcal A$-free fields, as discussed in the work of Murat and Tartar \cite{M1,M2,M3,MT,T2,T1,T3,T4}
\begin{align}\label{Lambdacone}
\Lambda_{\mathcal A} =
\bigcup_{\xi \in \mathbb R^N \setminus \{0\}}
{\rm ker} \mathbb A(\xi ) \subset W.
\end{align}
 When $\mathcal B$ is $\mathbb R$-elliptic, by Proposition \ref{prop:ExistenceofAandB}, the above {\it wave cone} coincides with the {\it image cone} $\mathcal I_{\mathcal B}$ defined as
\begin{align}\label{imagecone}
\mathcal I_{\mathcal B}:=
\bigcup_{\xi \in \mathbb R^N\setminus \{0\}}
{\rm Im}( {\mathbb B}(\xi) ) \subset W,
\end{align}
introduced in \cite[page 248]{AR21} (equivalently the $\mathbb B$-rank-one cone defined in \cite{BDG20} for the first-order operators). The equality between $\Lambda_\mathcal A$ and $\mathcal I_{\mathcal B}$ is due to \eqref{4R19}.  In many cases this span is the entire vector space, but in general this is not the case as remarked in \cite{FM99}.

\subsection{Function Spaces}\label{FS}

We recall here the definition and the main properties of the function spaces $W^{\B, p}$ and $BV^{\B}$, associated to  a $k$-th order linear partial differential operator $\B$ of the form \eqref{Budef}, introduced in \cite{BDG20, GR19}. 
Let $\Omega$ be an open and bounded subset of $\R^N$ with Lipschitz boundary.
For $1\leq p\leq \infty$, the Sobolev space $W^{\B, p}(\Omega)$ is defined by
     \begin{equation}
         \notag
         W^{\B, p}(\Omega) := \left\{ u\in L^p(\Omega; V) : \B u\in L^p(\Omega; W)   \right\}.
     \end{equation}
It is a Banach space endowed with the norm
        \begin{equation}
            \notag
            \|u\|_{W^{\B, p}(\Omega)}:= \|u\|_{L^p(\Omega; V)} + \|\B u\|_{L^p(\Omega; W)}.
        \end{equation}
In the case $\B = \nabla^k$,  we simply write $W^{k, p}(\Omega; \R^m)$. We also introduce the function space $V^{\B, p}(\Omega)$ 
        \begin{align}\nonumber
            V^{\mathcal B, p}(\Omega) := \{u \in W^{\mathcal B,p}(\Omega): \nabla^l u \in L^p(\Omega, V \odot^l \mathbb R^N), l = 1,\dots k - 1\},
        \end{align}
endowed with the norm
\begin{equation}
\notag
            \|u\|_{V^{\B, p}(\Omega)}:=\|u\|_{L^p(\Omega; V)} + \sum_{l=1}^{k-1}\|\nabla^l u\|_{L^p(\Omega;V \odot^l \mathbb R^N)}+\|\B u\|_{L^p(\Omega;W)}.
            \end{equation}    
    In other words, $V^{\mathcal B, p}(\Omega)= W^{\mathcal B, p}(\Omega)\cap W^{k-1, p}(\Omega; V)$.           
We gather now some approximation results proven in \cite{GR19} involving the above function spaces. 
To this end, we give a preliminary definition.

\begin{Definition}\label{defOmega}
A connected open set $\Omega \subset \mathbb R^N$ is called a $C^0$-domain if for any $x \in \partial \Omega$ there exists a neighborhood ${\mathcal N}(x)$ of $x$, relatively open in $\Omega$, a coordinate system in $\mathbb R^N$ and a continuous function $f$ such that, in the new coordinates $(x', x_N)$, $\mathcal N(x)=\{(x', x_N):0 <x_N<f(x'), x' \in B_1(0)\}.$ 
\end{Definition}

The following results are proven in \cite{GR19} as a unique theorem in which they appear as several steps, but we prefer to recall them separately because of their interest.

\begin{Theorem}[{\cite[Lemma 2.7]{GR19}}]
\label{thm:densitynormandarea}
    Let $1\leq p <+\infty$, let $\Omega$ be an open subset of $\R^N$ and $\mathcal B$ be a linear partial differential operator  of order $k$ as in \eqref{Budef}. We have that  
       \begin{itemize}
       \item[\rm (i)]$C^\infty(\Omega;V)\cap V^{\B, p}(\Omega)$ is dense in $V^{\B, p}(\Omega)$ with respect to the norm topology;
           \item[{\rm (ii)}] $C^\infty(\Omega; V)\cap W^{\B, p}(\Omega)$ is dense in $W^{\B, p}(\Omega)$ with respect to the norm topology. 
       \item[\rm (iii)] if $\Omega$ is a  bounded $C^0$-domain, then $C^\infty(\overline \Omega;V)$
is dense in $W^{\B,p}(\Omega)$. The same holds true for $V^{\mathcal B,p}(\Omega)$.
       \end{itemize}
\end{Theorem}

It is worth noting that in Theorem \ref{thm:densitynormandarea} we do not require that $\B$ is $\mathbb C$-elliptic. Such a condition comes into play in \cite{BDG20} when dealing with traces.

The following embedding proposition enables us to have a control of the derivatives of  $u\in W^{\B, p}(\Omega)$ of any order.
  \begin{Proposition}[{\cite[Lemma 2.8]{GR19}}]\label{propWB=VB}
      Let $\B$ be given by \eqref{Budef} satisfying the $\mathbb C$-ellipticity property. Let $\Omega$ be a bounded Lipschitz domain in $\R^N$. Then, $W^{\B, p}(\Omega)\simeq V^{\B, p}(\Omega)$, for any $1\leq p\leq \infty$. 
  \end{Proposition}

The space $BV^\B(\Omega)$ of functions of bounded $\B$-variations is given by
      \begin{equation}\notag
          BV^\B (\Omega):=\{u\in L^1(\Omega; V) : \B u\in \Mcal(\Omega; W) \},
      \end{equation}   
which is a Banach space endowed with the following norm
      \begin{equation}
          \notag
          \|u\|_{BV^\B}:= \|u\|_{L^1(\Omega; V)} + |\B u|(\Omega),
      \end{equation}

\noindent      where $|\cdot|:\mathcal M(\Omega;W)\to [0,+\infty)$ stands for the total variation of a measure.
If $\B$ is a $k$-th order $\mathbb C$-elliptic differential operator, the space $BV^\B$ is then a generalization of $\nabla^k B$ spaces introduced in \cite[eq. (2.45) and (2.46)]{DT} (see also the spaces $BV^k$ introduced in \cite{ADC}, and recall that, when $k=2$, $\nabla^2 B= BV^2 = BH$ where $BH$ is the space of functions with bounded hessian, see \cite{D}). 

Observe also that, when $\mathcal B$ is $\mathbb C$-elliptic, by \cite[Theorem 3.7]{DG24}, we have that 
\begin{align}\label{BVBequiv}BV^\B(\Omega)= W^{k-1, 1}(\Omega; V) \cap BV^\B(\Omega).
\end{align}
More precisely, following \cite{GR19} we recall the subsequent result (cf. \cite[Theorem 2.3]{PPRV}).
\begin{Theorem}\label{thm2.3}
Let $\Omega \subset \mathbb R^N$
be a bounded Lipschitz domain. A  $k$-th order, homogeneous
and linear differential operator $\mathcal B$ with constant coefficients  as in \eqref{Budef} is $\mathbb C$-elliptic if
and only if there exists a positive constant $C$ such that
\begin{equation}\notag
\|u\|_{W^{k-1, {\frac{N}{ N-1}}}(\Omega; V)}
\leq C(|\mathcal B u|(\Omega) + \|u\|_{L^1(\Omega)}), \;\;\;\hbox{for } u \in BV^{\mathcal B}(\Omega).
\end{equation}
\end{Theorem}
In principle, one could define the {\it $\B$-variation }of $u\in L^1(\Omega; V)$ via duality by 
                  \begin{equation}\label{VBu}
                      V_{\mathcal B}(u,\Omega) := \sup\biggl\{ \int_\Omega u\cdot \B^\ast \varphi dx :  \varphi \in C^1_c(\Omega; W), \|\varphi\|_{L^\infty(\Omega;W)}\leq 1 \biggr\}.
                  \end{equation}
As in the standard case of functions of bounded variation, one can show that $u\in BV^\B(\Omega)$ if and only if \eqref{VBu} is finite and that it coincides with $|\mathcal Bu|$, where the latter is interpreted as the total variation of the measure $\mathcal B u$. The proof follows the same arguments as \cite[Proposition 3.6]{AFP00}. Indeed, if $u\in BV^\B(\Omega)$, we easily deduce that 
       \begin{align}
           \notag
           \left|\int_{\Omega} u\cdot \B^\ast \varphi dx   \right|\leq |\B u|(\Omega).
       \end{align}
Taking the supremum over $\varphi\in C^1_c(\Omega; W)$ with $\|\varphi\|_{L^\infty(\Omega; W)}\leq 1$, we conclude that $V_{\mathcal B}(u; \Omega)<\infty$. On the other hand, if $V_{\mathcal B}(u; \Omega)<\infty$, an application of the Riesz representation theorem alongside the density of $C^1_{\rm c}(\Omega; W)$ in $C^0(\Omega; W)$ leads us to conclude that $u\in BV^\B(\Omega)$.     

We now introduce several notions of convergence in $BV^\B$. In many applications, the norm convergence turns out to be too strong hence other weaker notions of convergence in $BV^\B$ come into play. Let $u\in BV^\B(\Omega)$ and $(u_j)_j\subset BV^\B(\Omega)$. We say that 
\begin{itemize}
\item[(i)] $(u_j)_j$ converges to $u$ in the weak*–sense (in symbols $u_j\overset{\ast}{\rightharpoonup}
u$) if and only if $u_j \to u$
strongly in $L^1(\Omega;V)$  and 
$\mathcal B u_j \overset{\ast}{\rightharpoonup}
\mathcal B u $ in the weak*–sense of $W$–valued Radon
measures on $\Omega$
 as $j \to +\infty$;
\item[(ii)] $(u_j)_j$ converges to $u$ in the strict sense (in symbols $u_j \overset{s}{\to} u$)
 if and only if $u_j\to u$ in $L^1(\Omega; V)$  and $|\B u_j|(\Omega) \to |\B u|(\Omega)$;

\item[(iii)]  $(u_j)_j$ converges to $u$ in the area-strict sense if and only if  $u_j\to u$ in $L^1(\Omega; V)$ and 
           \begin{equation}
               \notag
               \lim_{j\to\infty} \langle \B u_j\rangle(\Omega) = \langle \B u\rangle (\Omega),
           \end{equation}
            where $\langle\mu\rangle: \mathcal M(\Omega;W)\times \mathcal B(\Omega)\to [0,\infty)$ denotes the area functional, i.e.
    \begin{equation}\label{areafunctional}
    \langle\mu \rangle(U):= \int_{U}\sqrt{1+|\mu^a(x)|^2}dx + |\mu^s|(U),
    \end{equation}
    for every $\mu \in \mathcal M(\Omega;W)$ and Borel set $U\subset\Omega$. 
\end{itemize}

\begin{Remark}\label{compactnessBVA}
\begin{itemize}
    \item[\rm (i)] As already observed in {\rm \cite{BDG20}}, following the same argument of {\rm \cite[Theorem 3.23]{AFP00}}, for a bounded extension domain $\Omega$ (or a bounded open set $\Omega$ with Lipschitz boundary) a bounded sequence in $BV^\B$ norm admits always a subsequence which converges in the sense of $(i)$, whenever $\B$ is a first order operator. Furthermore, arguing as in {\rm \cite[Corollary 3.49]{AFP00}}, the strong convergence holds in $L^q$ for every $1\leq q< 1^*$.
The Sobolev inequality given by Theorem \ref{thm2.3} and the above argument allow us to deduce that when $\B$ is an operator of order $k$, then bounded sequences in $BV^\B$ are strongly converging in $W^{k-1,q}$, with $1\leq q< 1^*$.
\item[\rm (ii)] For a bounded Lipschitz domain $\Omega$, in view of Theorem \ref{thm2.3}, the strong $L^1$ convergence in (i), (ii), (iii) above can be replaced by the strong convergence in $W^{k-1, 1}(\Omega; V)$. It is also easily seen that $(iii)\Rightarrow(ii) \Rightarrow (i)$. Recall that the $\mathcal B$-variation is sequentially lower semicontinuous with respect to weak*-convergence (for measures). More precisely, if $u_j \overset{\ast}{\rightharpoonup} u$ then
$|\mathcal Bu|(\Omega)\leq \liminf_{j\to +\infty}
|\mathcal B u_j|(\Omega).$
 Moreover,
if $u_j \in BV^{\B}(\Omega)$ is a bounded sequence with $u_j \rightharpoonup u$ in $L^1(\Omega;\R^d)$
then $u_j \overset{\ast}{\rightharpoonup} u$ in the sense (i). 
\item[\rm (iii)] If $\Omega$
 is open and bounded with Lipschitz boundary, then it is
easy to conclude by Banach–Alaoglu's theorem that if $(u_j) \subset BV^{\B}(\Omega)$
is uniformly bounded
$BV^{\B}$–norm, i.e. $\sup_{j} \{\|u_j\|_{L^1}+ |\mathcal B u_j|\}<+\infty$, then there exists $u \in BV^\B(\Omega)$
 and a subsequence $(u_{j_s})_s$ of $(u_j)_j$ such that
$u_{j_s} \overset{\ast}{\rightharpoonup} u$ as $s \to +\infty$, i.e. the so-called
weak*–compactness principle (for $BV^\B$) holds.
\end{itemize}

\end{Remark}

We state a first approximation result which holds only for operators of first order.

   \begin{Theorem}[{\cite[Theorem 2.8]{BDG20}}]\label{thm2.8}
       Let $\B$ be the operator given by \eqref{Budef} of the first order. Then $(C^\infty \cap BV^\B)(\Omega)$ is dense in $BV^\B(\Omega)$ with respect to the area-strict topology.
   \end{Theorem}
The next proposition provides a smooth approximation result in area strict sense up to the boundary.
\begin{Proposition} [{\cite[Lemma 4.15]{BDG20}}]
\label{prop:densityuptoboundary}
    Let $\B$ be a first-order, $\mathbb C$-elliptic operator as in \eqref{Budef}. Let $u\in BV^{\B}(\Omega)$. Then, there exists a sequence $(u_j)_j\subset C^\infty(\overline{\Omega})$ such that $u_j$ converges to $u$ in the area-strict sense.
    The result trivially extends to vector-valued fields.
\end{Proposition}

The following result shows a strictly-density result for $BV^\B$-functions, with $\B$ an operator of arbitrary order.
\begin{Proposition}[{\cite[Lemma 6]{RS20}}] 
\label{Prop:densityinstrictsense}
    Let $\B$ be a $k$-th order operator as in \eqref{Budef}. Let $\Omega$ be an open set of $\R^N$ and $u\in BV^\B(\Omega)$. Then, there exists a sequence $(u_j)_j\subset C^\infty(\Omega;V)\cap V^{\B,1}(\Omega)$  converging $\B$-strictly to $u$  and such that $u_j \to u$ in $W^{k-1,1}(\Omega)$.
\end{Proposition}

\begin{Remark}
\label{remRSlemma6}
The strong convergence of $(u_j)_j\subset W^{k-1,1}(\Omega)$ is explicitly proven in \cite[Lemma 6]{RS20}. Clearly, this latter convergence could be obtained from $\mathcal B$-strict convergence when $\mathcal B$ is $\mathbb C$-elliptic, as a consequence of \eqref{BVBequiv}.
\end{Remark}

The next result improves the statement of Proposition \ref{Prop:densityinstrictsense}, showing an approximation result for $BV^\B$ in area-strictly sense. The proof follows a direct argument as in \cite[Theorem 2.2]{DT}, thus being valid for operators of any order. 

\begin{Remark}\label{remAAR}
The result below explicitly proves an approximation at the level of the fields and not only in terms of their $\B$ gradients as in \cite[Theorem 1.6]{AR21}. Other possible strategies to deduce the same result could rely on  backwards arguments based on \cite[Theorem 1.6]{AR21} itself or iterative procedures in view of the same result for first order operators (cf. Theorem \ref{thm2.8}), taking into account that any operator of order $k$ may be viewed as a first-order operator (cf \cite[Theorem 20]{ARS25}).
\end{Remark}

    \begin{Proposition}
    \label{prop:densityareastrictly}
        Let $\Omega$ be a bounded open set of $\R^N$. Let $\B$ a $k$-th order operator as in \eqref{Budef} and $u\in BV^\B(\Omega)$.  There exists a sequence $(u_j)_j\subset C^\infty(\Omega;V)\cap V^{\B,1}(\Omega)$  such that $u_j$ converges to $u$ in area-strict sense on $\Omega$ and $u_j \to u$ in $W^{k-1,1}(\Omega;V)$.
    \end{Proposition}
\begin{proof}
 We are going to refine the same arguments adopted in the proof of Proposition \ref{Prop:densityinstrictsense} hence it suffices to show that for any $\delta>0$ there exists $u_\delta\in C^\infty(\Omega;V)\cap V^{\B,1}(\Omega)$  such that 
             \begin{equation*}
                 \langle \B u_\delta \mathcal L^N\rangle \leq \langle \B u\rangle +\delta,
             \end{equation*}
             where the area functional $\langle \cdot \rangle$ is evaluated in $\Omega$. 
    \\
    Fix $\delta>0$. Let $r>0$ be fixed such that $\Omega_0 \coloneqq \{x\in\Omega: \hbox{dist}(x, \partial\Omega)\geq 1/r \}$  is such that 
             \begin{equation}
                 \label{eq:intBuonboundary}
                 \int_{\partial\Omega_0} |\B u| =0,
             \end{equation}
            \begin{equation}
            \label{eq:intomega-omega0bu}
                \hbox{meas}(\Omega\setminus\Omega_0) \leq\delta \;\;\; \hbox{and}\;\;\; \int_{\Omega\setminus\Omega_0} |\B u|\leq \delta.
            \end{equation}
    Define $\Omega_{-1}\coloneqq\emptyset$ and for any integer $j$
         \begin{align}
             \Omega_{j}&\coloneqq\left\{x\in\Omega : \hbox{dist}(x, \partial\Omega) > {1\over r+ j} \right\}, \notag\\
             A_1&\coloneqq \Omega_2 \;\;\;\hbox{and}\;\;\; A_j\coloneqq\Omega_{j+1}\setminus\Omega_{j-1}, \;\;\; \hbox{for } j\geq 2.\notag
         \end{align}
    The family of set $(A_j)_j$ is an open cover of $\Omega$ with the property that any point of $\Omega$ belongs to at most two of the sets $A_j$. Let $(\varphi_j)_j$ be a $C^\infty$-smooth partition of unity subordinate to the cover $A_j$, i.e., $\varphi_j\in C^\infty_{\rm c}(A_j)$, $0\leq \varphi_j\leq 1$ and $\sum_{j=1}^\infty \varphi_j = 1$. Let 
    $\rho\in C^\infty_{\rm c}(B)$ ($B$ denotes the unitary ball in $\R^N$) be such that $\rho\geq 0$ and $\int_{\R^N}\rho(x) dx=1$. Set $\rho_{\e}\coloneqq \e^{-N} \rho(x/\e)$. Define $u_\delta$ as
         \begin{equation}
             \notag
             u_\delta \coloneqq \sum_{j=1}^\infty \rho_{\e_j}\ast (\varphi_j u),
         \end{equation}
      where $\e_j$ is a decreasing sequence going to $0$ as $j\to\infty$ and is chosen sufficiently small such that 
                 \begin{align}
                     &\left| \int_A |\rho_{\e_1}\ast \varphi_1 \B u|dx - |\varphi_1 \B u| (A)   \right| \leq \delta, \;\;\; \hbox{for }  A=\{\Omega, \Omega_0\},\label{eq:estsatisfiedbyepsilon1}\\
                     &\|\rho_{\e_j}\ast(u\varphi_j)-u\varphi_j\|_{W^{k-1,1}(\Omega;V)} \leq 2^{-j}\delta, \label{2.18bis}\\
                     &\int_\Omega |\rho_{\e_j}\ast (\B(\varphi_j u) -\varphi_j\B u ) dx - (\B(\varphi_j u) - \varphi_j\B u )  | dx\leq 2^{-j}\delta,\label{2.18tris}\\
                     &\left|  \int_{\Omega_0} \left(\langle \rho_{\e_1}\ast (\varphi_1\B u)\rangle dx - \langle \varphi_1\B u\rangle(\Omega_0)\right)  \right| \leq \delta,\label{eq:estsatisfiedbyepsilon}
                  \end{align}
                   where in \eqref{2.18bis} it has been exploited the fact that the commutator $[\mathcal B,\varphi_i ]u  :=\B(\varphi_j u) -\varphi_j\B u= 0$.
         Clearly $u_\delta \in C^\infty(\Omega;V)$.
 Furthermore, the sequence $\e_j$ satisfies $(\Omega_2\setminus\overline{\Omega}_1) + B(0, \e_1)\subset \Omega_3\setminus\overline{\Omega}_0,$ $A_2+ B(0, \e_2) \subset \Omega_4\setminus \overline{\Omega}_1$ and $ A_j+ B(0, \e_j) \subset A_{j-1}\cup A_j\cup A_{j+1},$ for $j\geq 1\notag$, where $B(0, \e_j)$ denotes the ball of center $0$ and radius $\e_j$.
       Moreover, as in \cite{RS20}, by \eqref{2.18bis}, we have
    $\|u_\delta - u\|_{W^{k-1,1}(\Omega;V)}\leq \delta,$ and
 since  mollification decreases the total variation, it results, by \eqref{2.18tris} that
\begin{align}\label{estProp2.11}
|\B u_\delta(\Omega)|\leq \sum_i 
|\varphi_i\B u\ast\varrho_{\varepsilon_j}|(\Omega)+
\sum_i
|[\B,\varphi_ ]u\ast \varrho_{\varepsilon_i}- [\B,\varphi_i ]u|(\Omega)|\leq |\B u|(\Omega) + \delta,
\end{align}
        which indeed recovers the statement of Proposition \ref{Prop:densityinstrictsense}. Since $\langle\cdot\rangle$ is a positive measure as observed in \cite[(1.8)]{DT}, we deduce that 
             \begin{align}
                 \left|  \langle \B u\rangle - \langle \B u_\delta\rangle     \right| &=\left|\langle\B u\rangle(\Omega) - \int_\Omega \langle\B u_\delta \mathcal L^N\rangle \right|\notag\\
                 &\leq \left| \langle \B u \rangle(\overline{\Omega}_0) -\int_{\overline{\Omega}_0} \langle\B u_\delta \mathcal N \rangle  \right| + |\langle \B u \rangle|(\Omega\setminus\overline{\Omega}_0) + \int_{\Omega\setminus\overline{\Omega}_0} \langle \B u_\delta \mathcal L^N\rangle.\label{estareafunc}
             \end{align}
        We estimate the first integral on the right-hand side of \eqref{estareafunc}. To this end, note that $\B u= \varphi_1 \B u$ on $\Omega_0$ and hence $\langle \B u \rangle= \langle\varphi_1 \B u\rangle$ on $\Omega_0$. In addition, $\B u_\delta= \rho_{\e_1}\ast \B (\varphi_1 u) = \rho_{\e_1}\ast (\varphi_1\B u)$ on $\Omega_0$ and $\langle \B u_\delta \mathcal L^N\rangle= \langle(\rho_{\e_1}\ast(\varphi_1\B u)\rangle$ on $\Omega_0$. Therefore, thanks to \eqref{eq:intBuonboundary} as well as \eqref{eq:estsatisfiedbyepsilon}, it follows that
                 \begin{align}
                    & \left| \langle \B u \rangle(\overline{\Omega}_0) - \langle\B u_\delta \mathcal L^N\rangle\right|(\overline{\Omega}_0) \notag = \left| \langle\varphi_1\B u \mathcal L^N \rangle(\overline{\Omega}_0) - \langle(\rho_{\e_1}\ast (\varphi_1\B u)\mathcal L^N\rangle ( \Omega_0 )  \right|\leq\delta.\notag
                 \end{align}
        Since $0\leq \langle\mu\rangle\leq k_1(\mathcal L^N+|\mu|)$  and exploiting \eqref{eq:intomega-omega0bu}, we get that
                \begin{align}
                  |\langle \B u \rangle|(\Omega\setminus\overline{\Omega}_0)&\leq 
  k_1\int_{\Omega\setminus \overline{\Omega}_0} 1 dx + |\B u|(\Omega\setminus \overline{\Omega}_0)\leq  2 k_1\delta,\notag\\
                   \ |\langle \B u_\delta\rangle (\Omega\setminus\overline{\Omega}_0) &\leq k_1\int_{\Omega\setminus \overline{\Omega}_0} (1+ |\B u_\delta|)dx.\notag
                \end{align}
        Due to \eqref{estProp2.11}, we know that 
               \begin{align*} 
                   \left||\B u|(\Omega) -\int_{\Omega} |\B u_\delta| dx   \right|\leq \delta.
               \end{align*}
        This along with \eqref{eq:estsatisfiedbyepsilon1} and \eqref{eq:intomega-omega0bu} leads us to deduce that 
              \begin{align}
                  \int_{\Omega\setminus\overline{\Omega}_0} |\B u_\delta| dx = \int_{\Omega} |\B u_\delta| dx - \int_{\Omega_0} |\rho_{\e_1}\ast(\varphi_1(\B u))| dx
                  \leq |\B u| (\Omega) - |\varphi_1 \B u| (\Omega)+2 \delta
                  = \int_{\Omega\setminus \Omega_0} |\B u| +2\delta\leq 3\delta.\notag
              \end{align}
        Finally, 
            \begin{equation}
                \notag
                \left| \langle \B u \rangle (\Omega)- \langle\B u_\delta\rangle (\Omega)  \right| \leq C \delta,
            \end{equation}
        concluding the proof. 
\end{proof}

To conclude this section, we prove, in the spirit of \cite[Theorem 2.4]{T83}, a compactness result for the space $BV^\B$, with $\B$ being of first order. 
\begin{Theorem}
\label{thm:cptembedding}
    Let $\Omega$ be an open and bounded subset of $\R^N$ with Lipschitz boundary. Let $\B$ be  a $\mathbb{C}$-elliptic operator as in \eqref{Budef} with $k=1$. The injection $BV^\B(\Omega)$ into $L^p(\Omega; V)$ is compact for any $p\in [1, {N\over N-1})$.
\end{Theorem}
\begin{proof}
    Let $(u_j)_j$ be a bounded sequence in $BV^\B(\Omega)$. Extending by zero to $\R^N$ and denoting the new sequence by $(\widetilde{u}_j)_j$, we know that $(\widetilde{u}_j)_j\subset BV^\B(\R^N)$ is still bounded (see \cite[Corollary 4.21]{BDG20}).  Let $\rho\in C_c^\infty(\R^N)$ be such that $\rho\geq 0$ and $\int_{\R^N} \rho(x)dx =1$. Set $\rho_n(x)= n^{-N}\rho(x/n)$. Then,  $\rho_n\ast\widetilde{u}_j\in W^{\B, 1}(\Omega)$ and is bounded independently of $j$ and $n$. Since $\B$ is a $\mathbb{C}$-elliptic operator, the Sobolev embedding theorem (see \cite[Theorem 1.1]{GR19}) guarantees that $\rho_n\ast\widetilde{u}_j$ is relatively compact in $L^p(\Omega; V)$ for any $p\in[1, {N\over N-1})$. Here, note that the support of $\rho_n\ast\widetilde{u}_j$ is a fixed compact set. Moreover, since the injection $W^{\B, 1}(\Omega)$ into $L^{{N\over N-1}}(\Omega)$ is continuous (cf. \cite[Theorem 1.1]{GR19}), we have that $\rho_n\ast\widetilde{u}_j\in L^{{N\over N-1}} (\R^N)$ and $\rho_n\ast\widetilde{u}_j \to u_j$ in $L^{{N\over N-1}}$, as $n \to +\infty$. This along the relatively compactness of $\rho_n\ast \widetilde{u}_j$ in $L^p$ leads us to conclude the result for $u_j$ after extracting a diagonal subsequence.
\end{proof}

\begin{Corollary}\label{embeddingk}
Let $\Omega$ be an open and bounded subset of $\R^N$ with Lipschitz boundary. Let $\B$ be  a $k$- homogeneous, $\mathbb{C}$-elliptic operator as in \eqref{Budef} with $k>1$. The injection of $BV^\B(\Omega)$ into $W^{k-1,p}(\Omega; V)$ is compact for any $p\in [1, {N\over N-1})$.
\end{Corollary}
\begin{proof}[Proof]
By \cite[Theorem 3.7]{DG24} (see also Theorem \ref{thm2.3}), it follows that bounded sequences in $BV^{\mathcal B}(\Omega)$ are bounded in $W^{k-1,1}(\Omega;W)$. Hence, classical Sobolev's embedding theorems guarantee compactness in $W^{k-2,p}(\Omega,W)$, for $p \in [1,\frac{N-1}{N})$ and a bound on the $(k-1)$-th derivatives.
To conclude, it suffices to invoke \cite[Theorem 20]{ARS25}, which ensures that every linear operator $\mathcal B$ of order $k$, with $k>1$, can be seen as a first order linear $\mathbb C$-elliptic operator $\mathcal C$, acting on the space of $(k-1)$-th distributional derivatives. Consequently, in view of Theorem \ref{thm:cptembedding}, the claim is proved. \end{proof}

\subsection{Integrands}\label{subsecint}

The following definitions and properties are recalled for the readers' convenience and we refer to \cite{AR21} for more details. 
Let $\Omega$ be an open subset of $\R^N$.  

\begin{Definition}
    \label{recfuncts}
    For a Borel integrand $f:\Omega \times W \to [0,\infty)$ for which 
there exists a modulus of continuity $\omega : [0,\infty) \to [0,\infty),$  
 such that 
 \begin{align}\label{modcont}
 |f (x, \xi) - f (y, \xi)| \leq \omega(|x - y|)(1 + |\xi|)\;\;\; \hbox{ for all } x, y \in \Omega \hbox{ and }\xi\in W,
 \end{align} 
 we define
    \begin{itemize}
    \item[(i)] the {\it generalized} lower (upper) recession function $f_{\#}$,  ($f^{\#}$) as
    \begin{align}\label{fdown}
    f_{\#}(x, \xi):=\liminf_{\substack{x' \to x,\\
    t\to +\infty,\\
    \xi' \to \xi}} \frac{f(x', t \xi')}{t},  \;\;\;\;\;  \;\;\; \left(f^{\#}(x, \xi):=\limsup_{\substack{x' \to x,\\
    t\to +\infty,\\
    \xi' \to \xi}} \frac{f(x', t \xi')}{t}\right).
    \end{align}
     \item[(ii)] the strong recession function $f^\infty(x, \xi)$ as
    \begin{align}\label{finfty}
    f^\infty(x, \xi):=\lim_{\substack{x' \to x,\\
    t\to +\infty,\\
    \xi' \to \xi} } \frac{f(x', t \xi')}{t},
    \end{align}
    whenever this limit exists.
\end{itemize}\end{Definition}

For $f \in C(\Omega \times W)$, we recall the transformation
$(S f )(x,\hat \xi) := (1 -|\hat \xi|) f
\left(
x,\frac{\hat \xi}
{1 - |\hat \xi|}\right) $ for every  $(x,\hat \xi) \in  \Omega \times B_W$,
where $B_W$ denotes the open unit ball in $W$. Note that $S f \in C(\Omega \times B_W)$. 

Recall now the sets $E(\Omega; W) := \{f \in C(\Omega \times W : S f \hbox{ extends to }C(\overline{\Omega \times B_W})\}$ and  $E(W) :=\{f \in C(W) : S f \hbox{ extends to } C(\overline{B_W})\}$ which have been introduced in \cite{KR10}. Heuristically, $E(W)$ is isomorphic to the continuous functions on the compactification
of $W$ that adheres to it each direction at infinity. In particular, all functions $f\in E(\Omega; W)$ have a linear growth at infinity, i.e., there exists a positive constant $M$ such that
$| f (x, \xi)| \leq M(1 + |\xi|)$ for all $x \in \Omega$ and all $\xi \in W$. Endowed with the norm $ \|f \|_{E(\Omega;W)} := \|S f \|_{\infty}$, $E(\Omega; W)$ is a Banach space and $S$ is an isometry with inverse $(Tg)(x, \xi) := (1 + |\xi|)g\left(x,\frac{\xi}{1 + |\xi|}\right),$ $(x, \xi) \in \Omega \times W.$ For any  $f \in E(\Omega; W)$, the strong recession function $f^\infty(x, \xi)$ defined as in \eqref{finfty}, exists. Moreover,  any $f \in E(\Omega; W)$ satisfies $ f (x, \xi) = (1 + |\xi|)S f\left(x,\frac{\xi}{1 + |\xi|}\right),$ for all $x\in\Omega, \xi\in W.$ In particular, there exists a modulus of continuity $\omega : [0,\infty) \to [0,\infty),$ depending solely on the uniform continuity of $S f$, such that \eqref{modcont} holds.

\subsection{$\mathcal A$ ($\mathcal B$)- quasi-convexification}\label{subsecABquasiconvex}

We are going to recall a crucial result proved in \cite[Corollary 1]{R19}.

\begin{Proposition}[{\cite[Corollary 1]{R19}}]
    \label{cor1R19}
Let $\mathcal A$ be a linear homogeneous differential operator on $\R^N$ from $W$ to $X$ satisfying the constant rank property. Let $\B$ be the operator given by Proposition \ref{prop:ExistenceofAandB} (b) 
and let $f : W \to \mathbb R$ be Borel measurable and locally bounded. Then
\begin{equation}\label{QAf}
{\mathcal Q}_{\mathcal A} f (\eta):= 
\inf\left\{
\int_{\mathbb T_n}f (\eta + w(x)) dx : w \in C^\infty(\mathbb T_N; W): \mathcal A w=0, \int_{\mathbb T_n}w(x)dx=0,\right\}
\end{equation}
and
\begin{equation}\notag
{\mathcal Q}_{\mathcal B} f (\eta):=
\inf\left\{
\int_{[0,1]^N}f (\eta + \mathcal B u(x)) dx : u \in C^\infty_C((0,1)^N; V)\right\}.
\end{equation}
are equal for all $\eta \in W$.
\end{Proposition}

\begin{Remark}
\label{remAquasi}
\begin{itemize}
\item[i)] By \cite[Lemma 3]{R19} the $\mathcal A$ (or equivalently $\mathcal B$)-quasiconvexity is independent on the domain of integration.
\item[ii)] By \cite[Theorem 3.6]{FM99}, \cite[Remark 1.3]{APR20} and \cite[Theorem A.1]{AR21} (see also \cite[Theorem 1.2]{APR20}), the function ${\mathcal Q}_{\mathcal A}f$ is $\mathcal A$-quasiconvex, namely
\begin{align}\label{Aqcxdef}
{\mathcal Q}_{\mathcal A} f (\eta) &=\inf\left\{
\int_{\mathbb T_n}{\mathcal Q}_{\mathcal A}f (\eta + w(x)) dx : w \in C^\infty(\mathbb T_N; W): \mathcal A w=0, \int_{\mathbb T_n}w(x)dx=0,\right\}\nonumber\\
&= \inf\left\{
\int_{[0,1]^N}\mathcal Q_{\mathcal B}f (\eta + \mathcal B u(x)) dx : u \in C^\infty_C((0,1)^N; V)\right\}.
\end{align}
\item[iii)] From \eqref{QAf}, it trivially follows that
\begin{equation}\label{QAfleqf}\mathcal Q_{\mathcal A}f(\eta)\leq f(\eta), \hbox{ for every }\eta \in W.
\end{equation}
\item[iv)] 
For a upper semicontinuous and locally bounded $f:W \to \mathbb R$ and for constant rank operators $\mathcal A$, $\mathcal Q_{\mathcal A}f$, defined as in \eqref{Aqcxdef}, is the greatest $\mathcal A$-quasiconvex function, such that \eqref{QAfleqf} holds, namely the {\it $\mathcal A$-quasiconvex envelope} of $f$. Indeed, this is proven in \cite[Proposition 3.4]{FM99} when $\mathcal A$ is a first-order operator, but a careful inspection of that proof guarantees that the same arguments could be used for operators of any order, indeed see \cite[Lemma 2.18]{APR20}.

\item[v)] {\rm \cite[Theorems 1.1 and 3.6]{BFL00}} are still valid for higher-order operators $\A$, namely the superlinear case counterparts of {\rm \cite[Theorem 1.2 or 1.10]{APR20}} or {\rm \cite[Theorem A.1]{AR20}} holds, also when $\mathcal A$ is an operator of order higher than one, see Section \ref{relBFLhigh} for details.
\end{itemize}
\end{Remark}

If $f:W \to \mathbb R$ has linear growth from above and below, then the strong recession function of ${\mathcal Q}_{\mathcal B} f$ is defined by \eqref{finfty}, i.e. for every $\xi \in \mathcal R(\mathcal B):={\rm span}(\mathcal I_\mathcal B)$, (with $\mathcal I_{\mathcal B}$ {\it image cone} of $\mathbb B$; Defined in \eqref{imagecone})
\begin{align}\label{strec}
({\mathcal Q}_{\mathcal B} f)^\infty(\xi): =\lim_{\substack{t \to +\infty\\\xi' \to \xi}} \frac{{\mathcal Q}_{\mathcal B} f(t \xi')}{t},
\end{align}
whenever the limit (which is taken in $\mathcal R(\mathcal B)$) exists. 

Observe that the same arguments of \cite[Lemma A.1]{BDG20} (despite the fact that here $\mathcal B$ is a $k$ homogeneous  partial differential operator), guarantee that the limit in \eqref{strec} exists for every $\xi \in {\mathcal I}_{\mathbb B}$ and
\begin{align}\label{exiQBfinfty}
({\mathcal Q}_{\mathcal B} f)^\infty(\xi)=\lim_{\substack{t \to +\infty\\\xi' \to \xi}} \frac{{\mathcal Q}_{\mathcal B} f(t \xi')}{t}
= \lim_{\substack{t \to +\infty}} \frac{{\mathcal Q}_{\mathcal B} f(t \xi)}{t}=\sup_{t \geq 0}\frac{{\mathcal Q}_{\mathcal B} f(t \xi)}{t}.
\end{align} 
Indeed the key point in the \cite[Lemma A.1]{BDG20} is the fact that $\mathcal Q_{\mathcal B} f$, being $\mathcal A$-quasiconvex as observed in Remark \ref{remAquasi}-(iv), 
 is also rank-one convex in the directions of the $\Lambda_\mathcal A$ cone (see \cite[Proposition 3.4]{FM99}, where the result is stated for first-order operator but the proof does not use the order of $\mathcal A$, and \cite{KK16}).
 This latter convexity property, together with the linear growth of $f$, inherited by $\mathcal Q_{\mathcal B}f$, via \eqref{QAfleqf}, guarantees also that $\mathcal Q_{\mathcal B}f$ is Lipschitz (see \cite{Da08}).

\section{Relaxation}\label{secRel}

In this section, we show a relaxation result for linear-growth integral functionals defined on $W^{\B, 1}(\Omega)$. Throughout this section, we assume that $\B$ is a $k$-th order linear differential, $\mathbb C$-elliptic operator as in \eqref{Budef}. Note that, by definition, $\B$ satisfies the constant rank properties, too.

Let $\Omega$ be a bounded open set of $\R^N$ with Lipschitz boundary. Let $f:\Omega\times W\to [0, \infty)$ be a continuous integrand satisfying 
Assumptions $ (H_{1 \B})$-$(H_{4 \B})$ of Theorem \ref{thm:RelaxationB}. Let $G_{\mathcal B}: W^{\B, 1}(\Omega)\to [0, \infty)$ be a linear-growth integral functional of the form \eqref{GB} with $g_2$ therein replaced by the continuous integrand $f:\Omega\times W\to [0, \infty)$ satisfying Assumptions $(H_{1 \B})$-$(H_{3 \B})$. With an abuse of notation, for any $u\in W^{\B, 1}(\Omega)$, the measure $\B u \ll\mathcal L^N$ is identified with its density $\frac{d \B u}{d \Leb^N}$.
We are going to provide an integral representation of the weakly$^*$ lower semicontinuous envelope $\overline{\mathcal{G}}_{\B}$ of $G_\B$ given by \eqref{GBcal}.

\begin{Remark}\label{equivembed}
By Corollary \ref{embeddingk}, $\overline{{\mathcal G}}_{\B}$ in \eqref{GBcal} is equivalent to
\begin{align}
        \nonumber
        \mathcal{G}_{\B}(u):= \inf \biggl\{ \liminf_{j\to\infty} G_{\B}(u_j) :(u_j)\subset W^{1,\mathcal B}(\Omega),\hspace{0.2cm}  u_j\to u \hbox{ in } L^1(\Omega; V),
        \B u_j\mathcal{L}^n \weakcon \B u \hspace{0.3cm} \mbox{in } \Mcal(\Omega;W)    \biggr\}.
    \end{align}
\end{Remark}
We are in position of proving Theorem \ref{thm:RelaxationB}. For the readers' convenience, we briefly restate it. 

\begin{Theorem}
\label{thm:RelaxationBbis} Let $\B$ be a $k$-th order, $\mathbb C$-elliptic, linear partial differential operator as in \eqref{Budef}.
    Let $f:\Omega\times W\to[0,\infty)$ be a continous integrand satisfying Assumptions $(H_{1 \B})$-$(H_{4 \B})$. Then, $\overline{\mathcal{G}}_{\B}$ defined in \eqref{GBcal} is characterized by 
    \begin{equation}
              \notag
    \overline{\mathcal{G}}_{\B}(u) = \int_\Omega Q_{\B}f\left( x, {d\B u\over d\Leb^N}(x) \right)dx + \int_\Omega (Q_{\B}f)^\#\left( x, {d\B^s u\over d|\B^s u|}(x)\right)d|\B^s u|(x),
          \end{equation}
\end{Theorem}

Before proving Theorem \ref{thm:RelaxationB}, we recall a useful result which is crucial for our proof.

\begin{Theorem}[{\cite[Theorem 1.8]{APR20}}]
\label{Theorem1.8AAPR20}
Let $\A$ be a $l$-th order linear partial differential operator. Let $\Omega\subset\R^N$ be an open and bounded subset with Lipschitz boundary. Let $f: \Omega\times W\to [0, \infty)$ be a continuous integrand such that assumptions $(H_{1 \B})$-$(H_{3 \B})$  and that $(H_{4 \B})$ holds for any $(x, \xi)\in\Omega\times {\rm span}\Lambda_{\A}$, $\Lambda_\A$ given by \eqref{Lambdacone}. Further, suppose that for any $\mu\in \Mcal(\Omega; W)\cap \ker\A$ there exists a sequence $(v_j)_j\subset L^1(\Omega; W)\cap\ker\A$ such that 
\begin{equation}
    \label{eq1.5APR20}
    v_j\Lcal^N\weaklystar\mu \quad \mbox{in } \Mcal(\Omega; W) \quad\mbox{and}\quad \langle v_j\Lcal^N\rangle(\Omega)\to \langle\mu\rangle(\Omega),
\end{equation}
where $\langle\cdot\rangle$ is the area functional in \eqref{areafunctional}. 
Let $\mathcal G$ be the functional defined as 
       \begin{equation}
           \notag
           \mathcal G_\A(v):= \int_\Omega f(x, v(x))dx, \qquad v\in L^1(\Omega; W)\cap \ker\A,
       \end{equation}
and let $\overline {\mathcal G}_\A$ be its  weakly$^\ast$ lower semicontinuous envelope, i.e.
\begin{equation}
    \notag
    \overline{\mathcal G}_\A(v):= \inf\left\{ \liminf_{j\to\infty} \mathcal G_\A(v_j) : (v_j)\subset L^1(\Omega; W)\cap\ker\A, v_j\Lcal^N\weaklystar \mu \quad \mbox{in } \Mcal(\Omega; W)     \right\}.
\end{equation}
Then,
    \begin{equation}
        \notag
        \overline{\mathcal G}_\A(\mu)= \int_\Omega Q_{\A}f\left(x, {d\mu\over d\Lcal^N}(x)\right)dx + \int_\Omega (Q_\A f)^\#\left(x, {d\mu^s\over d|\mu^s|}(x)\right)d|\mu^s|(x).
    \end{equation}
\end{Theorem}
\smallskip

\begin{proof}[Proof of Theorem \ref{thm:RelaxationB}]
    We split the proof into a lower and an upper bound. \\
    {\bf\underline{Lower bound.}} First, recall that $\mathcal B$ is $\mathbb C$-elliptic. Hence, Proposition \ref{prop:ExistenceofAandB} (a) guarantees that there exists an $l$-th order linear differential operator $\mathcal A$ such that the complex symbol involving $\mathbb A$ and $\mathbb B$ is exact. Then, noting that \cite[Theorem 1.3]{AR21} guarantees the validity of \eqref{eq1.5APR20}, we can make use of Theorem \ref{Theorem1.8AAPR20} with $\mu=\B u$. Therefore, 
                            \begin{align}
                                 \overline{\mathcal G}_\B(u)&\geq \overline{\mathcal G}_\A(\B u) \notag\\
                                 &=\int_{\Omega} Q_{\A}f\left(x, {d\B u\over d\Leb^N}(x) \right)dx +\int_{\Omega}(Q_{\A}f)^\#\left( x, {d\B^s u\over d|\B^s u|}(x)\right)d|\B^s u|(x)\notag\\
            &=\int_{\Omega} Q_{\B}f\left(x, {d\B u\over d\Leb^N}(x) \right)dx +\int_{\Omega}(Q_{\B}f)^\#\left( x, {d\B^s u\over d|\B^s u|}(x)\right)d|\B^s u|(x),\notag
                            \end{align}
        where in the last equality, we have used that $Q_{\A}f=Q_{\B}f$ as stated in Proposition \ref{cor1R19}.  This concludes the proof of the lower bound.
        \smallskip

    {\bf\underline{Upper bound.}}  The proof of the upper bound relies on a density argument, making use  of Theorem \ref{thm:densitynormandarea} and Proposition \ref{prop:densityuptoboundary}. We split our arguments  into several steps. 
    \medskip 
    
    \underline{Step 1.} Consider $u\in C^\infty(\Omega; W)\cap W^{\B, 1}(\Omega)$. In this case, the construction of the recovery sequence is strongly based on the arguments carried out in Step 3 of \cite[Theorem 1.8]{APR20} (see also \cite[Step 3 of Theorem 1.7 and Remark 5.1]{APR20}); Due to the fact that $\B u\in (C^\infty\cap L^{1})(\Omega, W)$ along with $\A(\B u)=0$.\\
    For the readers' convenience, we recall the main steps.  
    Fix $m \in \mathbb N$ and consider a partition of $\mathbb R^N$ of cubes of side length $m$. Let $\{Q_i^m\}_{i=1}^{L(m)}$
be the maximal collection of those cubes whose centers $\{x_i^m\}_{i=1}^{L(m)}$ are compactly contained in $\Omega$. It results that
$\mathcal L^N(\Omega) =\sum_{i=1}^{L(m)}
\mathcal L^N(Q_i^m) + o_m(1)$
with $o_m(1)\to 0$  as $m \to +\infty$. 
Without loss of generality, in view of rescaling arguments, we can further assume that $\Omega \subset \subset Q$, where $Q$ is the unit cube.
Following the same arguments as in \cite[Proof of Theorem 1.7]{APR20}, using the cubes $Q_i^m$, for every $m$, it is possible to build a sequence $(V^m_j)_j\subset C^\infty_{\rm per}(Q;W)\cap\ker \A$ such that $V^m_j\weakcon 0$ in $\Mcal(Q; W)$, as $j\to +\infty$. 
 Note that in view of \cite[Lemma 2.15]{APR20}, we can also assume that $\int_Q V^m_j dx=0$. For a fixed $m, j \in \Nb$, we can apply \cite[Lemma 2]{R19} to $V^m_j$ to deduce the existence of $(U^m_j)\in C^\infty_{\rm per}( Q; W)$ such that $V_j^m = \B U_j^m$. Hence, we take the sequence  $u+ U_j^m$, observing that $\B u + V_j^m = \B u + \B U_j^m = \B (u+ U_j^m)$ and arguing as in Step 3 of \cite[Theorem 1.7]{APR20}, we build a suitable recovery sequence and we deduce
          \begin{align*}
             \overline{\mathcal{G}}_{\B}(u) \leq \inf_{m>0}\liminf_{j\to\infty} \int_\Omega f(x, \B(u+ U_j^m)(x))dx 
              \leq \int_\Omega {\mathcal Q}_{\mathcal A} f(x, \B u(x)) dx = \int_\Omega {\mathcal Q}_{\B} f (x, \B u(x))dx, 
          \end{align*}
    where in the last equality we have used  Proposition \ref{cor1R19}. 
    To conclude the proof, we need to prove that $u+U_j^m\to u$ in $W^{k-1, 1}(\Omega; V)$. To this end, it suffices to show that for any $m$, $U_j^m\to 0$ in $W^{k-1, 1}(\Omega; V)$. Without loss of generality, up to subtracting its average, for every $m$ the sequence $U^m_j \rightharpoonup 0$ in $L^{\frac{N}{N-1}}(\Omega;V)$, as $j \to +\infty$ and for every $m$, hence the convergence holds also in $L^1(\Omega;V)$, as $j\to +\infty$ for every $m$. Thus, by the bounds on $|\mathcal B U^m_j|(\Omega)$ and the convergence in $L^1$, it suffices to invoke Corollary \ref{embeddingk} to conclude the proof.
    
 \underline{Step 2.} Let $u\in W^{\B, 1}(\Omega)$. Thanks to Theorem \ref{thm:densitynormandarea}[(ii)], there exists a sequence $(u_j)_j\in C^\infty(\Omega; W)\cap W^{\B, 1}(\Omega)$ such that $u_j$ strongly converges to $u$ in $W^{\B,1}(\Omega)$. In particular, $u_j\weakcon u$ in $W^{\B, 1}(\Omega)$ and, by Proposition \ref{propWB=VB}, the convergence holds also in $W^{k-1,1}(\Omega;V)$. By the lower semicontinuity of $\overline{\mathcal{G}}_{\B}$ together with Step 1, we deduce 
              \begin{align*}
            \overline{\mathcal{G}}_{\B}(u)&\leq \liminf_{j\to\infty} \overline{\mathcal{G}}_{\B}(u_j) 
                  \leq \liminf_{j\to\infty} \int_\Omega {\mathcal Q}_{\B} f (x, \B u_j(x))dx.
              \end{align*}
    Since $\Q_{\B}f(x,\cdot)$ is a Lipschitz continuous function in $\mathcal R(\B)$ (see \cite[Proposition 4.5]{AR21}), we immediately deduce via Dominated Convergence theorem that
       \begin{align}
           \overline{\mathcal{G}}_{\B}(u)&\leq \limsup_{j\to\infty} \int_\Omega {\mathcal Q}_{\B} f (x, \B u_j(x))dx 
            \leq \int_\Omega {\mathcal Q}_{\B} f (x, \B u(x))dx,\label{step3relformula}
       \end{align}
     concluding the proof of Step 2.
  
\underline{Step 3.} Let $u\in BV^\B(\Omega)$. Proposition \ref{prop:densityareastrictly} guarantees the existence of a sequence $(u_j)_j$ in $(C^\infty\cap V^{\B, 1})(\Omega)$ such that $u_j \to u$ in area-strict sense on $\Omega$. 
Invoking the lower semicontinuity of $\overline{\mathcal{G}}_\B$ and Step 2, it follows that 
     \begin{align*}
         \overline{\mathcal{G}}_\B(u) &\leq \liminf_{j\to\infty} \overline{\mathcal{G}}_\B(u_j)
          \leq\liminf_{j\to\infty} \int_\Omega \mathcal{Q}_{\B} f(x, \B u_j(x))dx.
     \end{align*}
It is worth noting that $\mathcal{Q}_{\B}f(x,\cdot)$ is $\B$-quasiconvex (see Remark \ref{remAquasi}, (iv), which relies on \cite{FM99} and \cite{R19}) and  satisfies Assumption $(H_{1\B})$ in $\Omega \times \mathcal R(\B)$, as remarked in \cite[Proposition 4.5]{AR21}. It is also immediately seen that $\mathcal{Q}_{\B}f(x,\cdot)$  satisfies $(H_{2 \B})$. Furthermore,  Assumption $(H_{3 \B})$ also holds on $\Omega \times \mathcal R(\B)$ as shown in \cite[formula (5.11)]{APR20} (see also \cite{BDG20}). In particular, putting together assumptions $(H_{1\B})$ and $(H_{3 \B})$, we get that $\mathcal{Q}_{\B}f$ is continuous on $\Omega \times \mathcal R(\B)$. 
Exploiting the continuity of $\mathcal{Q}_\B f$ along with the bound from below in Assumption $(H_{2\B})$,  we can apply \cite[Lemmma 2.3]{AB} to $\mathcal{Q}_{\B}f$, and we can find a sequence $(g_n)_n$ of integrands in $E(\Omega; \mathcal R(\B))$ such that
\begin{align}\nonumber\sup_n
g_n (x, \xi) = \mathcal{Q}_\B f(x, \xi) \hbox{ and  } \sup_n
(g_n)^\infty(x, \xi) = (\mathcal{Q}_\B f)_{\#}(x, \xi).
\end{align}
Therefore, the same arguments of \cite[Proposition 5.1]{BDG20}, relying on \cite[Lemma 2.3]{AB}, (see also \cite[proof of Theorem 4]{KR10}) 
allow us to obtain that 
\begin{align}\label{1sharp}
\liminf_j \overline{\mathcal G}_{\mathcal B}(u_j) &=\liminf_j \int_\Omega {\mathcal Q}_{\mathcal B} f(x,\mathcal B u_j(x))dx 
\\
&\geq \int_\Omega g_n (x, \mathcal B u(x))dx + \int_\Omega g_n^\infty\left(x,\frac{d \mathcal B u}{d |\mathcal B^s u|}(x)\right) d |\mathcal B^s u|(x)\notag\\
&\geq
\int_\Omega {\mathcal Q}_\B f(x,\B u(x))dx +\int_\Omega ({\mathcal Q}_\B f)_{\#}\left(x,\frac{d \mathcal B u}{d |\mathcal B^s u|}(x)\right) d |\mathcal B^s u|(x)\notag.
\end{align}
along the sequence $(u_j)_j$ converging strictly in area to $u$. Applying the same argument to $-\mathcal{Q}_{\B}f$, which also satisfies $(H_{1 \B})$-$(H_{3 \B})$ in $\Omega \times \mathcal R (\B)$ and taking into account the same arguments as in the proof of \cite[Theorem 4]{KR10}, from \eqref{1sharp} applied to $-\mathcal{Q}_{\B }f$, we get
\begin{align}\label{2sharp}
\limsup_j \overline{\mathcal G}_{\B}( u_j) &=\limsup_j \int_\Omega \mathcal{Q}_\B f(x,\B u_j(x))dx 
\\
&\leq \int_\Omega \mathcal{Q}_\B f(x,\B u(x))dx +\int_\Omega (\mathcal{Q}_\B f)^{\#}\left(x,\frac{d \B u}{d |\B^s u|}(x)\right) d |\B^s u|(x)\notag,
\end{align}
where we have exploited the following inequality 
     \begin{equation}
         (-\mathcal{Q}_{\B}f)_{\#}(x, \xi):=\liminf_{\substack{x' \to x,\\
    t\to +\infty,\\
    \xi' \to \xi}} \frac{- \mathcal{Q}_{\B}f(x', t \xi')}{t}=
    -\limsup_{\substack{x' \to x,\\
    t\to +\infty,\\
    \xi' \to \xi}} \frac{ \mathcal{Q}_{\B}f(x', t \xi')}{t}
=-(\mathcal{Q}_{\B}f)^{\#}(x, \xi). \label{limitsuplowrec}
     \end{equation}
Combining \eqref{1sharp} and \eqref{2sharp}, it follows that
$$
\lim_{j\to\infty} \int_\Omega \mathcal{Q}_{\B}f(x,\B u_j(x))dx = \int_\Omega \mathcal{Q}_{\B }f(x,\B u(x))dx + \int_\Omega (\mathcal{Q}_{\B}f)^{\#}\left(x,\frac{d \B u}{d |\B^s u|}(x)\right)d |\B^s u|(x),$$
whenever $u_j \to u$ area strictly in $BV^\B(\Omega)$. 
This along with \eqref{step3relformula} concludes the proof.
\end{proof}

\begin{Remark}\label{remOmegachiuso}
\begin{itemize}
     \item[\rm (i)] The properties of $\Q_{\B}$ follow from those of $\Q_{\A}$. For instance the continuity is a consequence of {\rm \cite[Proposition 4.5]{AR21}} (see also {\rm \cite[Remark A.1]{AR21}}).
    \item [\rm (ii)] As pointed out in {\rm \cite{APR20}}, in view of the fact that $\mathcal{Q}_\B f$ satisfies Assumptions $(H_{1 \B})$-$(H_{3 \B})$ (see for instance {\rm \cite{AR21}}) the limits in \eqref{limitsuplowrec} can be taken just along $t \to +\infty$. Hence, 
 $$ (-\mathcal{Q}_{\B}f)_{\#}(x, \xi)=\liminf_{t\to +\infty} \frac{- \mathcal{Q}_{\B}f(x, t \xi')}{t}=
    -\limsup_{
    t\to +\infty} \frac{ \mathcal{Q}_{\B}f(x, t \xi')}{t}
=-(\mathcal{Q}_{\B}f)^{\#}(x, \xi).$$

    \item[\rm (iii)] One may replace $(\mathcal{Q}_{\B}f)^\#$ with $(\mathcal{Q}_{\B}f)^\infty$ in the proof of lower bound provided that 
    $f:\overline \Omega \times W \to[0,\infty)$ is a continuous integrand such that $(H_{2 \B})$ holds and
               \begin{itemize}  
                   \item[ $(H^{'}_{1 \B})$] $f:\overline \Omega \times \mathcal R(B) \to [0, \infty)$ is Lipschitz in the second variable uniformly w.r.t. the first one;
                    \item[{ $(H_{3 \B}^{'})$}]  Assumption $(H_{3 \B})$ holds for $\mathcal{Q}_{\B}f$  restricted to  $\Omega \times \mathcal R(\B)$.                    
               \end{itemize}
Indeed, first recall that $\mathcal{Q}_{\A}f=\mathcal{Q}_{\B}f$ and $\mathcal{Q}_{\B}f(x,\cdot)$ is $\B$-quasiconvex (see, e.g., {\rm \cite{FM99}} and {\rm \cite{R19}}). In addition, $\mathcal{Q}_{\B}f$ satisfies Assumptions $(H_{1 \B})$-$(H_{2 \B})$ on $\overline \Omega \times \mathcal R(\B)$ (see for instance {\rm \cite{BDG20}} and  {\rm \cite[formula (5.11)]{APR20}}), as well as assumptions $(H_{3 \B}^{'})$ on $\overline \Omega \times \mathcal R(\B)$ (see, e.g., {\rm \cite{FM99}} and {\rm \cite[Lemma 6.1]{BDG20}}, in which proof the order of the operator does not play any role). Thus, {\rm \cite[Proposition 5.1]{BDG20}} may be applied to prove the lower semicontinuity in $BV^\B(\Omega)$ of $\int_\Omega \mathcal{Q}_\B f(x,\frac{d\B u}{d \mathcal L^N})dx + \int_\Omega (\mathcal{Q}_\B f)^\infty \left(x, \frac{d \B^s u}{d |\B^s u|}\right)d |\B^su|$. Observe that the argument in {\rm \cite{BDG20}} uses the area strict density of $W^{\mathcal B,1}(\Omega)$ in $BV^{\mathcal B}(\Omega)$, which by Proposition \ref{prop:densityareastrictly} is known to be true for operators $\mathcal B$ of any order.

\item[\rm (iv)] Under assumptions $ (H_{1 \B}^{'}), (H_{2 \B})$ and  $(H_{3 \B}^{'})$ via a careful inspection of the proof of {\rm \cite[Theorem 5.1]{BDG20}} we can conclude that the proof of Step 3 in the upper bound of Theorem \ref{thm:RelaxationB} can be obtained by a different argument, once again relying on the density of $W^{\mathcal B, 1}(\Omega)$ in $BV^{\mathcal B}(\Omega)$, obtained in Proposition \ref{prop:densityareastrictly} . Indeed, if $f:\overline{\Omega}\times W\to [0, \infty)$ is a continuous integrand satisfying assumptions $(H_{1\B}^{'}), (H_{2 \B}), (H_{3 \B}^{'})$, {\rm \cite[Proposition 5.1]{BDG20}} guarantees that $\int_\Omega \mathcal{Q}_\B f(x,\frac{d\B u}{d \mathcal L^N})dx + \int_\Omega (\mathcal{Q}_\B f)^{\#} \left(x, \frac{d \B^s u}{d |\B^s u|}\right)d |\B^su|$ is the area strict continuous extension of $\int_\Omega Q_{\B}f(x,\B u)dx$ from $W^{\B,1}(\Omega)$ to $BV^{\B}(\Omega)$.  In particular, the latter functional is area strictly continuous and moreover $(\mathcal{Q}_\B f)^{\#}= (\mathcal{Q}_\B f)^{\infty}$. 
\end{itemize}

\end{Remark}

\section{Interaction with measures in the relaxation process}
\label{Sect:Applications}

In this section, we will apply Theorem \ref{thm:RelaxationB} to extend \cite[Theorems 1.1 and 1.2]{KKZ23} in particular to the context of operators $\mathcal B$ coinciding either with $\mathcal E$ or $D^2$. More generally, the result we are going to prove can be treated within the framework of the Global Method for Relaxation introduced in \cite{BFMglob} and later generalized to many contexts (see, e.g., \cite{BMZNoDEA, BMZGaeta, CFVG20, FHP, H} among a much wider literature).  Thus, we introduce the following definition.
\begin{Definition}\label{GM}
A set function $\mathcal F:BV^{\mathcal B}(\Omega) \times L^1(\Omega;\mathbb R^l)\times \mathcal A(\Omega)\to [0,+\infty)$, which is the trace of a Radon measure absolutely continuous with respect to $|\mathcal B u|+ \mathcal L^N+ |v|$ belongs to the class
${\mathcal G \mathcal M}(\Omega)$, if for any positive Radon measure $\mu$,
it results that $\frac{d \mathcal F}{d \mu}(u(x_0),v(x_0))= \lim_{\varepsilon \to 0} \frac{m(u,v, Q(x_0,\varepsilon))}{\mu(Q(x_0,\varepsilon))}$, for $\mu$ -a.e. $x_0 \in \Omega$
where $Q(x_0,\varepsilon)\subset \mathbb R^N$ denotes a cube centered at $x_0$ with side-length $2 \varepsilon$, and
\[
m(u,v, A):=\inf \left\{\mathcal F(w,\mathfrak V,A): {\rm supp}(u-w)\subset \subset A, \int_A v dx= \int_A \mathfrak V dx, (w,\mathfrak V) \in BV^\B(\Omega)\times L^1(\Omega)\right\}.
\]   
\end{Definition}
In particular, when $\mathcal B$ has at least order $3$, then 
\begin{align}\label{GBmathcalF}
\mathcal F[(u,v; A)]=\int_A \mathfrak{f}(x, u(x), \nabla u(x), \nabla^2 u(x),\dots, \nabla^{k-1}u(x),\B u(x),v(x))d x + \lambda (A),
\end{align}
where $\lambda$ is a Radon measure singular with respect to $\mathcal L^N$,
and 
\begin{align}\label{fgothic}
\mathfrak{f}(x_0, a, \xi,  H, \dots, B, G)= \lim_{\varepsilon \to 0}\frac{m(a + \xi(x-x_0)+ \frac{1}{2} H(x-x_0)\cdot (x-x_0)+ \dots+  B, G, Q(x_0, \varepsilon))}{\varepsilon^N}
\end{align}
where
$a \in \mathbb R^d$, $\xi \in \mathbb R^{d \times N}$, $H \in \mathbb R^{d\times d \times N}$, $G \in \mathbb R^l$ and $B$ coincides with the linear operator $B_\beta$, where $\beta$ is such that $|\beta|=k$ in Definition \eqref{Budef}. For homogeneous linear operators $\B$ of order $k<3$, some derivatives are not present in the formula \eqref{GBmathcalF}.

Let $\Omega$ be a bounded open set of $\R^N$ with Lipschitz boundary. Let $\B$ be a $k$-th order, $\mathbb{C}$-elliptic, linear partial differential operator given by \eqref{Budef}. For any $(u,v)\in W^{\B, 1}(\Omega)\times L^1(\Omega;\R^l)$, we consider the functional
		\begin{equation}
		\label{originalfunctional}
		 F(u,v):= G_1(u,v)=\int_\Omega f_1(u(x))f_2(v(x))dx+ \int_\Omega h(\B u(x))dx,
		\end{equation}
where  $f_1, f_2$ and $h$ are continuous functions satisfying Assumptions $\rm (H_1), \rm (H_2q), \rm(H_3p), \rm (H_4)$, with $p=q=1$. We are going to provide an integral representation of the sequentially lower semicontinuous envelope of $F$ with respect to weak$^\ast$ convergence:
		\begin{align}\label{Frelax}
		\begin{aligned}
		\overline{\mathcal F}(u,v):=
		\inf\left\{\,\liminf_{j\to +\infty}
		F(u_j,v_j)
		\,\left|\,
		\begin{array}{l}
		(u_j,v_j)\in W^{\B,1}(\Omega)\times L^1(\Omega;\R^l),\\
		(u_j, 
  v_j)\overset{\ast}{\rightharpoonup} (u,v)\hbox{ in } BV^{\B}(\Omega)\times  \mathcal M(\bar\Omega;\R^l),
		\end{array}
		\right.\right\}.
		\end{aligned}
		\end{align}
As in \cite{KKZ23}, we will distinguish the cases $\Omega \subset \mathbb R$ and $\Omega \subset \mathbb R^N$, with $N \geq 2$.

\subsection{Interactions in the relaxation process with emerging measures in dimension $N\geq 2$}
\label{sec:relNgeq2}

The following result summarizes several results proven in \cite{KKZ23} and \cite{BDG20}
\begin{Proposition}[{\cite[Proposition 5.1]{BDG20}} and {\cite[Proposition 2.1]{KKZ23}}]

    \label{prop:areastrictcont}
    Let $f_2$ and $h$ be Borel functions satisfying Assumptions ${\rm (H_2q)}, {\rm (H_3p)}$ for $p=q=1$, and ${\rm (H_4)}$. Then, the functionals 
    \begin{align*}
    &	v\mapsto \int_{\overline\Omega} df^{\ast\ast}_2(v)(x),\quad \Mcal(\overline\Omega;\R^l)\to \R,\\
    &	v\mapsto \int_{\Omega} df^{\ast\ast}_2(v)(x),\quad \Mcal(\Omega;\R^l)\to \R,~~\text{and}\\
    &   u\mapsto \int_{\Omega} d\Q_{\B}h(\B u)(x), \quad BV^{\B}(\Omega)\to \R
    \end{align*}
are sequentially continuous with respect to area-strict convergence in 
$\Mcal(\overline\Omega;\R^l)$, $\Mcal(\Omega;\R^l)$ and $BV^{\B}(\Omega)$, respectively,  where $df^{\ast\ast}_2(v)(x)$ and $d\Q_{\B}h(\B u)(x)$ are defined as in \eqref{def:df}.
\end{Proposition}
We refer to \cite[Subsection 2.3]{KKZ23} for comments regarding the recession function of $f_2^{**}$.
\begin{Remark}
    \label{remBorel}
    There is no continuity imposed on $f_2$ and $h$, since their envelopes $(f_2)^{**}$ and $Q_\B h$ are Lipschitz functions in view of their convexity and $\B$-quasiconvexity (see {\rm \cite{FM99}} and {\rm \cite{R19}}), together with  $(H_{2}q)$ ($q=1$) and $(H_{3}p)$ ($p=1$), (see {\rm \cite{Da08}}).
\end{Remark}
For any $A\in O_r(\overline\Omega)$, we introduce the localized energy $ \Fcal(\cdot,\cdot,A): W^{\B, 1}(\Omega\cap A)\times \mathcal{M}(A;\mathbb R^l)\to [0,+\infty]$; Defined as
			\begin{align}
			\begin{aligned}
			\mathcal F(u,v; A):=\left\{\!\!\!\!
			\begin{array}{ll}
			\int_{\Omega\cap A} \big( f_1(u(x))f_2(v(x))+ h(\B u(x))\big)\,dx, &\hbox{ if } (u,v)\in W^{\B,1}(\Omega\cap A)\times L^1(\Omega\cap A;\mathbb R^l),\\ 
			\\
			+\infty,  &\hbox{ otherwise.}
			\end{array}
			\right.
			\end{aligned}
			\label{Flocalized}
			\end{align}
           
		Its \emph{relaxation}, i.e., its lower semicontinuous envelope with respect to weak$^*$ convergence, is the functional
		$\overline{\mathcal F}(\cdot,\cdot,A): BV^{\B}(\Omega\cap A)\times \mathcal M(A;\mathbb R^l)\to [0,+\infty]$
		defined as	
		\begin{align}\label{FlocrelaxB}
		\begin{aligned}
		\overline{\mathcal F}(u,v; A):=
		\inf\left\{\,\liminf_{j\to +\infty}
		\mathcal F(u_j,v_j,A)
		\left|
		\begin{array}{l}
		(u_j,v_j)\in W^{\B,1}(\Omega\cap A)\times L^1(\Omega\cap A;\mathbb R^l),\\
		(u_j,v_j)\overset{\ast}{\rightharpoonup} (u,v)\hbox{ in } BV^{\B}(\Omega\cap A)\times \mathcal M(A;\mathbb R^l)
		\end{array}
		\right.\!\!\!\right\},
		\end{aligned}
		\end{align}

\begin{Proposition}\label{proprel**}
    Let $\B$ be a linear partial differential operator, $\mathbb C$-elliptic, of order $k$. Let $f_1: V \to \R$, $f_2:\mathbb R^l \to \mathbb R$, $h: W \to \mathbb R$ be continuous functions satisfying ${\rm (H_1)}, {\rm (H_2q)}, {\rm(H_3p)}$, with $p=q=1$ and ${\rm (H_4)}$. Let $A\in\mathcal{O}_r(\overline{\Omega})$. Let $\mathcal{F}(\cdot, \cdot, A)$ and its relaxation $\overline{\mathcal{F}}(\cdot, \cdot, A)$ be given by  \eqref{Flocalized} and \eqref{FlocrelaxB}. Define
        \begin{align}\label{F**locrelaxB}
	\begin{aligned}
{\mathcal F}^{\ast \ast}(u,v; A):=&\inf\Big\{\liminf_{j\to +\infty}\int_{\Omega\cap A}f_1(u_j(x))f_2^{\ast\ast}(v_j(x)) + \int_{\Omega\cap A}d \Q_{\B} h(\B u_j)(x):\\
&\quad
BV^{\B}(\Omega\cap A)\times L^1(\Omega\cap A;\mathbb R^l)\ni (u_j,v_j) \overset{\ast}{\rightharpoonup} (u,v)\hbox{ in } BV^{\B}(\Omega\cap A)\times \mathcal M(A;\mathbb R^l)
\Big\}.
\end{aligned}
	\end{align}
Then, for every $A\in \mathcal O_r(\overline \Omega)$, such that
 \begin{align}\label{lippartbound}
 \partial (\Omega \cap A) \hbox{ is Lipschitz,}
 \end{align}
 it holds $\overline{\mathcal F}(u,v; A)= \mathcal F^{**}(u,v; A)$ for every $(u,v)\in BV^\B(\Omega \cap A) \times \mathcal M(A;\mathbb R^l)$.
\end{Proposition}
The proof is provided in the Appendix, for the readers' convenience, since it is not entirely similar to \cite[proof of Proposition 2.7]{KKZ23}.

The next lemma states that, for every partial differential operator, $\mathbb C$-elliptic of order $k$, the localized relaxed functional \eqref{F**locrelaxB} is the restriction of a Radon measure on $\mathcal{O}_r(\overline{\Omega})$. 
The proof is omitted relying on a suitable variant of the so-called De Giorgi-Letta criterion, see \cite[Lemma 2.5]{ABF} (\cite[Lemma A.1]{KKZ23}). The only difference relying on standard estimates possibly involving higher order derivatives due to the fact that $\mathcal B$ might be an operator of order greater than one. 
\begin{Lemma}\label{FcalRadon}
    Let $\Omega$ be a bounded open set of $\R^N$ with Lipschitz boundary. Let $\B$ be a linear partial differential operator of order $k$.  Let $f_1$, $f_2$ and $h$ be continuoud functions satisfying Assumptions ${\rm (H_1)}, {\rm (H_2q)}$, and ${\rm (H_3p)}$, with $p=q=1$. For any $(u, v) \in BV^{\B}(\Omega)\times \Mcal(\overline{\Omega}; \R^l)$, the set function $\overline{\mathcal{F}}(u, v, \cdot)$ given in \eqref{FlocrelaxB} is the trace of a Radon measure absolutely continuous with respect to $\mathcal{L}^N+ |\B u| + |v|$.
\end{Lemma}

The main result of this subsection is Theorem \ref{thm:relk<N,N>=2}, which states an integral representation for any relaxed functional $\overline{\mathcal F}(u,v)$ given by \eqref{Frelax} belonging to the class $\mathcal G \mathcal M(\Omega)$. It is worth emphasizing that this technical assumption is necessary only to compute the upper bound and it is satisfied by some classical operators such as $\B= \mathcal E$ and $\B=D^2$. Indeed, it suffices to invoke the integral representation in \cite{FMZ}, in turn relying on \cite{CFVG20} and \cite{FHP}.
We restate Theorem \ref{thm:relk<N,N>=2}, for the readers' convenience.

   \begin{Theorem}
   Let $N\geq 2$ and $\Omega$ be a bounded Lipschitz domain in $\R^N$.  Let $\B$ be a linear partial differential operator, $\mathbb{C}$-elliptic of order $k$, with $k< N$. Let $f_1, f_2, h$ be continuous functions satisfying Assumptions ${\rm (H_1)}, (H_{2 q}), (H_{3 p})$, with $p=q=1$ and ${\rm (H_4)}$. Assume also that $\overline {\mathcal F}\in \mathcal G \mathcal M(\Omega)$. Then, for any $u\in BV^{\B}(\Omega)$ and $v\in\Mcal(\Omega; \R^l)$, 
        \begin{align}\label{repk<N}
            \overline{\mathcal F}(u,v)=
            \int_{\Omega} d\mathcal{Q}_\B h(\B u)(x)				
				+\int_\Omega  g\left(u(x), \frac{d v^a}{d {\mathcal L}^N}(x)\right)\,dx
				+\int_{\bar\Omega} \f1min \,
				(f_2^{\ast\ast})^\infty\left(\frac{d v^s }{d|v^s|}(x)\right)\,d|v^s|(x),
        \end{align}
         with $g(a,b) := \min\left\{ f_1(a)f_2^{\ast\ast}(b_1) + f_1^{\min}(f_2^{\ast\ast})^\infty(b_2) : b_1, b_2\in\R^N, b_1 + b_2=b  \right\}$
           and  $f_1^{\min} := \inf_{a\in\R^N} f_1(a)$.
   \end{Theorem}

   \begin{proof}
       First note that, thanks to Assumptions $(H_1), (H_{2 q})$-$(H_{3p})$, ($p=q=1$) and \cite[Theorems 1.2 and 3.7]{DG24}, for  any sequences $(u_j, v_j)_j\subset W^{\B, 1}(\Omega)\times L^1(\Omega; \R^l)$ with $\sup_j F(u_j, v_j)<\infty$, we can extract (non-relabeled) subsequences such that $u_j \overset{\ast}{\rightharpoonup} u$ in $BV^{\B}(\Omega)$, and $v_j\overset{\ast}{\rightharpoonup} v$ in $\Mcal(\overline{\Omega}; \R^l)$.
        Hence, the convergence in definition \eqref{Frelax} is well posed.  
       Its integral representation will be achieved in two steps, first we will show the lower bound and then the upper bound proving that the lower bound is sharp.
       \smallskip
       
       {\bf\underline{Lower bound.}} Let $(u_j, v_j)_j\subset W^{\B, 1}(\Omega)\times L^1(\Omega; \R^l)$ be such that $u_j\weaklystar u$ in $BV^{\B}(\Omega)$ and $v_j\weaklystar v$ in $\Mcal(\Omega; \R^l)$. Due to the structure of the functional, we have that 
          \begin{align}
              \liminf_{j\to \infty} \mathcal F(u_j, v_j) \geq \liminf_{j\to \infty} \int_\Omega f_1(u_j(x))f_2(v_j(x))dx + \liminf_{j\to \infty} \int_\Omega h(\B u_j(x))dx. \label{spliliminf}
          \end{align}
        An application of Theorem \ref{thm:RelaxationB} and Remark \eqref{remOmegachiuso} lead us to estimate the second integral on the right-hand side of \eqref{spliliminf}. Indeed, 
            \begin{align*}
                \liminf_{j\to \infty} \int_\Omega h(\B u_j(x))dx &\geq \int_\Omega Q_{\B}h\left({d\B u\over d\L^N }(x)\right)dx + \int_\Omega (Q_{\B}h)^{\#}\left({d\B^s u\over d|\B^su|}(x)\right)d|\B^su|(x)
                \\
                &=
                \int_\Omega Q_{\B}h\left({d\B u\over d\L^N }(x)\right)dx + \int_\Omega (Q_{\B}h)^{\infty}\left({d\B^s u\over d|\B^su|}(x)\right)d|\B^su|(x).
            \end{align*}
        Regarding the first integral on the right-hand side of \eqref{spliliminf}, we observe that by \eqref{defg} and \cite[proof of the lower bound in Theorem 1.1, (which we refer to for details)]{KKZ23} guarantees that
             \begin{equation}
             \notag
                 \liminf_{j\to \infty} \int_\Omega f_1(u_j(x))f_2(v_j(x))dx\geq \int_\Omega g\left(u(x), {dv^a\over d{\mathcal L}^N}(x)  \right)dx + \int_{\overline{\Omega}} f_1^{\min} (f_2^{\ast\ast})\left( 
                {dv^s\over d |v^s|}(x) \right) d|v^s|(x).
             \end{equation}

         {\bf\underline{Upper bound.}} The proof of the upper bound is split into two steps.\\
         {\underline{Step 1}}   Consider $u\in BV^\B (\Omega)$ and $v\in L^1(\Omega; \R^l)$.
         First, we recall that $\mathcal F$ satisfies the assumptions of Lemma \ref{FcalRadon}, hence for every $(u,v) \in BV^\B (\Omega) \times \mathcal M(\overline \Omega)$ it is a measure absolutely continuous with respect to $|\B u|+ |v|+ \mathcal L^N$. Moreover, by assumption, we know that  $\mathcal F \in \mathcal G \mathcal M(\Omega)$, hence there exists a function $\mathfrak{f}$ such that it is the density of $\mathcal F$ with respect to the Lebesgue measure.
         From the proof of the lower bound, we know that for a.e. $x_0 \in \Omega$, every $a \in \mathbb R^d$, $\xi \in \mathbb R^{d\times N}, G,H \in \mathbb R^{d\times d \times N}_{sym}$ and $B$ as in \eqref{fgothic}, it results that 
          \begin{align}\label{thisone}g(a,G) + \mathcal Q_\B h(B)\leq \mathfrak{f}(x_0,a,\xi, H, G,\dots, B).
          \end{align}
        To represent $\overline{\mathcal F}$, with respect to the $\mathcal L^N$ measure, it will suffice to prove the reverse inequality of \eqref{thisone}. To this end, we observe that 
        \begin{align}
        &m\left(a+ \xi(x-x_0)+ \frac{1}{2}H(x-x_0)\cdot (x-x_0)+ \dots+ p_k(x-x_0), G, Q(x_0,\varepsilon)\right)\leq \notag\\
        &\quad\mathcal F\left(a+ \xi(x-x_0)+ \frac{1}{2}H(x-x_0)\cdot (x-x_0)+ \dots+ p_k(x-x_0), G, Q(x_0,\varepsilon))\right), \label{mleqFcal}
        \end{align}
        where $p_k(x-x_0)$ is a $k$-order form whose $\mathcal B p_k(x-x_0)$ coincide with $\mathbb B$. 
        Since the function $u(x):=a+ \xi(x-x_0)+ \frac{1}{2}H(x-x_0)\cdot (x-x_0)+ \dots+ p_k(x-x_0)$ is bounded in $Q(x_0,\varepsilon)$ and $v(x):=G \in \mathbb R^l$, we can argue as in the proof of \cite[Theorem 1.1, upper bound, step 1]{KKZ23}.

From the definition of $g$, for any fixed $\eta>0$, for $\mathcal{L}^N$-a.e. $x\in\Omega$, we have
             \begin{equation}
                 \label{almoptg}
                 f_1(u(x))f_2(v^\eta_{\rm osc}) + f_1^{\rm min} (f_2^{\ast\ast})^\infty(v^\eta_{\rm conc})\leq  g(u(x), v(x)) +\eta,
             \end{equation}
        where, as in Lemma \ref{decompositionlemma}, $v$ has been decomposed as $v(x) = v^\eta_{\rm osc}(x) + v^\eta_{\rm conc}(x)$. Note that for each $x$, the set of admissible choice for $v^\eta_{\rm osc}(x)$ is always non-empty and open. In addition, $u$ and $v$ are measurable and $g$ and $f_1$ are continuous. As a consequence, it is possible to choose $v^\eta_{\rm osc}$ as a measurable function. By coercivity of $f_2$, relation \eqref{almoptg} implies that $v^\eta_{\rm osc}, v^\eta_{\rm conc}\in L^1(\Omega; \R^l)$. \\
        Thanks to \cite[Lemma 2.4]{KKZ23} along with the fact that $v^\eta_{\rm conc}\in L^1(\Omega; \R^l)$, there exists $(v^{\eta, \delta}_{\rm conc})_{\delta>0} \subset L^1(\Omega; \R^l)$ such that 
            \begin{equation}\label{3.8}
                v^{\eta, \delta}_{\rm conc} \xrightarrow{\delta\to 0} v^{\eta}_{\rm conc}\quad \mbox{strictly and hence weakly*}\hspace{0.1cm} \mbox{in } \Mcal(\overline{\Omega}; \R^l),
            \end{equation}
        as well as
             \begin{equation}
                 \notag
                 v^{\eta, \delta}_{\rm conc} \to 0 \quad \mbox{in measure.}
             \end{equation}
        Let $\gamma>0$ and choose $u^\gamma_{\rm min}\in\R^N$ such that $f_1^{\rm min}\leq f_1(u^\gamma_{\rm min})\leq f_1^{\rm min} + \gamma.$
        Now, our aim is to modify the target function $u$ near the (small) sets where $v^{\eta, \delta}_{\rm conc}\neq 0$ and replace it with $u^\gamma_{\rm min}$. To that end, we heavily use the fact that we are in dimension $N\geq 2$ and that the order of the operator is $k <N$.        
     It is possible to find a function $\varphi\in C^\infty(\overline{B_1(0)}\setminus\{0\}) \cap W^{k,1}(B_1(0))$ such that $\varphi(x)\to\infty$ as $x\to 0.$        For $s>0$, let
             \begin{equation}\label{def:varphi}
                 \varphi_s(y):= 
                 \begin{cases}
                     0 & \varphi(y)\leq {1\over s};\\
                     s\varphi(y) - 1 & {1\over s}<\varphi(y) < {2\over s};\\
                     1 & \varphi(y)\geq {2\over s}.
                 \end{cases}
             \end{equation}
        Note that $0\leq \varphi_s\leq 1$ and the support of $\varphi_s$ shrinks to $0$ as $s\to 0$, since $\varphi(x)\to \infty$ as $x\to 0$. In particular, according to the order $k$ of the operator $\mathcal B$ one could choose $\varphi$ going sufficiently slowly to $+\infty$ as $x \to 0$ so that the same happens to its $k-$th order derivatives. In particular, one has that 
            \begin{equation}
                \notag
             \|\nabla^m \varphi_s\|_{L^1}\to 0 \hbox{ and }   \|\B \varphi_s\|_{L^1}\to 0 \quad\mbox{as } s\to 0,
            \end{equation}
        for all $m\in \{0,\dots, k-1\}$.    
        At this point, we observe that we can avoid the truncation argument presented in the
    \cite[First Step of Subsection 3.2]{KKZ23}. Hence, we define 
             \begin{equation}
                 \notag
                 \widetilde{u}_{\delta} (x):= (1-h_\delta(x)) u(x) + \Psi_\delta(x)u^\gamma_{\rm min}, 
             \end{equation}
      where      \begin{equation}
                \notag
                \Psi_\delta(x):= \sum_{1\leq j\leq J(1/\delta)} \varphi_{s(\delta)}(x-x_j^\delta),
            \end{equation}
        with $s(\delta)>0$, $(x_j^\delta)_{j, \delta} \subset Q(x_0,\varepsilon)$, $J(1/\delta)\in\Nb$ and $r(\delta)$ being given by \cite[Lemma 2.4]{KKZ23} in order to obtain \eqref{3.8}. We can choose $s(\delta)>0$ such that $s(\delta)\to 0$ as $\delta\to 0$ but still slow enough so that $\varphi(x)\geq {2/s(\delta)}$ for any $|x|\leq r(\delta)$. This implies that $\varphi_{s(\delta)}(x)=1$ for any $|x|\leq r(\delta).$
        For $\delta>0$ small enough, $\Psi_\delta$ has support contained in a union of disjoint balls centered at $x_j^\delta$ with vanishing radii as $\delta\to 0$. Furthermore, $0\leq \Psi_\delta(x)\leq 1$, $\|\nabla^m \Psi_\delta\|_{L^1}\to 0$ for all $m\in \{0,\dots k-1\}$ and $\|\B\Psi_{\delta} \|_{L^1}\to 0$ as $\delta \to 0$.  Therefore, we get 
            \begin{align*}
                \|\widetilde{u}_{\delta} - u\|_{L^1}\xrightarrow{\delta\to 0} \hbox{ and }\lim_{\delta\to 0} \int_{Q(x_0,\varepsilon)} {\rm d}\mathcal Q_\B h(\B \widetilde{u}_{\delta})(x) = \int_{Q(x_0,\varepsilon)}{\rm d}\mathcal Q_\B h(\B u)(x).
                \end{align*}           
        Note also that 
             \begin{align*}
                 \B \widetilde{u}_{\delta} &= (1-\Psi_\delta) \B u + (u^\gamma_{\rm min} - u)\B \Psi_\e+ \sum_{|\alpha|=k} \sum_{|j|\leq \alpha}B_\alpha 
                 \begin{pmatrix}
                     \alpha\\
                     j
                 \end{pmatrix}
                 \partial^{\alpha-j}\Psi_\e (\partial^j u).
             \end{align*}
        Moreover, 
            \begin{equation}
                \notag
                v^{\eta, \delta} := v^\eta_{\rm osc} + v^{\eta, \delta}_{\rm conc} \weakcon v=G \quad \mbox{as } \delta\to 0 \quad \mbox{in } \Mcal(\Omega; \R^l).
            \end{equation}
        Hence, 
            \begin{align}
                \overline{\Fcal}(u, v, Q(x_0,\varepsilon)) &\leq \liminf_{\delta\to 0} {\mathcal F}(\tilde u_\delta, v^{\eta, \delta}, Q(x_0,\varepsilon))\notag\\
                &= \liminf_{\delta\to 0} \int_{Q(x_0,\varepsilon)} {\rm d}\mathcal Q_\B h(\B \tilde u_\delta)(x) + \int_{Q(x_0,\varepsilon)} f_1(\tilde u_\delta(x))f_2^{\ast\ast}(v^{\eta, \delta}(x))dx\notag\\
                & = \int_{Q(x_0,\varepsilon)} d\mathcal Q_\B h(\B u)(x)+ \liminf_{\delta\to 0}\int_{Q(x_0,\varepsilon)} f_1(\tilde u_\delta(x))f_2^{\ast\ast}(v^{\eta, \delta}(x))dx\notag\\
                &= {\mathcal Q}_\B h(B) \varepsilon^N +  \liminf_{\delta\to 0}\int_{Q(x_0,\varepsilon)} f_1(\tilde u_\delta(x))f_2^{\ast\ast}(v^{\eta, \delta}(x))dx. \notag
            \end{align}
Then the proof proceeds identically to the one of 
to \cite[Step 1 of Section 3.2]{KKZ23}, hence it is omitted.
     Thus, recalling that $v=G$ and it is constant and that $u$ is a polynomial, we can conclude that
              \begin{align}
              \label{todivide}
                  \overline{\Fcal}(u, v, Q(x_0,\varepsilon)) \leq \int_{Q(x_0,\varepsilon)} g(u,G)dx +  \mathcal Q_\B h(B ) \varepsilon^N.
              \end{align}
        Finally recalling the definition of $\mathfrak{f}$ in \eqref{fgothic} and the inequality in \eqref{mleqFcal}, we can conclude that
        \[
        \mathfrak{f}(x_0,a,\xi, H, G, \dots, B) \leq g(a, G) + \mathcal Q_\B h(B),
        \]
        for $\mathcal L^N$-a.e. $x_0\in \Omega$. Indeed, it suffices to divide the right-hand side of \eqref{todivide} by $\varepsilon^N$ and take the limit as $\varepsilon \to 0$.
       Therefore, 
        $$
        \overline{\mathcal F}[(u,v)] \leq \int_{\Omega} g(u,v){\rm d}x + \int_\Omega \mathcal Q_\B h(\mathcal B u) {\rm d}x.$$

        To achieve the representation on $BV^\B(\Omega)\times L^1(\Omega; \R^l)$ it suffices to consider a sequence $(u_\varepsilon)_\e$ which area strictly converges to $u$, and apply the reverse Fatou's lemma in view of the (upper semi-) continuity of $g(\cdot, v)$. For the existence of such a sequence we refer to arguments analogous to those in Theorem \ref{Theorem1.8AAPR20}. 

        So the step 1 is concluded observing that 
        \[
        \overline{\mathcal F}(u,v,\overline \Omega)\leq \int_\Omega g(u,v){\rm d}x + \int_\Omega d \mathcal Q_\B h(\mathcal B u), 
        \]
        for every $u \in BV^\B(\Omega)\times L^1(\Omega;\mathbb R^l)$.
        \smallskip

          \underline{Step 2.} When  $u\in BV^{\B}(\Omega)$ and $v\in\Mcal(\overline{\Omega}; \R^l)$ using the approximation lemma \cite[Lemma 2.5]{KKZ23} and the Step 1 together with the lower semicontinuity of $\overline{\Fcal}$; Dominated convergence's theorem and Reshetnyak's continuity theorem, one can conclude that  
            \begin{equation}\notag
                \overline{\Fcal}(u, v) \leq \int_\Omega g(u(x), v^{\rm a}(x))dx + \int_{\overline{\Omega}} f_1^{\rm min}(f_2^{\ast\ast})^\infty\left( {dv^s\over d|v^s|}(x)  \right)d|v^s|(x) + \int_\Omega {\rm d}\mathcal Q_\B h(\B u)(x),
            \end{equation}
    as desired.
   \end{proof}

\subsection{Case with no concentration effects}
We are in position to state the result concerning the integral representation of $\Fcal$ when the operator $\B$ has order $k$, with $k\geq N$, for $N\geq 2$.
Since the order of the operator $\mathcal B$ is greater than or equal to $N$,  we can freeze the value of $f_1$ at $u(x_0)$, for $x_0 \in \overline \Omega$, thus prohibiting concentration effects. This is due to the fact that $u$ is continuous also on $\overline \Omega$. We restate Theorem \ref{thm:relKgeq N}, for the readers' convenience.

 \begin{Theorem}
      Let $\Omega$ be a bounded open set with Lipschitz boundary.  Let $\B$ be the $\mathbb C$-elliptic linear partial differential operator of order $k$ as in \eqref{Budef}, with $k \geq N$ ($k>N$ if $N\leq 2$). Let $f_1, f_2, h$ be continuous functions satisfying Assumptions {\rm $(H_1), {\rm (H_2q)}, {\rm (H_3p)}$}, with $p=q=1$ and {\rm $(H_4)$}. Let $ \Fcal(u,v)$ and $\overline{\Fcal}(u,v)$ be given by \eqref{Flocalized} and \eqref{FlocrelaxB}. Then,
                 \begin{align}\label{repk>1bis}
                     \overline{\Fcal}(u,v) &= \int_{\overline{\Omega}} f_1(u(x))df_2^{\ast\ast}(v)(x)+\int_\Omega d\mathcal Q_{\B}h (\B u)(x). 
                 \end{align}
  \end{Theorem}

\begin{proof}[Proof] 
The proof follows, also in this case, the same steps of \cite[Theorem 1.2]{KKZ23}: we will prove first a lower bound and then an upper bound.  Let $u\in BV^\B(\Omega)$ and $v\in\Mcal(\overline{\Omega}; \R^l)$.
 
  {\bf\underline{Lower Bound.}} Without loss of generality, we can assume that $f_2$ and $h$ are convex and $\B$-quasiconvex, respectively, since $f_2\geq f_2^{\ast\ast}$ and $h\geq Q_{\B} h$.
  Let $(u_j, v_j)_j\in W^{\B, 1}(\Omega)\times L^1(\Omega; \R^l)$ such that $u_j\to u$ in $L^1(\Omega;\mathbb R^d)$ and $v_j\weaklystar v$ in $\mathcal M(\overline{\Omega}; \R^l)$. Assume that, up to a (not relabeled) subsequence, the limit $\liminf_{j\to\infty} \int_\Omega f_1(u_j)f_2(v_j) + h(\B u_j) dx <+\infty$ exists. 
       The proof of the lower bound is immediate. Indeed, it suffices to exploit the uniform convergence of $u_k\to u$, which is a consequence of Sobolev's and Rellich's type embedding theorems.
  So
  \begin{align*} &\liminf_{j\to\infty} \int_\Omega f_1(u_j(x))f_2(v_j(x)) + h(\B u_j(x)) dx \\
  &\quad\geq \liminf_{j\to\infty} \int_\Omega (f_1(u_j(x))-f_1(u(x)))f_2(v_j(x))dx + \liminf_{j\to +\infty}\int_\Omega f_1(u(x))f_2(v_j(x)) + h(\B u_j(x)) dx\\
  &\quad\geq \liminf_{j\to +\infty} \int_\Omega f_1(u(x))f_2(v_j(x)) dx +\liminf_{j\to +\infty}\int_\Omega h(\B u_j(x))dx\\
  &\quad\geq \int_\Omega f_1(u(x))d f_2(v)(x) + \int_\Omega d Q_{\B}h(\B u)(x). 
  \end{align*}
       Here, we have used in the first inequality the superadditivity of the liminf, in the second one the uniform convergence of $f_1(u_n) $ to $f_1(u)$ and the boundedness of $\int_\Omega f_2(v_n)dx$ and in the third one the lower semicontinuity result provided by Theorem \ref{thm:RelaxationB}.
       taking into account that $f_1(u(\cdot))$ is uniformly continuous in $\overline \Omega$.

       {\bf\underline{Upper bound.}} In this case by Proposition \ref{proprel**}, it is enough to fix $u \in BV^\B(\Omega)$ since the representation of $\overline{\mathcal F}(u,v,\Omega)$ amounts at providing an integral representation for the functional $\mathcal F^{\ast \ast}(u,v,\overline{\Omega})$ in \eqref{F**locrelaxB}.
In fact, relying on the fact that one can extend $u$ outside $\Omega$ without charging $\partial \Omega$ up to $(k-1)$-th derivatives,
\begin{align*}
\mathcal F^{\ast \ast}(u,v;\Omega)
 &\leq \inf\left\{\liminf_{j\to +\infty} \int_{\overline{\Omega}} f_1(u(x))f_2^{\ast \ast}(v_j(x))dx + \int_\Omega d Q_{\B} h(\B u)(x): v_j \overset{\ast}{\rightharpoonup} v \hbox{ in } \mathcal M(\overline \Omega ;\mathbb R^l)\right\} \\
 & \leq \int_\Omega d Q_{\B} h(\B u)(x)+ \inf\left\{\liminf_{j\to +\infty} \int_{\overline{\Omega}} f_1(u(x))f_2^{\ast \ast}(v_j(x))dx: v_j \overset{\ast}{\rightharpoonup} v \hbox{ in } \mathcal M(\overline \Omega ;\mathbb R^l)\right\}.
 \end{align*}

 Now the relaxation of the latter functional again can be obtained as a particular case of Theorem \ref{thm:RelaxationB}, applied to $\mathcal B= {\rm Id}$ and density $f:\Omega \times \mathbb R^l$ defined as $f(\cdot, \cdot):= f_1(u(\cdot))f_2^{\ast \ast}(\cdot)$ taking into account 
 the uniform continuity of $f_1(u(\cdot))$.
 This concludes the proof.
\end{proof}

\begin{Remark}\label{k=n=2}
When $k=N=2$ and the operator $\B= D^2$, and $\Omega$ is of class $C^2$, up to a finite number of points (see {\rm \cite{D}}), then there is a continuous embedding of $BH(\Omega;\mathbb R^d)$ into $C(\overline{\Omega};\mathbb R^d)$ and then in our relaxation procedure there will be no concentration effects.
In general, the continuous embeddings in the case $k=N=2$ might not hold and the concentration effects detected in Theorem \ref{thm:relk<N,N>=2} can appear.
\end{Remark}

\subsection{Relaxation result in one-dimensional case}
Throughout this section, we assume that $N=1$. As pointed out in \cite[comments on Theorem 1.1]{DG20}, for $N=1$, the only $\mathbb C$-elliptic operators $\B$ of $k$-th order are of the form  \eqref{def:ellipticoperN1}, i.e.,
            \begin{equation}
                \notag
                \B u = B{d^ku\over dx^k}, \;\;\; u:(\alpha,\beta)\subset\R \to V
            \end{equation}
for some $B\in {\rm Lin}(V; V)$. Therefore, we have that $BV^\B(\alpha,\beta) = BV^k((\alpha,\beta); V),$ where $BV^k((\alpha,\beta); V)$, introduced in \cite{ADC}, is the space of function $u\in L^1((\alpha,\beta); V)$ whose $k$-th derivative is the sense of distribution is a measure with total bounded variation. Equivalently, it is the space $\{u \in W^{k-1,1}((\alpha,\beta);V): D^k u \in \mathcal M((\alpha,\beta);W)\}$, see \eqref{BVBequiv}. As in \cite{KKZ23}, the blow-up technique introduced in \cite{FM92} can be adopted 
to obtain a different characterization of integral functionals defined on interval $\Omega:=(\alpha,\beta)$. Indeed,  any finite measure $\nu$ on $\Omega\subset\R$ can be written as
    \begin{equation}
        \label{atomic+diffmeasure}
        \nu= \nu^0 + \nu^{\rm diff},
    \end{equation}
where $\nu^{\rm diff}$ is the diffuse part and  $\nu^0$ is the zero-dimensional part of $\nu$ charging point, i.e., $\nu^0:= \nu \res S^0,$ with $S^0 = S^0(\nu) := \{x\in\Omega: |\nu|(\{x\})>0 \}.$ In particular, $\nu^0({x_0}) = \nu (\{x_0\})$ for any $x_0\in\Omega$ and zero for all but countably many, $|\nu^0|= |\nu|^0$ and $|\nu^{\rm diff}|=|\nu|^{\rm diff}$. 
Moreover, thanks to \eqref{def:ellipticoperN1}, the following decomposition follows
            \begin{equation}
             \notag
             \B u= \B u\mathcal L^1 + (\B u)^s =(\B u) \mathcal{L}^1 + (\B u)^{j} + (\B u)^{c} = B\frac{d^k u}{d x^k}\mathcal{L}^1 + B\left[\frac{d^{k-1}u}{d x^{k-1}}\right] \mathcal{H}^0\res S_{\frac{d^{k-1}u}{d x^{k-1}}} + B \left(\frac{d^{k}u}{d x^{k}}\right)^{c},
         \end{equation}
where for $\mathcal{H}^0$-a.e. $x\in S_{\frac{d^{k-1}u}{d x^{k-1}}}$, $[\frac{d^{k-1}u}{d x^{k-1}}](x):= \frac{d^{k-1}u}{d x^{k-1}}(x^{+})-\frac{d^{k-1}u}{d x^{k-1}}(x^{-})$, with $\frac{d^{k-1}u}{d x^{k-1}}(x^{\pm}) := \lim_{y\to x^{\pm}} \frac{d^{k-1}u}{d x^{k-1}}(y)$ and $\frac{d^ku}{dx^k}$ is the approximate derivative of $u$. 
Comparing it with \eqref{atomic+diffmeasure}, the zero-dimensional part $(\B u)^0$ of $\B u$ is given by $(\B u)^{j}$ and  the diffuse part $(\B u)^{\rm diff}$ is  $B\frac{d^{k}u}{d x^{k}}u \mathcal{L}^1+ (\B u)^{c}$. 
The reminder of this subsection is devoted to the proof of Theorem \ref{thm:relN=1}.

  \begin{Remark}\label{rem1=N<K}
    When $N=1$ and the order of the operator $\B$ is $k>N$, from Theorem \ref{thm:relKgeq N}, we have no concentration phenomena and formula \eqref{repk>1bis} turns into the following one
  \begin{align}\notag
                     \overline{\Fcal}(u,v) &= \int_{\Omega} f_1(u(x))df_2^{\ast\ast}(v^{\rm diff})(x)+\int_\Omega dh^{\ast\ast}(\B u)(x) + \sum_{x\in S^0}f_1(u(x))(f_2^{**})^\infty(v^0( \{x\}))\notag\\
                     &\quad + f_1\left(u(\alpha^{+})\right) (f_2^{**})^\infty(v^0(\alpha^+))+ f_1(u(\beta))(f_2^{**})^\infty(v^0(\beta^{-})).\notag
                 \end{align}
  \end{Remark}
  \begin{proof}[Proof of Theorem \ref{thm:relN=1}]
    Let $u\in BV^\B(\Omega)$ and $v\in\Mcal(\Omega; \R^l)$. In view of Lemma \ref{FcalRadon}, it suffices to show our claim for appropriate localized measure $\overline{\Fcal}$, which is a signed Radon measure defined as in  Lemma \ref{FcalRadon}. 
    
    Since $k=N=1$, $\B u= B(\nabla u)\mathcal{L}^1 + \symB(\nu_u)[u] \mathcal{H}^0\res S_u + \B u^{c},$
where for $\mathcal{H}^0$-a.e. $x\in S_u$, $[u](x):= u(x^{+})-u(x^{-})$, with $u^{\pm} := \lim_{y\to x^{\pm}} u(y)$ and $\nabla u$ is the approximate gradient of $u$ and $B\in \Lin(V; V)$ (see also \cite[Lemma 2.4]{R19_1}). Note that since in one-dimension the measure-theoretic normal vector is trivial ($\nu_u=1$), we get that $\symB(\nu_u)[u]=B[u]$. The proof follows the same steps of \cite[Theorem 1.2]{KKZ23} and we summarize it, for the reader's convenience. \\
{\bf\underline{Lower Bound.}}  Setting
                 \begin{equation}
                     \label{thetadef}
                     \theta\coloneqq \Lcal^1+|v| + |\B u|, \;\;\; \theta\in\Mcal(\overline{\Omega}),
                 \end{equation}
    we are going to prove via blow-up technique the following
                \begin{align}\label{2}
                    {d\overline{\Fcal}^{\rm diff} \over d\theta^{\rm diff}}(x_0)\geq f_1(u(x_0)){df_2^{\ast\ast}\over d\theta^{\rm diff}}(x_0) + {dh^{\ast\ast}(\B u^{\rm diff})\over d\theta^{\rm diff}}(x_0) &\;\;\; \mbox{for } \theta^{\rm diff}\hbox{-a.e. } x_0\in\Omega,\notag\\
\end{align}
where $\theta^{\rm diff}$ denotes the diffuse part of $\theta$.
Without loss of generality, we can assume that $f_2$ and $h$ are convex (quasiconvexity reduces to convexity since $N=1$), since $f_2\geq f_2^{\ast\ast}$ and $h\geq h^{\ast\ast}$ Let $(u_j, v_j)_j\in W^{\B, 1}(\Omega)\times L^1(\Omega; \R^l)$ such that $u_j\to u$ in $L^1(\Omega;\mathbb R^d)$ and $v_j\weaklystar v$ in $\mathcal M(\overline{\Omega}; \R^l)$. Assume that, up to a (not relabeled) subsequence, the limit $\lim_{j\to\infty} \int_\Omega f_1(u_j(x))f_2(v_j(x)) + h(\B u_j(x)) dx <+\infty$ exists. For any Borel set $A\subset\Omega$, define the sequence of measures $(\mu_j)_j$ by
            \begin{equation*}
                \mu_j (A) := \int_A f_1(u_j(x))f_2(v_j(x)) + h(\B u_j(x)) dx.
            \end{equation*}
      Since $f_1$, $f_2$ and $h$ are nonnegative functions and they satisfy the growth conditions given by Assumptions $(H_1)$, $(H_{2}q)$, $(H_{3}q)$ (with $p=q=1$) and $(H_4)$, up to a subsequence, we obtain 
                \begin{equation}
                    \label{conveofmeasure1}
                    \mu_j\weaklystar \mu \quad \mbox{and}\quad |v_j|+|\B u_j|\weaklystar \lambda \quad \mbox{in } \mathcal{M}(\overline{\Omega}). 
                \end{equation}
    In view of \eqref{atomic+diffmeasure}, $\mu$ can be split as $\mu=\mu^0+\mu^{\rm diff}$. Set
           \begin{equation*}
               I_{x_0, \e}:= (x_0-\e, x_0+\e) \quad \mbox{and} \quad I:= (-1,1),
           \end{equation*}
    with $x_0\in\Omega$ and $\e>0$ being small enough so that $I_{x_0, \e}\subset\Omega$.
    \smallskip

 \underline{{Diffuse contribution.}} Let $x_0$ be a point with $\lambda(\{x_0\})=0$ such that the measure $\mu$, $f_2(v)$ and $h(\B u)$ have a finite density with respect to $\theta^{\rm diff}$ at $x_0$, i.e,
$$
\frac{d\mu}{d\theta^\diff}(x_0)=\lim_{\varepsilon \rightarrow 0}\frac{\mu(I_{x_0,\varepsilon})}{\theta^\diff (I_{x_0,\varepsilon})}<+\infty,
$$
$$
	\frac{d f_2(v)}{d\theta^\diff}(x_0)=\lim_{\varepsilon \rightarrow 0}\frac{f_2(v)(I_{x_0,\varepsilon})}{\theta^\diff (I_{x_0,\varepsilon})}<+\infty,
\quad\text{ and }\quad
\frac{d h(\B u)}{d\theta^\diff}(x_0)=\lim_{\varepsilon \rightarrow 0}\frac{h(\B u)(I_{x_0,\varepsilon})}{\theta^\diff (I_{x_0,\varepsilon})}<+\infty.
$$
In particular, for $\theta^\diff$-a.e. $x_0\in\Omega$,  $|\mu|(\{x_0\})=|\B u|(\{x_0\})=|v|(\{x_0\})=0$ and $h(\B u)(\{x_0\})=f_2(v)(\{x_0\})=0.$
In addition, for $\theta^\diff$-a.e. $x_0\in\Omega$,
    \begin{equation}
        \notag
        \frac{d f_2(v)}{d\theta^\diff}(x_0)=\frac{d f_2(v^\diff)}{d\theta^\diff}(x_0) \quad \mbox{and}\quad \frac{d h(\B u)}{d\theta^\diff}(x_0)=\frac{d h(\B u^\diff)}{d\theta^\diff}(x_0).
    \end{equation}
Set $ \vartheta_{x_0,\eps}:=\theta^\diff(I_{{x_0},\eps}),$ and note that $\vartheta_{x_0,\eps}\geq \Lcal^1(I_{{x_0},\eps})=2\eps.$ Consider a sequence $\e\to 0$ such that  $\mu\left(  \{  x_{0}-\varepsilon, x_0+\varepsilon\} \right)= \lambda(\{x_0-\varepsilon, x_0+\varepsilon\}) =0$. 
Therefore,   $\mu(\partial I_{x_0,\varepsilon})=0$. This along with the change of variables $y=x_0+\theta_{x_0, \e}x$ implies that  
\begin{align}\nonumber
\begin{aligned}
	&\frac{d|\mu|}{d\theta^\diff}(x_0) =  
	\lim_{\eps\to 0}\frac{\mu(I_{x_0,\varepsilon})}{\vartheta_{x_0,\eps}}
\\
	&=  \lim_{\eps\to 0}\lim_{j\to +\infty}\frac{1}{\vartheta_{x_0,\eps}}\int_{I_{x_0,\eps}}f_1(u_j(y))f_2(v_j(y))+ h(\B u_j(y))\,dy
\\
	&=  \lim_{\eps\rightarrow 0}\lim_{j\rightarrow +\infty}\int_{I_{0,\eps\vartheta_{x_0,\eps}^{-1}}}
	f_1(u_j(x_0)+\vartheta^k_{x_0,\eps} w_{j,\eps}(x))f_2(\eta_{j,\eps}(x))+ h(\B w_{j,\eps}(x))\,dx,
\end{aligned}
\end{align}
where we have set $w_{j, \e}\coloneqq {u_{j}(x_0+\vartheta_{x_0, \e}x)-u_j(x_0)\over \vartheta_{x_0, \e}^k}$ and $\eta_{j, \e} \coloneqq v_j(x_0+\vartheta_{x_0, \e}x).$ By coercivity of $h$ combined with the fact that ${d|\mu|\over d\theta^{\rm diff}}$ is finite, we have that $\int_{I_{0, \e(j)\vartheta^{-1}_{x_0. \e(j)} }} |\B w_{j, \e(j)}|dx$ is bounded and $w_{j, \e(j)}(0)=0$. This implies that $w_{j, \e(j)}$ is bounded on $I_{0, \e(j)\vartheta^{-1}_{x_0. \e(j)}}$ by a constant independent of $j$. Hence
     \begin{equation}
         \notag
         \vartheta^k_{x_0,\eps} w_{j,\eps}(x) \to 0\;\;\; \hbox{uniformly in } x\in  I_{0, \e(j)\vartheta^{-1}_{x_0. \e(j)}}.
     \end{equation}
In addition, we have that 
   \begin{equation}
       \notag
       u_j(x_0) \to u(x_0).
   \end{equation}
Indeed, since $u_j, u\in C^{k-1}$,
   \begin{align}
     |u_j(x_0)-u(x_0)|&\leq |u_j(x_n)-u(x_n)|+|u_j(x_0)-u_j(x_n)|+ |u(x_0)-u(x_n)|\nonumber\\
     &\leq |u_j(x_n)-u(x_n)| + \int_{x_0}^{x_n} u'_j(x)-u'(x)dx\nonumber\\
    & \leq |u_j(x_n)-u(x_n)| +  (|u_j'|+ |u'|)(x_n-x_0) \to 0
     \;\;\;\hbox{as } n\to\infty, \nonumber
   \end{align}
   since the last summand due to the $\mathbb C$-ellipticity of $\B$ is equivalent to  $(|\B u_j|+ |\B u|)(x_n-x_0) \to 0
     \;\;\;\hbox{as } n\to\infty$. Then, the proof proceeds identically to the one of the diffuse part of the lower bound in \cite[Theorem 1.2]{KKZ23}, so we can conclude that \eqref{2} holds.
\smallskip

  \underline{Contributions charging  points in the interior and on the boundary.}
 We prove the following
                \begin{align}
                    \overline{\Fcal}^0(\{x_0\})& \geq f^0_h (u(\{x_0^{+}\}), u(\{x_0^{-}\}), v^0(\{x_0\}))\;\;\;&&\mbox{for $\theta^0$-a.e. } x_0\in\Omega;\label{lbk1N1}\\
                    \overline{\Fcal}^0(\{x_0\})& \geq \inf_{z\in V} f^0_h(u(\{x_0^{+}\}), z, v^0(\{x_0\}))\;\;\;&&\mbox{for } x_0=\inf\Omega;\label{lbk1N1-1}\\
                    \overline{\Fcal}^0(\{x_0\})& \geq \inf_{z\in V} f^0_h(z, u(\{x_0^{-}\}), v^0(\{x_0\}))\;\;\;&&\mbox{for } x_0=\sup\Omega.\label{lbk1N1-2}
                \end{align}
                
                Let $x_0\in S^0(\mu)$.  Since $u\in BV^\B(\alpha, \beta)$, $u(x_0^+)$ and $u(x_0^-)$ exist in $\mathbb R^m$. 
In the following, we denote by $(\varepsilon)$ a sequence of positive reals converging to $0$ such that  
\begin{align}\label{eps-dontchargeboundary}
	|v|(\{x_0\pm \varepsilon\})=|\B u|(\{x_0\pm \varepsilon\})= \mu(\{x_0\pm \varepsilon\})=0.
\end{align}
By dominated convergence, $\mu(\{x_0\})=\lim_{\varepsilon\to 0^+} \mu(I_{x_0,\varepsilon})$, and by
\eqref{conveofmeasure1}, using the definition of $\mu_k$ and a change of variables, we deduce that
	\begin{align}\label{pointlb-1}
	 \mu(\{x_0\})	\geq \lim_{\varepsilon \to 0^+}\lim_{j\to +\infty}\varepsilon\int_I \left(f_1(u_j(x_0+\varepsilon y))f_2(v_j(x_0+ \varepsilon y)) + h( \B u_j(x_0 +\varepsilon y))\right)dy,
	\end{align}
	with $I=(-1,1)$.
Set $v_{j,\varepsilon}(y):=v_j(x_0+ \varepsilon y).$ Since $v_j \overset{\ast}{\rightharpoonup} v$ in $\mathcal M(\bar\Omega;\R^d)$
and $|v|(\{x_0\pm \varepsilon\})=0$ by \eqref{eps-dontchargeboundary}, for fixed $\varepsilon$,
\[
	\int_{I}\varepsilon v_{j,\varepsilon}(y) dy=\int_{I_{x_0,\varepsilon}} v_j(y) \,dx
	\underset{j\to\infty}{\longrightarrow} v(I_{x_0,\varepsilon})
\]
As $\varepsilon\to 0^+$, we conclude that
	 $\displaystyle{\lim_{\varepsilon\to 0^+}\lim_{j \to +\infty}\int_I\varepsilon v_{j,\varepsilon}(y)dy=v(\{x_0\}).}$
A diagonalization argument now ensures the existence of a sequence $\overline v_\varepsilon:=v_{j(\varepsilon),\varepsilon}$ such that
\begin{align}\label{rescalconv}
	 \lim_{\varepsilon \to 0^+}\int_{I}\varepsilon \overline v_\varepsilon(y)dy=v(\{x_0\}).
\end{align}	
On the other hand, setting $u_{j,\varepsilon}(y):=u_j(x_0+\varepsilon y)$, it follows that 
\begin{align}\label{upm}
	 \lim_{\varepsilon\to 0^+}\lim_{j\to +\infty}\|u_{j,\varepsilon}-u^{\pm}\|_{L^1(I)}=0,
\quad\text{where}\quad
	 u^{\pm}(y)=\left\{
	 \begin{array}{ll}
	 u(x_0^+), \hbox{ if } y \geq 0,\\
	 u(x_0^-), \hbox{ if } y  < 0.
	 \end{array}
	 \right.
\end{align}
In addition, setting $\overline u_\varepsilon:= u_{j(\varepsilon),\varepsilon}$ and $\overline v_\varepsilon:= v_{j(\varepsilon),\varepsilon},$ with a suitable $j(\eps)\to+\infty$ (fast enough), \eqref{pointlb-1} turns into
\begin{align}\nonumber
	 \mu(\{x_0\})\geq\lim_{\varepsilon \to 0^+}\varepsilon\left(\int_I f_1(\overline{u}_\varepsilon(y)) f_2(\overline{v}_\varepsilon(y)) 
		+h\left(\frac{\B \overline u_\varepsilon(y)}{\varepsilon}\right)\right)dy.
\end{align} 
Now, setting $ \hat{v}_\varepsilon(y):= \varepsilon\overline{v}_\varepsilon(y)-\varepsilon\int_I\overline v_\varepsilon dy + v(\{x_0\}),$ 
by \eqref{rescalconv}, it follows that	
\begin{align}	\nonumber
	 |\hat v_\varepsilon(y)-\varepsilon\overline{v}_\varepsilon(y)|
	= \left|\int_I \varepsilon\overline{v}_\varepsilon(y)dy - v(\{x_0\})\right| \to 0 \hbox{ as }\varepsilon\to 0.
\end{align} 	
This along with the global Lipschitz continuity of $f_2$ leads us to 
\begin{align*}	 
	 \mu(\{x_0\})\geq\lim_{\varepsilon \to 0^+} \int_I\left(f_1(\overline{u}_\varepsilon(y)) \varepsilon f_2\Big(\frac{1}{\varepsilon}\hat v_\varepsilon(y)\Big) +\varepsilon h\left(\frac{\B \overline u_\varepsilon(y)}{\varepsilon}\right)\right)dy.
\end{align*} 
Above, we may replace $f_2$ and $h$ by their recession functions. Therefore,
\begin{align}	 \label{pointlb-3}
	\mu(\{x_0\})\geq
	\lim_{\varepsilon \to 0^+} \int_I\left(f_1(\overline{u}_\varepsilon(y))f_2^\infty(\hat v_\varepsilon(y)) +h^\infty( \B \overline u_\varepsilon(y))\right)dy.
\end{align} 
As $\int_I \hat v_\varepsilon=v(\{x_0\})$ by construction, $\hat v_\varepsilon$ is admissible in the infimum defining $f_h^0$ in \eqref{deff0W-RelN=1}. Applying Lemma \ref{boundaryc}  to the integrand  $f_1 f_2^\infty+ h^\infty$, we can modify $\overline{u}_\varepsilon$ into a function ${\hat u}_\varepsilon$ which has the same values as $u^{\pm}$ on $\{\pm 1\}$, cf.~\eqref{upm}. 
	 Thus, the new function $\hat u_\varepsilon$ is also admissible in \eqref{deff0W-RelN=1} (in place of $u$) and we obtain \eqref{lbk1N1}
	from \eqref{pointlb-3},
	 more precisely,
	 \begin{align*}
	 \mu(\{x_0\})\geq
	 \lim_{\varepsilon \to 0^+} \int_I\left(f_1({\hat u}_\varepsilon(y))f_2^\infty(\hat v_\varepsilon) +h^\infty( \B{\hat u}_\varepsilon  (y))\right)dy
	 \geq f_h^0\left(u(x_0^+), u(x_0^-), v(\{x_0\})\right).
	 \end{align*} 
The proof of \eqref{lbk1N1-1} and \eqref{lbk1N1-2} is omitted since it is identical to the one in \cite[subsection 4.4]{KKZ23}.
\smallskip

{\bf\underline{Upper Bound.}}
We will show the upper bound, i.e. we will construct a recovery sequence. To this end, let $(u,v)\in BV^\B(\Omega)\times \Mcal(\overline\Omega;\R^d)$.  
Without loss of generality, in view of Proposition \ref{proprel**}), $f_2$ and $h$ can be assumed to be convex and quasiconvex, respectively (recall that when $N=1$ quasiconvexity reduces to convexity) 
and we can use the representation \eqref{F**locrelaxB} of $\overline{\mathcal F}=\mathcal F^{\ast\ast}$. By Lemma \ref{FcalRadon}, it is enough to prove the upper bounds for the density	of $\overline{\mathcal F}$ with respect to the atomic and diffuse parts of $\theta=\theta^0+\theta^\diff$ in \eqref{thetadef}. We therefore have to show that 
\begin{align}
&\begin{aligned}\label{ub1d-diffuse}
			\frac{d \overline{\mathcal F}}{d\theta^\diff}(x_0)\leq 
		f_1(u(x_0))\frac{d f_2(v)}{d\theta^\diff}(x_0)+ \frac{d h(\B u)}{d\theta^\diff}(x_0),
        && \quad\text{for $\theta^\diff$-a.e. $x_0\in \Omega$,}
		\end{aligned} 
\end{align}
and
\begin{alignat}{2}
		&\begin{aligned}\nonumber
		\overline{\mathcal F}^0(\{x_0\})\leq f_h^0\left(u^+(x_0), u^-(x_0), v(\{x_0\})\right)
		\end{aligned}
		&&\quad\text{ for $\theta^0$-a.e. $x_0\in \Omega$,}\\
		&\begin{aligned}
		\nonumber
		\overline{\mathcal F}^0(\{x_0\})\leq 
		\inf_{z\in\R^m} f_h^0\left(z, u^-(x_0), v(\{x_0\})\right)
		\end{aligned}
		&&\quad\text{for $x_0= \sup\Omega$, and}\\
		&\begin{aligned}\nonumber 
		\overline{\mathcal F}^0(\{x_0\})\leq 
		\inf_{z\in\R^m} f_h^0\left(u^+(x_0),z, v(\{x_0\})\right)
		\end{aligned}
		&&\quad\text{for $x_0= \inf\Omega$}.
\end{alignat}	
Here, recall that $\theta=\Lcal^1+|v|+|\B u|\in \Mcal(\overline\Omega)$ (with $|\B u|(\partial\Omega):=0$), which is always possible taking into account.

\underline{Diffuse contributions.} 
We will show \eqref{ub1d-diffuse}. Recall that, as above, $I_{{x_0},\eps}:=(x_0-\eps,x_0+\eps)$ and
$\vartheta_{x_0,\eps}:=\theta^\diff(I_{{x_0},\eps})$. Let $x_0\in \Omega$ be such that the Besicovitch derivatives 
$\frac{d \overline{\mathcal F}(u,v,\cdot)}{d\theta^\diff}(x_0)$,
$\frac{d f_2(v)}{d\theta^\diff}(x_0)$ and $\frac{d W(\B u)}{d\theta^\diff}(x_0)$ exist with finite values. In particular, $|v|(\{x_0\})=|\B u|(\{x_0\})=0$, $x_0$ is not a jump point of $u$.
Since $\B$ takes the form \eqref{def:ellipticoperN1} and in view of the Sobolev embedding theorems, $u$ has a continuous representative at $x_0$. In addition,  $\vartheta_{x_0,\eps}\to 0$ as $\eps\to 0$.
By Proposition~\ref{proprel**} applied to $f_2=f_2^{\ast\ast}$ and $h=h^{**}$ (recall that since $N=1$, $Qh=h^{\ast \ast}$),  it suffices to 
find a sequence $(v_j)_j\subset L^1(\Omega;\R^d)$ 
such that $v_j\weaklystar v$ in $\Mcal(\overline\Omega;\R^d)$ and
\begin{align}\label{1drec-diff-0}
\begin{aligned}
	&\liminf_{\eps\to 0}\limsup_{j\to\infty}\frac{1}{\vartheta_{x_0,\eps}}\int_{I_{x_0,\eps}}
	\left(f_1(u(x))f_2(v_j(x))\,dx+dh(\B u)(x)\right) \\
	&\quad\leq \lim_{\eps\to 0} \frac{1}{\vartheta_{x_0,\eps}} \int_{I_{x_0,\eps}}
	\left(f_1(u(x_0))df_2(v)(x)+dh(\B u)(x)\right).
\end{aligned}
\end{align}
As the contribution of $dh(\B u)(x)$ appears on both sides, \eqref{1drec-diff-0} can be reduced to 
\begin{align}\nonumber
\begin{aligned}
	\liminf_{\eps\to 0}\limsup_{j\to\infty}\frac{1}{\vartheta_{x_0,\eps}}\int_{I_{x_0,\eps}}
	f_1(u(x))f_2(v_j(x))\,dx 
	&\leq f_1(u(x_0)) \lim_{\eps\to 0} \frac{1}{\vartheta_{x_0,\eps}} \int_{I_{x_0,\eps}}
	df_2(v).
\end{aligned}
\end{align}

Then, in view of \eqref{def:ellipticoperN1} and the Sobolev embedding theorems, the proof develops identically to the one of \cite[Proposition 4.6]{KKZ23},  hence it is omitted.

\underline{Contributions charging points in the interior and on the boundary.} The proof is identical to the one of \cite[Subsection 4.7]{KKZ23}: it suffices to replace $D$ by $\B$ and takes into account that $\B$ as in \eqref{def:ellipticoperN1}.
\end{proof}

\begin{Remark}\label{remvoletal}
All the relaxation theorems in this section can be easily proven in the space $BV^\B(\Omega)\times \mathcal M(\Omega;\mathbb R^l$). Indeed, similar integral representations hold: it suffices to drop all contributions on $\partial \Omega$ in the integrals. In addition, the energetic models as the ones in \eqref{G} can be complemented with mass constraints on $v$. In fact, they can be easily dealt with, up to slight modifications of the recovery sequences.
\end{Remark}

\subsection{Application of Theorem \ref{thm:RelaxationB} to the variable exponent setting}
Theorem \ref{thm:RelaxationB} allows us also to deduce the following lower semicontinuity result, which finds application in imaging problems (see, e.g., \cite{BKP10, CLR}) and whose proof, being identical to the one of \cite[Proposition 3.1]{BHHZ}, is omitted. 

Let $p \in C(\overline\Omega) $ be a continuous function such that $1\leq p(x)<+\infty$, let $Y:=\{x \in \Omega: p(x)=1\}$, and let $\B$ be a $k$-th homogeneous linear partial differential operator as in \eqref{Budef}. We introduce the space $BV^{\B,p(\cdot)}(\Omega):= BV^{\B}(\Omega) \cap W^{\B, p(\cdot)}(\Omega \setminus Y)$, where 
$ W^{\B, p(\cdot)}(\Omega \setminus Y):=\{u \in L^{p(\cdot)}(\Omega; V): \B u \in L^{p(\cdot)}(\Omega; W)\}$. Here a function $u$ belongs to $L^{p(\cdot)}(\Omega;V)$ if and only if $ \int_{\Omega} |v(x)|^{p(x)}dx <+\infty$ (see \cite{HHbook} for the properties of the variable exponent spaces). 

\begin{Theorem}\label{thm:asBHHZsci}
Let $f:W \to [0,+\infty)$ be a $\mathcal B$-quasiconvex function with linear growth from above and below, i.e. such that $(H_{3p})$ holds for $p=1$. Then,
the functional $$BV^{\B,p(\cdot)}(\Omega;\mathbb R^l) \ni u \mapsto \mathcal F_{\B, p(\cdot)}(u):=\int_\Omega f(\B u(x))^{p(x)}dx + \int_\Omega f^\infty\left(\frac{d \B^s u}{d |\B^s u|}(x)\right) d |\B^s u|(x)$$ is lower semicontinuos with respect to the strong convergence in $L^{p(\cdot)}(\Omega;W )$, i.e. along sequences $(u_j)_j$ such that $\lim_j \int_\Omega |u_j-u|^{p(x)}dx =0$.
\end{Theorem}

\section{Relaxation in the Sobolev setting}\label{relBFLhigh}

For $k$-th homogeneous linear partial differential $\mathbb C$-elliptic operators $\B$, we will prove the Sobolev counterpart of Theorem \ref{thm:RelaxationB} providing an integral representation of the relaxed functional in the superlinear case, in the same spirit of \cite[Theorems 1.1, 3.6]{BFL00}. Our analysis deals with relaxation of functionals depending on both $u$ and $v$, when $v$ is an $L^p$ field and not a measure.  It is also worth stressing that some of the cases we will present, constitute the weakly lower semicontinuous envelope of the functional appearing in Theorem \ref{thm:relk<N,N>=2}, when $(u,v) \in W^{\B,p}\times L^p$, $p>1$.
In addition, we will discuss the Sobolev relaxation, analogous to Theorems \ref{thm:relk<N,N>=2} and \ref{thm:relN=1}, i.e., when $u \in W^{\B.p}$, with $p>1$ and $v \in \mathcal{M}(\overline{\Omega};\mathbb R^l)$.

Let $\Omega$ be an open and bounded subset of $\R^N$. Let $\mathcal{O}(\Omega)$ be the space of all open subsets of $\Omega$. To perform our investigation, it is convenient to introduce the $k$-th $\mathbb C$-elliptic homogeneous linear differential operator $\hat{\B}$ on $\R^N$ from $V_1:=\mathbb R^N \times V$ to $Z$, and observe that when $Z= \R^l \times W$ and $\hat{\B}w := \hat{\B}(w_1, w_2) \coloneqq (\B w_1, {\rm Id} w_2)$, we can recast the model in \eqref{G} into this setting, considering a $\mathfrak g$ therein, coupling the last two variable, i.e.  $\mathfrak g(x,u,v,\B u )=\mathfrak g(x,u, (v,\B u))$.  For $p,r\in[1, \infty]$, let $F_{\hat \B}:  L^r(\Omega; \R^m)\times W^{\hat{\B}, p}(\Omega)\times \mathcal{O}(\Omega) \to [0, \infty)$ be the functional defined by 
         \begin{equation}
             \label{def:funcsuperlinear}
             F_{\hat{\B}}(u,w; D)\coloneqq \int_D f(x, u(x), \hat{\B} w(x))dx,
         \end{equation}
where the density $f:\Omega\times \mathbb R^m \times Z \to [0, \infty)$ is a  Carath\'{e}odory function satisfying one of the following growth conditions:
   \begin{itemize}
       \item[(T1)] if $p, r\in [1, \infty)$, then, there exists a postive constant $C$ such that, for a.e. $x\in\Omega$ and for any $(A, \Xi)\in \R^m \times Z$,
       \begin{equation}
           \notag
           0\leq f(x,A, \Xi)\leq C(1+|A|^r+|\Xi|^p);
       \end{equation}
    \item[(T2)] if $r=\infty, p\in [1, \infty)$, then, there exists $a\in L^\infty(\Omega\times\R^m; [0,\infty))$ such that  for a.e. $x\in\Omega$ and for any $(A,\Xi)\in \R^m \times Z$,
       \begin{equation}
           \notag
           0\leq f(x,A,\Xi)\leq a(x, \xi)(1+ |\Xi|^p);
       \end{equation}
     \item[(T3)] if $r\in [1, \infty), p=\infty$, then, there exists $b\in L^\infty(\Omega\times Z; [0,\infty))$ such that  for a.e. $x\in\Omega$ and for any $(A, \Xi)\in \R^m \times Z$,
       \begin{equation}
           \notag
           0\leq f(x,A,\Xi)\leq b(x, \Xi)(1+ |A|^r);
       \end{equation} 
       \item[(T4)] if $p, r=\infty$, then,
       \begin{equation}
           \notag
           \chi_{\Omega}f\in L^\infty(\R^N\times\R^m \times Z; [0;\infty)).
       \end{equation} 
   \end{itemize}

\begin{Remark}\label{remuvw}
 We emphasize again that, up to the choice of $w:= (u,v)\in $, with $u \in BV^\B(\Omega)$ and $v \in \mathcal M(\overline{\Omega};\mathbb R^l )$ and $\hat B= (\B, Id)$ 
 $Z\ni\hat{\B}w= \hat{\B}(u,v)= (\mathcal B u, {\rm Id} v)\in W \times \R^l$,
the functional $F_{\hat{\B}}$ is of the type of $F$ in \eqref{originalfunctional}, provided that the density $(u,v)$ have the same summability.
\end{Remark}

Before showing our relaxation results, we recall \cite[Theorems 1.1, 3.6]{BFL00} in a convenient way that will fit our framework. Note also that a careful inspection of the proof leads us to conclude that  \cite[Theorems 1.1, 3.6]{BFL00} holds for operators $\A$ of any order and not just for first-order ones as originally stated.

\begin{Theorem}
\label{the:relaxBFL}
     Let $\mathcal A$ be a linear, homogeneous differential operator of order $k$, with $k\in\mathbb N$. Let $\Omega$ be an open subset of $\mathbb R^N$. Let $f:\Omega\times \mathbb R^m \times Z \to [0, \infty)$ be a Carath\'{e}odory function satisfying one of Assumptions ${\rm (T1)}$-${\rm (T4)}$. Then, for all $D\in \mathcal{O}(\Omega), u\in L^r(\Omega; \mathbb R^m), U\in L^p(\Omega; Z)\cap{\rm ker}\A,$ 
            \begin{align}\notag
             \overline{F}_{\A}(u,U; D)\coloneqq\inf\biggl\{ &\liminf_{j\to\infty} \int_D   f(x, u_j(x), U_j(x))dx:  (u_j, U_j)\in L^r(D; \mathbb R^m)\times L^p(D; Z)\notag\\
               &\quad  u_j\to u\;\hbox{in } L^{r}(D; \mathbb R^m),\; U_j\rightharpoonup U\;\hbox{in } L^{p}(D; Z), \; \A U_j=0         \biggr\} 
               \notag\\
               &=\int_D \mathcal{Q}_{\A}f(x, u(x), U(x))dx,\notag
            \end{align}
    where $\mathcal Q_{\A}$ is the $\A$-quasiconvexification of $f(x, A, \cdot)$ defined by \eqref{QAf}, namely
           \begin{equation}
               \notag
               \mathcal{Q}_{A}f(x, A, \Xi)\coloneqq \inf\biggl\{\int_{Y^N} f(x, A, \Xi + \vartheta(y)) dy\; :\; \vartheta\in C^{\infty}_{\rm per}(Y^N; Z)\cap\ker\A, \;\; \int_{Y^N} \vartheta(y)dy=0     \biggr\},
           \end{equation}
    for all $\Xi\in W$.
\end{Theorem}

Deploying this result, we are now able to show an integral representation for the weakly lower semicontinuous envelope of $F_{\hat{\B}}$, i.e., 
\begin{align}
    \overline{F}_{\hat{\B}}(u, w; D)\coloneqq \inf\biggl\{& \liminf_{j\to\infty} F_{\hat{\B}}(u_j, w_j; D)  \; :\; (u_j, w_j)\in L^r(\Omega;\mathbb R^m) \times W^{\B, p}(D)\;\notag\\
    &\quad u_j\to u \hbox{ in }  L^r(\Omega; \R^m), \;\;   w_j\rightharpoonup w \hbox{ in }W^{\B, p}(D)  \biggr\}.\label{def:relenergysuperlinear}
\end{align}

\begin{Theorem}\label{thm:relBSob}
    Let $f:\Omega\times\R^m \times Z   \to [0,\infty)$ be a Carath\'{e}odory function such that one of Assumptions {\rm (T1)}-{\rm (T4)} is satisfied. Let $F_{\hat{\B}}$ and $\overline{F}_{\hat{\B}}$ be the functionals defined by \eqref{def:funcsuperlinear} and \eqref{def:relenergysuperlinear} respectively. Then, 
        \begin{align}\notag
          \overline{F}_{\hat{\B}}(u,w; D) =\int_\Omega Q_{\hat{\B}}f(x, u(x),\hat{\B} w(x))dx,
        \end{align}
    where $\mathcal{Q}_{\hat{\B}}$ is the $\hat{\B}$-quasiconvexification of $f(x, A, \cdot)$, i.e., 
        \begin{equation}
            \notag
            Q_{\hat{\B}}f(x,A, \Xi) = \inf\biggl\{\int_{Y^N} f(x, A, \Xi + \B\widetilde{w}(y))dy \; :\; \widetilde{w}\in C^{\infty}_{c}(Y^N; V_1)\biggr\},
        \end{equation}
    for any $\Xi\in Z$. 
\end{Theorem}

\begin{proof}
    The proof adopts the blow-up technique. As it is custom, first we need to show that $\overline{F}_{\hat{\B}}(u,w; \cdot)$ is the trace of a Radon measure absolutely continuous with respect to $\L^N\res \Omega$. To this end, one can proceed exactly as in \cite[Lemma A.1]{KKZ23}, with minor modifications. Then, one has to show that for $\L^N$-a.e. $x_0\in\Omega$, 
        \begin{equation}
              \notag
              {d\overline{F}_{\hat{\B}}(u,w; \cdot)\over d\L^N}(x_0)= \mathcal Q_{\B}f(x_0, u(x_0), \hat{\B} w(x_0)).
          \end{equation}    
    The proof of the lower bound relies on Proposition \ref{cor1R19} as well as Theorem \ref{the:relaxBFL}. Indeed, we get that
         \begin{align}
              {d\overline{F}_{\B}(u,w; \cdot)\over d\L^N}(x_0)&\geq  {d\overline{F}_{\A}(u,{\hat{\mathcal B} w}; \cdot)\over d\L^N}(x_0) = \mathcal Q_{\A}f(x_0, u(x_0), \hat{\B }w(x_0))= \mathcal Q_{\B}f(x_0, u(x_0), \hat{\B} w(x_0)).\notag
         \end{align}
    For the proof of the upper bound, we proceed exactly as in \cite[Lemma 3.5]{BFL00}, after noticing that we can build up the recovery sequence from the one in \cite[Lemma 3.5]{BFL00}, modifying eq (3.14)  therein.  Indeed, in \cite[Lemma 3.5]{BFL00}, for any fixed $r>0$, the recovery sequence $\tilde w_{n, r}$ is defined by $\tilde w_{j.r}(x)\coloneqq  \tilde w(j(x-x_0)/r)$, where $\tilde w\in C^\infty_{\rm per}(\R^N; V_1)\cap \ker\A$ with $\int_{Y^N} \tilde w(y)dy=0$. Hence, invoking \cite[Lemma 2]{R19}, there exists $w\in C^{\infty}_{\rm per}(Q; V_1)$ such that $\hat \B w=\tilde w$. This allows us to follow the same step as in \cite[Lemma 3.5]{BFL00}, concluding the proof. 
\end{proof}

\begin{Remark}
 With slight modifications, Theorem \ref{thm:relBSob} could be extended to the case $(u,v) \in W^{\B,p}\times L^q$, $\infty \geq p,q>1$, in the spirit of the results in \cite{CRZ1, CRZ2, RZ1}.   
\end{Remark}

\subsection{Interactions with emerging measures in the relaxation process within the Sobolev setting}\label{Sobolevconc}
For a bounded, open subset $\Omega$ of $\R^N$ with Lipschitz boundary, we will discuss the relaxation of $F$ in \eqref{originalfunctional} defined in $W^{\B,p}(\Omega)\times L^1(\overline{\Omega};\mathbb R^l)$, with $p>1$, both for first and $1<k$-th homogeneous linear partial differential $\mathbb C$-elliptic operators $\B$ as in \eqref{Budef}, assuming $ (H_1)$, $ (H_{2}q)$, ($q=1$) and $ (H_{3}p)$ for $p>1$. 

We deal with the following modification of $\mathcal F$ in \eqref{Frelax},
	\begin{align}\nonumber
		\begin{aligned}
		\overline{\mathcal F}[(u,v)]:=
		\inf\left\{\,\liminf_{j\to +\infty}
		F(u_j,v_j)
		\,\left|\,
		\begin{array}{l}
		(u_j,v_j)\in W^{\B,p}(\Omega)\times L^1(\Omega;\R^l),\\
		(u_j, 
  v_j)\overset{\ast}{\rightharpoonup} (u,v)\hbox{ in } W^{\B,p}(\Omega)\times  {\mathcal M}(\overline\Omega;\R^l)
		\end{array}
		\right.\right\}.
		\end{aligned}
		\end{align}

To this end, we start with the case $N\geq 2$, which will provide the Sobolev counterpart of Theorem \ref{thm:relk<N,N>=2}. Since $p>1$, our assumptions on $\B$ and $\Omega$ allow us to invoke the Korn inequality proved in \cite{S24}, which, in turn, ensures that we can essentially replace $W^{\B,p}$ by $W^{k,p}$. For $N \geq 2$, the proof of Theorem \ref{thm:relk<N,N>=2} can be easily modified to include the case $\frac{1}{p}-\frac{K}{N} >0$,  with $u \in W^{\B,p}$ and
using $W^{\B,p} \ni u_j \rightharpoonup u$ weakly in $W^{\B,p}$ (instead of $W^{\B,1} \ni u_j \overset{\ast}{\rightharpoonup} u$ weakly in $BV^\B$). In that setting,
the correct relaxed energy density is simply the restriction of \eqref{repk<N} to $u \in W^{\B,p}$, i.e., the term with $(\mathcal Q_{\B} h)^\infty$ is
removed as $\B^s u = 0$. Apart from easy modifications, natural in $W^{\B,p}$, the proof of the $BV^\B$ setting can essentially be followed
step by step. The only point where one has to be careful is the upper bound.  More precisely, the function $\varphi$ defined in \eqref{def:varphi}, which now has to be $W^{\B,p}$ still allows the cut-off arguments based on it. This is possible only
for $\frac{1}{p}-\frac{K}{N} >0$, and $\frac{1}{p}=\frac{1}{N}$, when $k=1$, in order to avoid the local boundedness which otherwise holds in this case by embedding. The relaxation for  $\frac{1}{p}-\frac{K}{N} <0$ is
a much more classic problem where $u_j \to u$ uniformly and the variable for $u$ in the relaxation formula simply
gets frozen instead of the more complex phenomena observed in \eqref{thm:relk<N,N>=2}.
Namely, a formula of the type in \eqref{rep1} holds,

    \begin{align}\label{repkpN}
                     \overline{\Fcal}(u,v) = \int_{\overline{\Omega}} f_1(u(x))df_2^{\ast\ast}(v)(x)+\int_\Omega {\mathcal Q }_{\B} h(\B u(x))d x, 
                 \end{align}
                 for every $u \in W^{\B, p}(\Omega)$ and $v \in \mathcal M(\overline{\Omega};\mathbb R^l)$.

The case $N=1$ can be treated as in this latter discussed case, since, for any $k \in \mathbb N$, $W^{\B,p}(\alpha, \beta)= W^{k,p}(\alpha,\beta)$, with $p>1$
and the convergence of $u_j \to u$ is uniform, even in the case of first-order operators $\B$. Such an uniform convergence is due to the compact embedding of $W^{1,p}((\alpha,\beta))$ into $C([\alpha, \beta])$. Thus, once again \eqref{repkpN} holds.

\section{Appendix}\label{appendix}

This section is devoted to collect some useful results used in Section \ref{Sect:Applications}. The first result provides a decomposition for functions in $L^1$. For a proof see \cite[Lemma 2.31]{FL}.\\

\begin{Lemma} [Decomposition Lemma in $L^1$]
	\label{decompositionlemma}
	Let $(z_j)_j \subset L^1(\Omega;\mathbb R^d)$ be bounded. Then  
	$z_j$ can be decomposed as	$z_j= z_j^{\rm osc}+ z_j^{\rm conc},$
	with two bounded sequences $(z_j^{\rm osc})_j,(z_j^{\rm conc})_j \subset L^1(\Omega;\mathbb R^d)$ such that 
	$|z_j^{\rm osc}|$ is equiintegrable and $z_j^{\rm conc}$ is the purely concentrated part in the sense that
	$z_j^{\rm conc}\to 0$ in measure. 
	In fact, one may even assume that
		\;\;	$\Lcal^n(\{z_j^{\rm conc}\neq 0\}) \to 0\quad\text{as $j\to +\infty$}. $
	
	Moreover, whenever $f:\R^d\to \R$ is globally Lipschitz, 
	\[
		\|f(z_j)-f(z^{\rm osc}_j)- f(z^{\rm conc}_j)+ f(0)\|_{L^1(\Omega)}\to 0,\quad\text{as $j \to +\infty$}.
	\]
\end{Lemma}

We also use the following lemma which allows us to manipulate boundary values.
\begin{Lemma}\label{boundaryc}
	Let $\Omega \subset \mathbb R^N$ be a  bounded open set, let $A \subset \subset\Omega$ be an open subset with Lipschitz boundary and let $\B$ be a linear partial differential operator of order $k$, with $k\in\mathbb N$. Let $f_1:V \to \R, f_2:\R^l \to \R $ and $h: W\to \R$ be continuous functions satisfying $(H_1), (H_{2}q)$, $(H_3p)$, with $q=p=1$. Consider
	$(u, v) \in BV^\B(\Omega) \times \mathcal M(\Omega;\mathbb R^l)$ such that $|v|(\partial A) = 0$ and assume that 
	$(u_j,v_j)_j \subset W^{k-1,1}(A, \R^m)\times L^1(A;\R^d)$ is a sequence
	satisfying $u_j\overset{*}{\rightharpoonup} u$ in $BV^\B(A)$, 
	$v_j \overset{*}{\rightharpoonup} v$ in $\mathcal M(A;\mathbb R^l)$ and
	$$
	\lim_{k \to +\infty} \int_A(f_1(u_j(x))f_2(v_j(x))+ h(\B u_j(x)))dx=\ell,$$
	for some $\ell < +\infty$. Then, there exists a sequence 
	$(\bar u_j, \bar v_j)_j \subset W^{1,\B}(A)\times L^1(A;\R^l)$
	such that 
	$\bar u_j = u$ on $\partial A$ (in the sense of traces), $\bar u_j \rightharpoonup^* u$ in $BV^\B(A)$, $\bar v_j \rightharpoonup^* v$ in $\mathcal M(A;\mathbb R^l)$, and
	$$
	\limsup_{j \to +\infty}\int_A(f_1 (\bar u_j(x))f_2(\bar v_j(x))+ h(\B\bar u_j(x)))dx \leq \ell.
	$$
\end{Lemma}
The proof is omitted, since it is entirely similar to \cite[Lemma 2.2]{BFMbend} and \cite[Lemmas 4.4 and 4.5]{BZZ}. We point out that the proof relies on the trace theorems \cite[Theorem 1.1]{BDG20} when $k=1$, and it can be deduced with small modifications taking into account the derivatives of order greater than one for the convex combination of functions, exploiting \eqref{BVBequiv}, Theorem \ref{thm2.3} and the trace theorem \cite[Theorems 1.2 and 4.7]{DG24} for operators of order $k>1$.
Eventually, we present the proof of Proposition \ref{proprel**}.

\begin{proof}[Proof of Proposition \ref{proprel**}]
Let $\mathcal F, \overline{\mathcal F}$ and ${\mathcal F}^{\star\star}$ be the functionals defined by \eqref{Flocalized}, \eqref{FlocrelaxB} and \eqref{F**locrelaxB}, respectively. The proof of this result passes through an intermediate inequality. To this end, let 
\begin{align}\nonumber
\begin{aligned}
\overline{\mathcal F}^{\star\star}(u,v;A):=	
\inf\Big\{&\liminf_{j\to +\infty}\int_{\Omega\cap A}\left( f_1(u_j(x))f_2^{\ast\ast}(v_j(x)) + \Q_{\B} h(\B u_j(x))\right)dx: \\
&
W^{\B,1}(\Omega\cap A)\times L^1(\Omega\cap A;\mathbb R^l)\ni (u_j,v_j) \overset{\ast}{\rightharpoonup} (u,v)\hbox{ in } BV^{\B}(\Omega\cap A)\times \mathcal M(A;\mathbb R^l)
\Big\}.
\end{aligned}
	\end{align}
The inequality $\overline{\mathcal{F}}\geq \overline{\mathcal F}^{\star \star}$ is trivial since $f_2^{\ast \ast} \leq f_2$ and $h \geq Q_{\mathcal B}h$, and consequently 
 $\overline{\mathcal{F}}\geq {\mathcal F}^{\ast \ast}$, since in the definition of ${\mathcal F}^{\ast \ast}$ there are more test sequences than in the definition of $\overline{\mathcal F}^{\star\star}$.

 Now, we prove that $\overline{\mathcal{F}}\leq {\mathcal F}^{\star \star}$.
 By the definition \eqref{F**locrelaxB} there exists $(u_j, v_j)\in BV^{\B}(\Omega\cap A)\times L^1(\Omega\cap A;\R^l) \overset{\ast}{\rightharpoonup} (u,v)\in BV^{\B}(\Omega\cap A)\times \mathcal M(A;\mathbb R^l)$ such that 
      \begin{align}
          \lim_{j\to+\infty} \int_{\Omega\cap A} f_1(u_j(x))f_2^{\ast \ast}(v_j(x))dx + \int_{\Omega\cap A}d\Q_{\B} h(\B u_j(x))dx= {\mathcal F}^{\ast \ast}(u,v; A) \label{f4}
      \end{align}
For every $j\in\Nb$, take $u_j \in BV^{\B}(\Omega\cap A)$ as above. By Theorem \ref{thm:RelaxationB}, applied to $\Omega \cap A$, in view of \eqref{lippartbound}, there exists $(w_{j,i})_i\subset W^{\B,1}(\Omega \cap A)$ such that $w_{j,i}$ weakly* converges to $u_j$ in $BV^{\B}(\Omega \cap A;\mathbb R^d)$.
Observe also that the convergence of $w_{j,i}$ to $u_j$ is also strictly in area and
     \begin{align}
         \notag
         \lim_{i\to\infty}\int_{\Omega\cap A} h(\B w_{k,i}(x))dx = \int_{\Omega\cap A} d\Q_{\B} h(\B u_k)(x)dx. 
     \end{align}
\noindent Observe that by \eqref{exiQBfinfty}, $(\Q_{\B}h)^\sharp$ in Theorem \ref{thm:RelaxationB} coincides with $(\Q_{\B}h)^\infty$. In view of growth condition $(H_3p)$, it follows that 
        $ \limsup_{i\to\infty} \|\B w_{j,i} \|_{L^1(\Omega \cap A)}\leq C(1+ |\B u_j|(\Omega \cap A)). $
Furthermore, by dominated convergence, for any $v_j\in L^1(A)$, we get
    \begin{align}
        \notag
        \lim_{i\to\infty}\int_{\Omega\cap A} f_1(w_{j,i}(x))f_2^{\ast \ast}(v_j(x))dx = \int_{\Omega\cap A} f_1(u_{j}(x))f_2^{\ast \ast}(v_j(x))dx.
    \end{align}
Thus, by \eqref{f4},
\begin{align*}
\lim_{j\to +\infty} \lim_{i\to\infty}\left(\int_{\Omega\cap A} f_1(w_{j,i}(x))f_2^{\ast \ast}(v_j(x))dx +\int_{\Omega\cap A} h(\B w_{j,i}(x))dx \right) =\\
\lim_{j\to +\infty} \left(\int_{\Omega\cap A} f_1(u_{j}(x))f_2^{\ast \ast}(v_j(x))dx+ \int_{\Omega\cap A} d\Q_{\B} h(\B u_j)(x)dx\right)= {\mathcal F}^{**}(u,v; A).
\end{align*}
Therefore, we make use of a diagonalization argument to find a sequence $\widetilde{u}_j:= w_{j, i(j)}$ such that 
     \begin{align}
         &\lim_{j\to\infty} \int_{\Omega\cap A} f_1(\widetilde{u}_j(x))f_2^{\ast \ast} (v_j(x))dx +\int_{\Omega\cap A} h(\B \widetilde{u}_j(x))dx\notag\\
         &= \lim_{j\to\infty} \int_{\Omega\cap A} f_1(u_j(x))f_2^{\ast \ast} (v_j(x))dx +\int_{\Omega\cap A} d\Q_\B h(\B u_j)(x)dx 
         = {\mathcal F}^{**}(u,v; A).\label{f5} 
     \end{align}
In addition, $\limsup_{j\to \infty}\|\B \widetilde{u}_j\|_{L^1(\Omega\cap A;\mathbb R^d)}\leq C(1+\limsup_{j\to\infty} |\B u_j|(\Omega\cap A))$ and $\widetilde{u}_j\to u$ in $L^1(\Omega\cap A;V)$. In particular, $\widetilde{u}_j\overset{\ast}{\rightharpoonup} u$ in $BV^{\B}(\Omega\cap A)$. 

Now, we fix $\widetilde{u}_j\in W^{\B, 1}(\Omega\cap A)$ and $v_j\in L^1(\Omega\cap A;\mathbb R^l)$ as in the left hand side of \eqref{f5}. Our aim is to find a suitable recovery sequence for the $v_j$ which allows us to switch from $f_2^{\ast\ast}$ to $f_2$. It is well-known (see \cite[Theorem 6.68 and Remark 6.69]{FL}) that 
     \begin{align}
         \inf\biggl\{ &\liminf_{n\to \infty} \int_{\Omega\cap A} f_{1}(\widetilde{u}_j(x))f_2(v_{j,n}(x))dx +\int_{\Omega\cap A} h(\B \widetilde{u}_j(x)) dx : (v_{j,n})_n\in L^1(\Omega\cap A:\mathbb R^l), v_{j,n}\rightharpoonup v_j \hbox{ in }L^1(\Omega\cap A;\mathbb R^l)\biggr\}\notag\\
         &\quad  = \int_{\Omega\cap A} f_{1}(\widetilde{u}_j(x))f_2^{\ast\ast}(v_{j}(x))dx +\int_{\Omega\cap A} h(\B \widetilde{u}_j(x)) dx. \label{f1}
     \end{align}
In particular,  up to a suitable extension of $v_j$ (still denoted with the same symbol) as a function of $L^1(A;\mathbb R^l)$, we have
      \begin{align}
         \inf\biggl\{ &\liminf_{n\to \infty} \int_{\Omega\cap A} f_{1}(\widetilde{u}_j(x))f_2(v_{j,n}(x))dx +\int_{\Omega\cap A} h(\B \widetilde{u}_j(x)) dx : (v_{j,n})_n\in L^1(\Omega\cap A;\mathbb R^l), v_{j,n}\overset{\ast}{\rightharpoonup} v_j \hbox{ in }\mathcal M(A;\mathbb R^l)\biggr\}\notag\\
         &\quad  \leq \int_{\Omega\cap A} f_{1}(\widetilde{u}_j(x))f_2^{\ast\ast}(v_{j}(x))dx +\int_{\Omega\cap A} h(\B \widetilde{u}_j(x)) dx. \notag
     \end{align}
From \eqref{f1}, we deduce that there exists $(v_{j, n})_n\in L^1(\Omega\cap A;\mathbb R^l)$ with $v_{j,n}\rightharpoonup v_j$ as $n\to \infty$ such that   
     \begin{align*}
         &\lim_{n\to\infty} \int_{\Omega\cap A} f_{1}(\widetilde{u}_j(x))f_2(v_{j,n}(x))dx +\int_{\Omega\cap A} h(\B \widetilde{u}_j(x)) dx =  \int_{\Omega\cap A} f_{1}(\widetilde{u}_j(x))f_2^{\ast \ast}(v_{j}(x))dx +\int_{\Omega\cap A} h(\B \widetilde{u}_j(x)) dx,
     \end{align*}
     for every $j \in \mathbb N$.
     Observe, also  that, by extending, for every $n$, $v_{j,n}$ as $v_j$ in $A \setminus \Omega$, it results that $v_{j,n}\weaklystar v_j$ in $\mathcal M( A;\mathbb R^l)$.
Taking the limit as $j\to \infty$, from \eqref{f5}, we get 
    \begin{align}
         &\lim_{j\to\infty}\lim_{n\to\infty} \left(\int_{\Omega\cap A} f_{1}(\widetilde{u}_j(x))f_2(v_{j,n}(x))dx +\int_{\Omega\cap A} h(\B \widetilde{u}_j(x)) dx \right)\notag\\
         &\quad \quad \leq (=)\lim_{j\to\infty}\left( \int_{\Omega\cap A} f_{1}(\widetilde{u}_j(x))f_2^{\ast \ast}(v_{j}(x))dx +\int_{\Omega\cap A} h(\B \widetilde{u}_j(x)) dx\right)=
         {\mathcal F}^{**}(u,v; A). \label{f2}
     \end{align}
The metrizability of weak* convergence in $\mathcal M(A;\mathbb R^l)$ on bounded sets, allows us to use a diagonalization argument which, in turn, leads us to extract a subsequence $\widetilde{v}_j:= v_{j, n(j)}$ such that $\widetilde{v}_j\overset{\ast}{\rightharpoonup} v$ in $\mathcal M(A;\mathbb R^l)$ and
      \begin{align*}
          \lim_{j\to\infty}\lim_{n\to\infty} \left(\int_{\Omega\cap A} f_{1}(\widetilde{u}_j(x))f_2(v_{j,n}(x))dx +\int_{\Omega\cap A} h(\B \widetilde{u}_j(x)) dx \right)
          = \lim_{j\to\infty} \left(\int_{\Omega\cap A} f_{1}(\widetilde{u}_j(x))f_2(\widetilde{v}_j(x))dx +\int_{\Omega\cap A} h(\B \widetilde{u}_j(x)) dx.\right)
      \end{align*}      
This combined with \eqref{f2} yields to
    \begin{align}
        &\lim_{j\to\infty} \left(\int_{\Omega\cap A} f_{1}(\widetilde{u}_j(x))f_2(\widetilde{v}_j(x))dx +\int_{\Omega\cap A} h(\B \widetilde{u}_j(x)) dx\right)\notag\\
        &\quad \quad\leq (=) \lim_{j\to\infty}\left( \int_{\Omega\cap A} f_{1}(\widetilde{u}_j(x))f_2^{\ast \ast}(v_j(x))dx +\int_{\Omega\cap A} h(\B \widetilde{u}_j(x)) dx \right)
        ={\mathcal F}^{**}(u,v; A). \label{f3}
    \end{align}
Recall that $W^{\B,1}(\Omega\cap A) \ni \tilde u_j \weaklystar u$ in $BV^\B(\Omega\cap A)$ and that (with an abuse of notation which allows us to identify any $\mathbb R^l$-valued $L^1$ function in $\Omega \cap A,$ with its extension as a measure in $\mathcal M(A;\mathbb R^l)$), $L^1(\Omega\cap A;\mathbb R^l)\ni \tilde v_j \weaklystar v$ in $\mathcal M(A;\mathbb R^l)$. 
Then, from \eqref{f3}, definitions \eqref{Flocalized} and \eqref{FlocrelaxB}, it follows that 
     \begin{align}
     \overline{\mathcal F}(u,v; A)] 
         \leq \lim_{j\to\infty} \int_{\Omega\cap A} f_{1}(\widetilde{u}_j(x))f_2(\widetilde{v}_j(x))dx +\int_{\Omega\cap A} h(\B \widetilde{u}_j(x)) dx
          \leq(=) {\mathcal F}_{\ast \ast}(u,v; A),\notag
     \end{align}
    which concludes the proof. 
\end{proof}

\subsection*{Acknowledgment}
LD acknowledges support of the Austrian Science Fund (FWF) projects 10.55776/ESP1887024 and 10.55776/Y1292.
EZ is a member of INdAM GNAMPA, whose support is gratefully acknowledged.
Her research has also been supported by PRIN 2022: Mathematical Modelling of Heterogeneous Systems (MMHS)
- Next Generation EU CUP B53D23009360006.
She is also indebted with professor G. Dal Maso for helpful discussions about the proof of Theorem \ref{thm:asBHHZsci}.



\begin{thebibliography}{Biblio}
\bibitem{ABF} E. Acerbi, G. Bouchitt\'{e} and I. Fonseca, Relaxation of convex functionals: the gap problem, Ann. Inst. H. Poincar\'e, 20 (2003),
pp. 359-390.

\bibitem{AB} J. J. Alibert and G. Bouchitt\'{e}, Non-Uniform Integrability
and Generalized Young Measures, J. Convex Anal., 4 (1997), pp. 129-147. 

\bibitem{ADC} M. Amar and V. De Cicco, 
Relaxation of quasi-convex integrals of arbitrary order, Proc. Royal Soc. Edinburg., 12A (1994), pp. 927-946.

\bibitem{AD}
L. Ambrosio and G. Dal Maso, On the Relaxation in $BV(\Omega;\R^m)$ of quasi-convex integrals. J. Funct. Anal., 109, (1992),  76--97. 
    
\bibitem{AFP00} L. Ambrosio, N. Fusco and D. Pallara, Functions of bounded variation and free discontinuity problems, Oxford Mathematical Monographs, New York, 2000.

\bibitem{AR20} A. Arroyo-Rabasa, Slicing and fine properties for functions with bounded $\mathcal{A}$-variation,  Preprint, 2020.

\bibitem{AR21} A. Arroyo-Rabasa,  Characterization of Generalized Young Measures Generated by $\A$-free Measures, Arch. Rational Mech. Anal.,  242 (2021), pp. 235–325.

\bibitem{APR20} A. Arroyo-Rabasa, G. De Philippis and F. Rindler,  Lower semicontinuity and relaxation of linear-growth integral functionals under PDE constraints, Adv. Calc. Var.,  13 (2020), pp. 219-255.

\bibitem{ARS25} A. Arroyo-Rabasa and A. Skorobogatova,  A look into some of the fine properties of functions with bounded $\mathcal{A}$-variation,  ESAIM: COCV,  31 (2025), n. 50,


\bibitem{BZZ} J.-F. Babadjian, E. Zappale and H. Zorgati,
Dimensional reduction for energies with linear growth involving the bending moment,
J. Math. Pures Appl.,  90 (2008), pp. 520-549.

\bibitem{BFT} A.C. Barroso, I. Fonseca and R. Toader,
A relaxation theorem in the space of functions of bounded
              deformation, Ann. Scuola Norm. Sup. Pisa Cl. Sci. (4), {29}, (2000), n. 1. pp. 19-49.
\bibitem{BMZNoDEA}  A. C. Barroso, J. Matias and E. Zappale,
A global method for relaxation for multi-levelled structured deformations,
Nonlinear Diff. Equ. Appl.,  31 (2024), n. 50.

\bibitem{BMZGaeta}
 A. C. Barroso, J. Matias and E. Zappale, Multi-level structured deformations: relaxation via a global method approach, to appear in  {\it Comm. Math. Anal. Appl}, special issue. 
 
\bibitem{BHHZ} G. Bertazzoni, P. Harjulehto, P. H\"asto and E. Zappale, in preparation.

\bibitem{BFMglob} G. Bouchitt\'e, I. Fonseca and M. L. Mascarenhas, A global method for relaxation, Arch. Rational Mech. Anal., 145 (1998), pp. 51-98.


\bibitem{BFMbend} 
	G. Bouchitt\'{e}, I. Fonseca and M. L. Mascarenhas,
    Bending moment in membrane theory,  J. Elasticity,
	 73 (2004), pp. 75-99.
	

\bibitem{BFL00} A. Braides, I. Fonseca and G. Leoni,
$\A$-quasiconvexity: relaxation and homogenization,
ESAIM: COCV,  5 (2000), pp. 539-577.

\bibitem{BKP10} K. Bredies, K. Kunisch, and T. Pock: Total generalized variation, SIAM J. Imaging Sci.  3. (3), (2010), pp. 492--526.
\bibitem{BDG20} D. Breit, L. Diening and F. Gmeineder,  On the trace operator for functions of bounded $\A$-variation, Anal. PDE,  13 (2020), pp. 559–594.

\bibitem{CL}  J.W. Cahn and  F. L\"{a}rch\'{e}, Surface stress and the chemical equilibrium of small crystals—II. Solid particles embedded in a solid matrix. Acta Metall. 30 (1981), pp. 51--56. 
	
\bibitem{CRZ1} G. Carita, A. M. Ribeiro and E. Zappale, Relaxation for some integral functionals in $W^{1,p}_w\times L^q_w$, Bol. Soc.
Port. Mat., (2010), pp. 47–53.

\bibitem{CRZ2}  G. Carita, A. M. Ribeiro and E. Zappale, An homogenization result in {$W^{1,p}\times L^q$},
J. Convex Anal., 18 (2011), pp. 1093-1126.
		
\bibitem{CZA}
	G. Carita and  E. Zappale,
	A relaxation result in {$BV\times L^p$} for integral
		functionals depending on chemical composition and elastic
		strain,  Asymptot. Anal., 100 (2016), pp. 1-20.

\bibitem{CZE}
G. Carita and E. Zappale,
Integral representation results in $BV\times
L^p$,  ESAIM Control Optim. Calc. Var.,  23 (2017),  pp. 1555-1599.

\bibitem{CFVG20}
M. Caroccia, M. Focardi and N. Van Goethem,
On the integral representation of variational functionals on $BD$, SIAM J. Math. Anal.,  52 (2020), pp. 4022-4067.

\bibitem{CLR} Y. Chen, S. Levine and M. Rao, Variable exponent, linear growth functionals in image restoration, SIAM J. Appl. Math. 66 (4), (2006) pp. 1383-1406.

\bibitem{CG22} S. Conti and F. Gmeineder,
$\mathcal A$-quasiconvexity and partial regularity. 
 Calc. Var. 61 (2022), n. 215.

\bibitem{DM}  G. Dal Maso, An introduction to  $\Gamma$ -convergence,
Progr. Nonlinear Differential Equations Appl., {\bf 8}
Birkhäuser Boston, Inc., Boston, MA, 1993.

\bibitem{DD} A. De Simone and G. Dolzmann, Existence of minimizers for a variational problem in two-dimensional nonlinear magnetoelasticity, Arch. Ration. Mech. Anal., 144 (1998), pp. 107-120. 

\bibitem{Da08} B. Dacorogna, Direct Methods in the Calculus of Variations, Springer, Berlin, 2008.
		
\bibitem{D} F. Demengel,
Fonctions à hessien borné
Ann, Inst. Fourier, 34 (1984), pp. 155-190.

\bibitem{DT} F. Demengel and R. Temam,
Convex Functions of a Measure and Applications,
 Ind. Univ. Math. J., 33 (1984), pp. 673-709.

\bibitem{DG20} L. Diening and F. Gmeineder,  Continuity Points Via Riesz Potentials for $\mathbb C$-Elliptic Operators, Q. J. Math., 71 (2020), pp. 1201–1218.


\bibitem{DG24} L. Diening and F. Gmeineder,   Sharp Trace and Korn Inequalities for Differential Operators, Potential Anal., (2024).

\bibitem{HHbook} L. Diening, P. Harjulehto, P. H\"ast\"o, M. Rů\v zi\v cka, Lebesgue and Sobolev spaces with variable exponents, in: Lecture Notes in Mathematics, vol. 2017, Springer, Heidelberg, (2011).


\bibitem{FMZ25} R. Ferreira, J. Matias, and E. Zappale, 
    Junction in a thin multi-domain for nonsimple grade two
              materials in BH, Nonlinear Anal. Real World Appl., 84 (2025), n. 104322.

\bibitem{FHP} I. Fonseca, A. Hagerty and R. Paroni,
Second-order structured deformations in the space of functions of bounded Hessian, J. Nonlinear Sci.,
 29, (2019), pp. 2699-2734.

 	\bibitem{FKP1}
			I. Fonseca, D. Kinderlehrer and P. Pedregal, Relaxation in $BV\times L^\infty$ of functionals depending on strain and composition. In: Boundary value problems for partial differential equations and applications.  RMA Res. Notes Appl. Math.   29,Masson, Paris, (1993) pp. 113-152.
			
			\bibitem{FKP2}
			I. Fonseca, D. Kinderlehrer and P. Pedregal, Energy functionals depending on elastic strain and chemical composition., Calc. Var. PDE,  2 (1994),    pp.  283-313.

\bibitem{FL} I. Fonseca and G. Leoni, Modern Methods in the Calculus of Variations: $L^{p}$ spaces, Springer, New York, 2007.

\bibitem{FM92} I. Fonseca and S. M\"{u}ller,   Quasiconvex integrands and Lower Semicontinuity in $L^1$, SIAM J. Math Anal.,  23 (1992), pp. 1081–1098.


\bibitem{FM99} I. Fonseca and S. M\"{u}ller,   $\A$-quasiconvexity, lower semicontinuity, and Young measures, SIAM J. Math Anal.,  30 (1999), pp. 1355–1390.
	
    \bibitem{FMARMA} I. Fonseca and S. M\"{u}ller, Relaxation of
		quasiconvex functionals in $BV(\Omega;\mathbb{R}^{p})$ for integrands $f(x,u,\nabla u)$,Arch. Ration. Mech. Anal., 123
	(1993), pp. 1-49.
	
	\bibitem{FMP} 
I. Fonseca, S. M\"{u}ller and P. Pedregal, Analysis of concentration and oscillation effects generated by gradients, SIAM J. Math. Anal., 29 (1998), pp. 736-756.


\bibitem{FMZ} M. Friedrich, J. Matias and E. Zappale, in preparation.

\bibitem{GR19} F. Gmeineder and B. Rai\textcommabelow{t}\u{a},  Embeddings for $\A$-weakly differentiable functions on domains,  J. Funct. Anal., 227 (2019), n. 108278.

\bibitem{H} A. Hagerty, Relaxation of functionals in the space of vector-valued
              functions of bounded {H}essian, Calc. Var. Partial Differential Equations, 58, (2019), n. 1, Paper No. 4, 38.

	\bibitem{HKW}
	D. Henrion, M. Kru\v{z}\'{i}k and T. Weisser, Optimal control problems with oscillations, concentrations and discontinuities, Automatica, 103 (2019), pp. 159-165.
	
	\bibitem{KKK}
	 A. Ka{\l} amajska,  S.~Kr\"{o}mer and M. Kru\v{z}\'{\i}k,   Weak lower semicontinuity by means of anisotropic parametrized measures.  Trends in Applications of Mathematics to Mechanics (eds.: E.~Rocca, U.~Stefanelli, L.~Truskinovsky, and A.~Visintin), Springer INdAM Series 27, Springer Cham, Switzerland (2018), pp. 23-52. 
		
		\bibitem{KK} A. Ka{\l} amajska and M. Kru\v{z}\'{\i}k,  Oscillation and concentrations in sequences of gradients,  ESAIM: COCV,  14 (2008), pp. 71–104.

\bibitem{KK16} B. Kirchheim and J. Kristensen, On rank one convex functions that are homogeneous
of degree one, Arch. Ration. Mech. Anal.,  221 (2016), pp. 527–558.

\bibitem{KR} C. Kreisbeck and F. Rindler,
Thin-film limits of functionals on {$\mathcal A$}-free vector
              fields, Indiana Univ. Math. J., 64 (2015), pp. 1383-1423.


\bibitem{KR10} J. Kristensen and F. Rindler,
Relaxation of signed integral functionals in $BV$, Calc. Var., 37 (2010), pp. 29-62.


\bibitem{KKZ23} S. Kr\"{o}mer, M. Kru\v{z}\'{i}k and E. Zappale, Relaxation of functionals with linear growth: Interactions of emerging measures and free discontinuities, Adv. Calc. Var.,  16 (2023), pp. 835-865.


\bibitem{LDR} H. Le Dret and A. Raoult, 
     Variational convergence for nonlinear shell models with
              directors and related semicontinuity and relaxation results,
   Arch. Ration. Mech. Anal., 154 (2000), pp. 101-134.
 
\bibitem{M1} F. Murat,
 Compacit\'{e} par compensation, Bull. Soc. Math. France M\'em. 60 (1979), pp. 125-127.
     

\bibitem{M2} F. Murat,
Compacit\'{e} par compensation. {II}, Proc. Intern. Meeting on Recent Methods in Nonlinear Analysis (Rome, 1978), Pitagora, Bologna, 1979, pp. 245-256.

 
 \bibitem{M3} F. Murat,
 Compacit\'{e} par compensation: condition n\'ecessaire et suffisante de continuit\'e{} faible sous une hypoth\`ese de rang constant,
Ann. Scuola Norm. Sup. Pisa Cl. Sci. (4),
 8 (1981), pp. 69-102.
 
\bibitem{MT}
F. Murat and L. Tartar, Optimality conditions and homogenization,
 Nonlinear variational problems (Isola d'Elba, 1983),Res. Notes in Math. 127, Pitman, Boston, MA, 1985, pp. 1-8.

\bibitem{PPRV} V. Pagliari, K. Papafitsoros, B. Rai\textcommabelow{t}\u{a} and A. Vikelis, Bilevel Training Schemes in Imaging for Total Variation–Type Functionals with Convex Integrands, SIIMS, 15 (2022), pp. 1690–1728.

\bibitem{R19} B. Rai\textcommabelow{t}\u{a},  Potentials for $\A$-quasiconvexity.  Calc. Var.,  58 (2019), n. 109.

\bibitem{R19_1} B. Rai\textcommabelow{t}\u{a}, Critical $L^p$-differentiability of $BV^\A$-maps and cancelling operators, Trans. Amer. Math. Soc., 372 (2024), pp. 7297-7326.

\bibitem{R24} B. Rai\textcommabelow{t}\u{a},
A simple construction of potential operators for compensated compactness, Quart. J. Math.,
75 (2024), pp. 451 – 456.


	
\bibitem{RS20} B. Rai\textcommabelow{t}\u{a} and A. Skorobogatova, Continuity and canceling operators of order $n$ on $\mathbb R^n$,  Calc. Var.,  59 (2020), n. 85.
\bibitem{RZ1} A. M. Ribeiro,  E. Zappale,
		Relaxation of certain integral functionals depending on strain and chemical composition.  Chin. Ann. Math. Ser. B, 34 (2013), 491--514. 
	
	\bibitem{RZ} A. M. Ribeiro  E. Zappale, Lower
		semicontinuous envelopes in $W^{1,1}\times L^{p},$ Banach Center
	Publications,  101 (2014), 187--206. Erratum: Lower
		semicontinuous envelopes in $W^{1,1}\times L^{p}$, 101 (2014), online.

	\bibitem{RS} F. Rindler, G. Shaw,  Liftings, Young measures, and lower semicontinuity.   Arch. Ration. Mech. Anal. 232, (2019),  1227-1328.

\bibitem{S24}  S. Schiffer,
On the complex constant rank condition and inequalities for differential operators,
Nonlinear Anal.,  239
(2024), n. 113435.

\bibitem{T2} L. Tartar,
Compensated compactness and applications to partial
              differential equations, Nonlinear analysis and mechanics: Heriot-Watt Symposium, IV, Res. Notes in Math. 39, Pitman, Boston, Mass.-London, 1979, pp. 136-212.

\bibitem{T1}
L. Tartar, Compacit\'{e} par compensation: r\'{e}sultats et perspectives,
Nonlinear partial differential equations and their
              applications. Coll\`ege de France Seminar, IV
              (Paris, 1981/1982), Res. Notes in Math. 84, Pitman, Boston, MA, 1983, pp. 350-369.
    
\bibitem{T3}
L. Tartar,
The Compensated Compactness Method Applied to Systems of Conservation Laws, Systems of Nonlinear Partial Differential Equations, J. M. Ball, eds. NATO Science Series C, Springer; Dordrecht, 1983, pp. 263-285. 

\bibitem{T4} L. Tartar,
Some remarks on separately convex functions,
Microstructure and phase transition, IMA Vol. Math. Appl., 54 (1993), pp. 191-204.

\bibitem{T83} R. Temam,  Probl\`{e}mes math\'{e}matiques en plasticit\'{e}, Gauthier-Villars, Paris, 1983.


\bibitem{V13} J. Van Schaftingen,   Limiting Sobolev inequalities for vector fields and canceling linear
differential operators, J. Eur. Math. Soc.,  15 (2013), pp. 877–921.
   
\end{thebibliography}
\end{document}